\pdfoutput=1
\documentclass[12pt, reqno]{amsart}

\usepackage{mathtools}
\mathtoolsset{showonlyrefs}
\numberwithin{equation}{section}

\usepackage{graphicx} 

\usepackage{mathtools}
\mathtoolsset{showonlyrefs}

\usepackage{amssymb}
\usepackage{amsthm}
\usepackage{amscd}
\usepackage{amsmath}
\usepackage{epic,eepic}
\usepackage{mathrsfs}
\usepackage{latexsym}
\usepackage{pifont}
\usepackage[all]{xy}
\usepackage {cancel}

\usepackage{tikz}

\usetikzlibrary{arrows.meta}

\usepackage{hyperref}

\theoremstyle{plain}
\newtheorem{lemma}{Lemma}[section]
\newtheorem{prop}[lemma]{Proposition}
\newtheorem{thm}[lemma]{Theorem}
\newtheorem{cor}[lemma]{Corollary}
\newtheorem{intthm}{Theorem}

\theoremstyle{definition}

\newtheorem{rem}[lemma]{Remark}
\newtheorem{defi}[lemma]{Definition}
\newtheorem{exa}[lemma]{Example}

\newcommand{\bde}{\begin{defi}}
\newcommand{\ede}{\end{defi}\vspace{1mm}}
\newcommand{\ble}{\begin{lemma}}
\newcommand{\ele}{\end{lemma}}
\newcommand{\bpr}{\begin{prop}}
\newcommand{\epr}{\end{prop}}
\newcommand{\bt}{\begin{thm}}
\newcommand{\et}{\end{thm}}
\newcommand{\bco}{\begin{cor}}
\newcommand{\eco}{\end{cor}}
\newcommand{\bre}{\begin{rem}}
\newcommand{\ere}{\end{rem}}
\newcommand{\bex}{\begin{exa}}
\newcommand{\eex}{\end{exa}}
\newcommand{\bpf}{\begin{proof}}
\newcommand{\epf}{\end{proof}}

\newcommand{\mcA}{\mathcal{A}}
\newcommand{\mcB}{\mathcal{B}}
\newcommand{\mcC}{\mathcal{C}}
\newcommand{\mcD}{\mathcal{D}}
\newcommand{\mcE}{\mathcal{E}}
\newcommand{\mcF}{\mathcal{F}}
\newcommand{\mcG}{\mathcal{G}}

\newcommand{\mcL}{\mathcal{L}}
\newcommand{\mcM}{\mathcal{M}}
\newcommand{\mcN}{\mathcal{N}}
\newcommand{\mcO}{\mathcal{O}}

\newcommand{\mcQ}{\mathcal{Q}}

\newcommand{\mcS}{\mathcal{S}}
\newcommand{\mcT}{\mathcal{T}}

\newcommand{\mcV}{\mathcal{V}}

\newcommand{\mbC}{\mathbb{C}}

\newcommand{\mbF}{\mathbb{F}}
\newcommand{\mbG}{\mathbb{G}}

\newcommand{\mbP}{\mathbb{P}}

\newcommand{\mbR}{\mathbb{R}}

\newcommand{\mbZ}{\mathbb{Z}}

\newcommand{\mfS}{\mathfrak{S}}

\newcommand{\mfb}{\mathfrak{b}}
\newcommand{\mfc}{\mathfrak{c}}

\newcommand{\mfg}{\mathfrak{g}}
\newcommand{\mfh}{\mathfrak{h}}

\newcommand{\mfl}{\mathfrak{l}}

\newcommand{\mfo}{\mathfrak{o}}

\newcommand{\mfs}{\mathfrak{s}}
\newcommand{\mft}{\mathfrak{t}}

\newcommand{\msD}{\mathscr{D}}
\newcommand{\msE}{\mathscr{E}}
\newcommand{\msF}{\mathscr{F}}
\newcommand{\msG}{\mathscr{G}}

\newcommand{\msL}{\mathscr{L}}

\newcommand{\msO}{\mathscr{O}}
\newcommand{\msP}{\mathscr{P}}

\newcommand{\msX}{\mathscr{X}}

\newcommand{\SSP}{\vspace{3mm}}
\newcommand{\LSP}{\vspace{5mm}}
\newcommand{\mr}{\mathrm}
\newcommand{\N}{N}
\newcommand{\LL}{L}
\newcommand{\MM}{M}
\newcommand{\Diag}{\rotatebox[origin=c]{45}{$\Leftarrow$}}

\newcommand{\A}{\mcA}

\title[The moduli space of dormant opers on elliptic curves]{The moduli space of dormant opers on elliptic curves}
\author{Naoka Karube}
\author{Yasuhiro Wakabayashi}

\address{\emph{Naoka Karube}
 \newline
 \textnormal{Graduate School of Information Science and Technology, Osaka University, Suita, Osaka 565-0871, JAPAN.}
 \newline
 \textnormal{\texttt{karube.naoka@ist.osaka-u.ac.jp}}}

 \address{\emph{Yasuhiro Wakabayashi}
 \newline
 \textnormal{Graduate School of Information Science and Technology, Osaka University, Suita, Osaka 565-0871, JAPAN.}
 \newline
  \textnormal{\texttt{wakabayashi@ist.osaka-u.ac.jp}}}

\date{}
\begin{document}
\maketitle

\begin{abstract}
A dormant oper is a specific type of principal bundle with a flat connection, defined on an algebraic curve in positive characteristic. The moduli spaces of dormant opers and  their variants,  known as dormant Miura opers, have been studied in various contexts. This paper focuses on the case where the underlying spaces are (possibly nodal) elliptic curves and provides a detailed examination of the geometric structures of their moduli spaces. In particular, we explicitly describe dormant (generic Miura) opers in terms of regular elements in an associated Lie algebra and establish the connectedness of these moduli spaces. We also explore generalizations to higher level and prime-power characteristic.
\end{abstract}
\tableofcontents 

\section{Introduction} \label{S1}

\LSP
\subsection{Opers on algebraic curves} \label{SSa}

 {\it $G$-opers} for a reductive group $G$ are defined as  $G$-bundles with flat connections on  algebraic curves,  subject to certain additional conditions.
 They were introduced by A. Beilinson and V. Drinfeld as fundamental objects in the geometric Langlands correspondence for constructing  Hecke eigensheaves on the moduli space of $G$-bundles by associating them with the base space of Hitchin's integrable system (cf. \cite{BeDr1}).
  $G$-opers and their variants, known as {\it (generic) Miura $G$-opers},
 play central roles in various research areas, including Teichm\"{u}ller theory, 
 and continue to be actively studied,  e.g., in relation to Gaiotto's conjecture
 (cf. ~\cite{DFKMMN}, ~\cite{Gai}).

For example, a $\mathrm{GL}_n$-oper on a smooth projective curve $X$ 
 consists of the following data:
\begin{align} \label{Eq21ke}
\msF := (\mcF, \nabla, \{ \mcF^j \}_{j=0}^n),
\end{align}
where $\mcF$ denotes  a rank $n$ vector bundle on $X$, $\nabla$ denotes  a flat connection $\mcF \rightarrow \Omega_X \otimes \mcF$ on $\mcF$ (i.e., $\nabla$ is a $k$-linear morphism satisfying the usual Leibniz rule), and $\{ \mcF^j \}_j$ denotes a decreasing filtration $0 = \mcF^n \subseteq \mcF^{n-1} \subseteq \cdots \subseteq \mcF^0 = \mcF$ on $\mcF$ such that each subquotient $\mcF^j/\mcF^{j+1}$ is a line bundle;
this collection   must satisfy the following conditions:
\begin{itemize}
\item
$\nabla (\mcF^j) \subseteq \Omega_X \otimes \mcF^{j-1}$ for every $j = 1, \cdots, n-1$;
\item
The $\mcO_X$-linear morphism $\mcF^j/\mcF^{j+1} \rightarrow \Omega_X \otimes (\mcF^{j-1}/\mcF^{j})$ induced by $\nabla$ is an isomorphism for every $j = 1, \cdots, n-1$.
\end{itemize}
After applying  a gauge transformation to the underlying vector bundle, each   $\mr{GL}_n$-oper can be transposed into a linear differential operator $D$ locally given by
\begin{align} \label{Eq12345}
D = \frac{d^n}{dx^n} + a_1 \frac{d^{n-1}}{dx^{n-1}} + \cdots + a_{n-1} \frac{d}{dx} +  a_n. 
\end{align}
Here, $x$ is 
 a local coordinate in $X$ and $a_1, \cdots, a_n$ are variable coefficients.

When $X$ is defined over 
the field of complex numbers $\mbC$, 
 a $\mathrm{PGL}_2$-oper on $X$ is equivalent to a {\it projective structure} on the  Riemann surface $X^\mr{an}$ associated to  $X$.
 That is, it corresponds to an atlas of coordinate charts  on $X^\mr{an}$ into the complex projective line $\mbP^1 (\mbC)$ such that the transition maps are M\"{o}bius transformations.
Similarly,
  a generic Miura $\mathrm{PGL}_2$-oper corresponds to an {\it affine structure}.

\LSP
\subsection{Opers in positive characteristic} \label{SSb}

The study of $G$-opers in {\it prime characteristic $p > 0$}, particularly for $G = \mr{PGL}_n$ (or $\mr{GL}_n$, $\mr{SL}_n$)
   originated in the context of the $p$-adic Teichm\"{u}ller theory, as  developed by S. Mochizuki (cf. ~\cite{Moc1},  \cite{Moc2}).
   Since then, this topic has been explored in various works, including   ~\cite{BeTr}, ~\cite{JRXY}, ~\cite{JoPa}, ~\cite{LaPa}, ~\cite{LiOs}, ~\cite{Oss1}, and ~\cite{Oss2}.
 A key ingredient common to  these developments is the  study of $p$-curvature, which serves as an  
 invariant measuring  the obstruction to the compatibility of $p$-power structures appearing in certain associated spaces of infinitesimal symmetries.
A $G$-oper is said to be  {\it dormant} if its $p$-curvature  vanishes.
 
 Returning to 
 the $\mr{GL}_n$-oper $\msF$ as in \eqref{Eq21ke},  the  $p$-curvature of $\nabla$ is given by 
 \begin{align}
 \psi (\nabla) : \mcT_X^{\otimes p} \rightarrow \mcE nd_{\mcO_X} (\mcF); \ \partial^{\otimes p} \mapsto (\nabla_{\partial})^p - \nabla_{\partial^p},
 \end{align}
 where, for any local section $\partial \in \mcT_X \left(= \Omega_X^\vee \right)$, we define  $\nabla_{\partial} := (\partial \otimes \mr{id}_\mcF)\circ \nabla$,  and  $\partial^p$ denotes  the $p$-th power iteration of $\partial$ (considered as a derivation).
According to a classical result by Cartier,   $\msF$  is dormant if and only if  the solution space of the  differential  equation $D f= 0$ associated with  the corresponding  differential operator $D$ attains  the highest possible rank.

K. Joshi and C. Pauly proved  that the number of  dormant $\mr{PGL}_n$-opers on a fixed hyperbolic curve  is finite (cf. ~\cite{JoPa}, or ~\cite{BeTr}, ~\cite{Wak5}).
Subsequently,  Joshi proposed a counting problem for these objects, along with an explicit conjectural formula.  
 To address this problem,
 the moduli theory of dormant $G$-opers and  dormant Miura $G$-opers, on (pointed) {\it hyperbolic} curves for general $G$ has been developed  by the second author  (cf. \cite{Wak2}, \cite{Wak7}, \cite{Wak6}, \cite{Wak5},   \cite{Wak8}, \cite{Wak9}, and \cite{Wak14}).
 
 One of the most significant  properties on the moduli space of dormant $G$-opers
 is the {\it generic \'{e}taleness}.
This property is particularly useful for analyzing  how the moduli space behaves  
 under 
deformations and degenerations of the underlying  curves.
Furthermore, it facilitates connections with other areas of enumerative geometry,  such as  Gromov-Witten theory, the conformal field theory of affine Lie algebras, and the combinatorics  of polytopes and graphs (cf. ~\cite{LiOs}, ~\cite{Wak5}, ~\cite{Wak8}).
 It is worth noting that the generic \'{e}taleness has been established  for classical algebraic groups of  types $A_n$, $B_n$, $C_n$ and $D_n$ when  $n$ is sufficiently small relative to $p$ (cf.  ~\cite{Wak14} and ~\cite{Wak5}).
 However,  a complete understanding of the geometric structure of the moduli space for  general $G$ remains an open problem   due to various complications arising from hyperbolicity.

\LSP
\subsection{Results of this paper I: Moduli space of dormant opers on elliptic curves} \label{SSc}

With this background in mind, this paper focuses on the case where the underlying curve is {\it elliptic}, and 
we examine the global structure of the relevant moduli spaces in a more explicit manner.

Let $k$ be an algebraically closed field of characteristic $p$, and 
let $\overline{\mathcal{M}}_{\mathrm{ell}}$ denote  the Deligne-Mumford compactification of the moduli stack classifying elliptic curves (equipped with a distinguished point) over $k$. 
This stack has an open (resp., a closed) substack
\begin{align}
\overline{\mathcal{M}}_{\mathrm{ell}}^{\mathrm{ord}} \quad  \left(\mathrm{resp.} \   \overline{\mathcal{M}}_{\mathrm{ell}}^{\mathrm{ss}}\right),
\end{align}
which classifies  ordinary (resp., supersingular) elliptic curves.

Suppose that $G$ is a simple algebraic group of adjoint type over $k$  whose Coxeter number is sufficiently small relative to $p$.
One can then consider   the moduli stack of dormant $G$-opers over $\overline{\mcM}_\mr{ell}$, given by 
\begin{align}
\mathcal{O}p_G^{^\mathrm{Zzz\ldots}}:=
\left\{(X, \mathscr{F}) \,\big|\, X \in \overline{\mathcal{M}}_{\mathrm{ell}},\, \mathscr{F} \text{ is a dormant } G\text{-oper on } X\right\}
\end{align}
(cf. \eqref{eqwq2}).
Similarly, the stack of dormant generic Miura $G$-opers over $\overline{\mcM}_\mr{ell}$ is defined  as
\begin{align}
\mathcal{M}\mathcal{O}p_G^{^\mathrm{Zzz\ldots}} :=
\left\{(X, \mathscr{F}) \,\big|\, X \in \overline{\mathcal{M}}_{\mathrm{ell}},\, \mathscr{F} \text{ is a dormant generic Miura } G\text{-oper on } X\right\}.
\end{align}
 By forgetting the (Miura) $G$-oper structures, we obtain  
  natural projections
\begin{align}
\Pi_G : \mathcal{O}p_G^{^\mathrm{Zzz\ldots}} \rightarrow \overline{\mathcal{M}}_{\mathrm{ell}},
\hspace{10mm}
\widehat{\Pi}_G : \mathcal{M}\mathcal{O}p_G^{^\mathrm{Zzz\ldots}}\rightarrow \overline{\mathcal{M}}_{\mathrm{ell}}
\end{align}
(cf. \eqref{Eq40}).

These moduli spaces and morphisms between them play a central role in our discussion, involving families of curves over fairly general base schemes.
Such relative formulations enable  us to investigate  deformations of various types of structures.
The following Theorem \ref{Th4992} describes a fundamental property of  $\mcO p^{^\mr{Zzz...}}_G$ and $\mcM \mcO p^{^\mr{Zzz...}}_G$.
(Corresponding results in different settings can be found in ~\cite{BeTr}, ~\cite{Moc2}, ~\cite{JoPa}, and ~\cite{Wak5}.
Additionally, the $p$-curvatures of $\mr{SL}_2$-opers on elliptic curves were examined in ~\cite{Chu}.)

\SSP
\begin{intthm}[cf. Theorem \ref{Th12}] \label{Th4992}
The  stack $\mcO p^{^\mr{Zzz...}}_G$ (resp., $\mcM \mcO p^{^\mr{Zzz...}}_G$) is finite over $\overline{\mcM}_\mr{ell}$ and thus  proper over $k$.
Moreover, the fiber of $\Pi_G$ (resp., $\widehat{\Pi}_G$) over any $k$-rational point in $\overline{\mcM}_\mr{ell}^\mr{ss}$ consists of   a single point.
\end{intthm}
\SSP

The next main result  concerns  the relationship between the ordinariness of elliptic curves and 
the \'{e}tale loci  of  $\mcO p^{^\mr{Zzz...}}_G$ and  $\mcM \mcO p^{^\mr{Zzz...}}_G$.
To state it, 
let $\mft$ denote  the  Lie algebra of
  a fixed  maximal torus $T \subset G$ that is split over $\mbF_p := \mbZ/p\mbZ$, and 
 let  $\mathfrak{t}_{\mathrm{reg}}$ be   its  Zariski open subset consisting of regular elements.
The set  of $\mbF_p$-rational points $\mathfrak{t}_{\mathrm{reg}}(\mbF_p)$  carries  a natural action of 
   the Weyl group $W$ of $(G,T)$, and
   its quotient $\mathfrak{t}_{\mathrm{reg}}(\mbF_p)/W$  is a finite set.

Additionally,  
there exists a natural morphism of $\overline{\mcM}_\mr{ell}$-stacks
\begin{align}
\Xi_G : \mathcal{M}\mathcal{O}p_G^{^\mathrm{Zzz\ldots}} \rightarrow \mathcal{O}p_G^{^\mathrm{Zzz\ldots}}.
\end{align}
For $G = \mr{PGL}_2$,
this morphism corresponds to the {\it Miura transformation}, which relates  
 solutions of the {\it KdV equation} to those of  its modified counterpart, the {\it mKdV equation}.
The following result describes     the symmetry of this morphism over $\overline{\mcM}_\mr{ell}^\mr{ord}$.

\SSP
\begin{intthm}[cf. Theorems \ref{Th7}] \label{Th77} 
\begin{itemize}
\item[(i)]
The projection $\Pi_G : \mcO p_{G}^{^\mr{Zzz...}} \!\rightarrow \overline{\mcM}_{\mr{ell}}$
restricts to a surjective, finite, and   \'{e}tale morphism
\begin{align} \label{Eq152}
\Pi_G^\mr{ord} : \overline{\mcM}_\mr{ell}^\mr{ord} \times_{\overline{\mcM}_\mr{ell}}\mcO p_{G}^{^\mr{Zzz...}} \!\rightarrow \overline{\mcM}_\mr{ell}^\mr{ord}.
\end{align}
Moreover,  its 
 degree coincides with  $\sharp (\mft_\mr{reg}(\mbF_p)/W) \left(=\sharp (\mft_\mr{reg}(\mbF_p))/\sharp (W)\right)$.
\item[(ii)]
The restriction
\begin{align}
\Xi^\mr{ord}_G  : \overline{\mcM}^\mr{ord}_\mr{ell} \times_{\overline{\mcM}_\mr{ell}}
\mcM \mcO p_G^{^\mr{Zzz...}} \rightarrow \overline{\mcM}^\mr{ord}_\mr{ell} \times_{\overline{\mcM}_\mr{ell}}\mcO p_G^{^\mr{Zzz...}}
\end{align}
of  the projection $\Xi_G : \mcM \mcO p_G^{^\mr{Zzz...}} \rightarrow \mcO p_G^{^\mr{Zzz...}}$  over $\overline{\mcM}^\mr{ord}_\mr{ell}$ is  an \'{e}tale Galois covering 
with Galois group $W$.
\end{itemize}
\end{intthm}
\SSP

By applying  Theorems \ref{Th4992} and \ref{Th77}, we obtain  the following assertion.

\SSP
\begin{intthm}[cf. Corollary \ref{Cor22}] \label{ThF}
Both $\mcO p_G^{^\mr{Zzz...}}$ and $\mcM\mcO p_G^{^\mr{Zzz...}}$ are connected.
\end{intthm}

\LSP
\subsection{Results of this  paper II: Canonical diagonal liftings for ordinary elliptic curves} \label{SSd}

This second part of this paper extends some of the previous discussions for $G = \mr{PGL}_n$ in two directions.
The first  extension  considers  the case of {\it prime-power characteristic}, while the second  deals with   a {\it finite-level} version, formulated using  the sheaf of differential operators of finite level (cf. ~\cite{PBer1}, ~\cite{PBer2}).
In particular, for a positive integer $\N$, we define the notion of a {\it dormant $\mr{PGL}_n^{(\N)}$-oper}  (cf. Definition \ref{Def79}), which coincides with the usual definition of a dormant oper when  $\N =1$.
Dormant $\mr{PGL}_n$-opers on hyperbolic curves in these settings, along with  their  moduli spaces, have been  studied in ~\cite{Wak18}, ~\cite{Wak12}.
 This paper extends such considerations to the case of elliptic curves.

Let $X_\N$ be an  elliptic curve   over the ring of Witt vectors $W_\N$ of length $\N$ over $k$.
Denote by $\overline{\tau}_\mr{reg}(W_\N)$ the equivalence classes  of $n$-tuples $(a_1, \cdots, a_n) \in W_\N^{\oplus n}$ whose components  modulo $p$ are mutually distinct, where two such $n$-tuples are considered equivalent if they differ by a translation  by an element of $W_\N$ (cf. \eqref{eqw22}).
The symmetric group of $n$ letters $\mfS_n$ acts on this set by permutation, and we obtain the quotient set $\overline{\tau}_\mr{reg}(W_\N)/\mfS_n$.
When $\N =1$,   this set coincides with $\mft_\mr{reg}(\mbF_p)/W$ for $G = \mr{PGL}_n$.

Next, we shall write
\begin{align}
\mcO p_{\N, 1}^{^\mr{Zzz...}} \ \left(\text{resp.,} \ \mcO p_{1, \N}^{^\mr{Zzz...}}  \right)
\end{align}
for the set of dormant $\mr{PGL}_n^{(1)}$-opers on $X_\N$  (resp., dormant $\mr{PGL}_n^{(\N)}$ on $X_1$).
A key technique to relate these objects 
 is {\it diagonal reduction}, introduced  in ~\cite{Wak12}, which allows us to construct dormant $\mr{PGL}_n^{(\N)}$-oper on $X_1$ from dormant $\mr{PGL}_n^{(1)}$-oper on $X_\N$.
This process defines a map 
\begin{align}
\Diag_{\!\!\N} : \mcO p_{\N, 1}^{^\mr{Zzz...}} \rightarrow \mcO p_{1, \N}^{^\mr{Zzz...}} 
\end{align}
(cf. \eqref{Eq1123}).
In ~\cite[Theorem D]{Wak12}, the second author proved that for $n=2$, this map is bijective when the underlying curve $X$ is replaced with a general hyperbolic curve.
We   establish an analogous result for elliptic curves, proving  
 the bijectivity for arbitrary  $n$ (sufficiently small relative to $p$) under the assumption that  $X_1$ is ordinary.
 The result is summarized as  follows (we also obtain a corresponding statement  for Miura $\mr{PGL}_n$-opers).

\SSP
 \begin{intthm}[cf. Propositions \ref{T39}, \ref{T34}, and Theorem \ref{T47}] \label{ThmB}
Suppose that $X_1$ is ordinary.
Then, there  exist canonical bijections
 \begin{align}
 \Theta_{\N, 1} : \overline{\tau}_\mr{reg} (W_\N (\mbF_p))/\mfS_n
 \xrightarrow{\sim} \mcO p_{\N, 1}^{^\mr{Zzz...}},
 \hspace{10mm}
 \Theta_{1, \N} : \overline{\tau}_\mr{reg} (W_\N (\mbF_p))/\mfS_n
 \xrightarrow{\sim} \mcO p_{1, \N}^{^\mr{Zzz...}}
 \end{align}
 such that  the following diagram commutes:
  \begin{align} \label{Eq1126}
\vcenter{\xymatrix@C=6pt@R=36pt{
 & \overline{\tau}_\mr{reg} (W_\N (\mbF_p))/\mfS_n\ar[ld]^-{\sim}_-{\Theta_{\N, 1}}\ar[rd]^-{\Theta_{1, \N}}_-{\sim}& \\
 \mcO p_{\N, 1}^{^\mr{Zzz...}}\ar[rr]_-{\Diag_{\!\!\N}}
 && \mcO p_{1, \N}^{^\mr{Zzz...}}.
 }}
\end{align}
In particular, the map $\Diag_{\!\!\N}$ is bijective, that is, each dormant $\mr{PGL}_n^{(\N)}$-oper on $X_1$  arises uniquely, via diagonal reduction,  from a  dormant  $\mr{PGL}_n^{(1)}$-oper on $X_\N$
 (which we refer to as  its  {\bf canonical diagonal lifting}).
 \end{intthm}
\SSP


Finally, we make a few  remarks on our results.
In \cite{Wak5},  
Joshi's conjectural formula  for computing  the total number of dormant $\mr{PGL}_n$-opers was proven.
  Extending this  formula to more general cases remains a central problem  in the enumerative geometry of dormant opers. 
  Theorems \ref{Th77} and \ref{ThmB}
 generalize this result   to general $G$ and $\N$  under the assumption that the underlying curves are elliptic.
 This is also expected to be useful in the hyperbolic case, as dormant opers on elliptic curves corresponds to those on $1$-pointed genus-$1$ curves  with  a specific {\it radius} at the  marked point (cf. ~\cite{Moc2}, ~\cite{Wak5}, or ~\cite{Wak12} for the definition of radius).

Furthermore,  we hope that this research will involve topics related to the geometric Langlands correspondence in positive characteristic, as explored in ~\cite{ChZh}, ~\cite{BeBr}, ~\cite{She}, and ~\cite{Tra}.

\vspace{5mm}

\begin{center}
\begin{tikzpicture}[>=Stealth, font=\small]

\node[draw, rounded corners=8pt, thick, minimum width=5.5cm, minimum height=2.0cm, align=center] (left) at (-2,0) {
  Dormant $\mathrm{PGL}_n^{(1)}$-opers on\\
  a lifting $X_N$ of $X_1$ to char.\,$p^N$
};

\node[draw, rounded corners=8pt, thick, minimum width=5.5cm, minimum height=2.0cm, align=center] (right) at (9,0) {
  Dormant $\mathrm{PGL}_n^{(N)}$-opers on\\
  an ordinary elliptic curve $X_1$\\
  in char.\,$p$
};

\draw[->, thick] (1.3,0.5) -- ++(4.2,0) node[midway, above]{Diagonal Reduction};

\draw[<-, thick] (1.3,-0.5) -- ++(4.2,0) node[midway, below]{Diagonal Lifting};

\end{tikzpicture}
\end{center}

\LSP
\subsection{Notation and Conventions} \label{SS0}

Throughout  this paper, we fix a prime number $p$ and an algebraically closed  field $k$  of characteristic $p$.
All schemes appearing in this paper are assumed to be locally noetherian.

Let $S$ be a scheme over 
 $\mbF_p := \mbZ/p\mbZ$
 and $f : X \rightarrow S$ a scheme over $S$.
 Denote by $F_S : S \rightarrow S$ (resp., $F_X : X \rightarrow X$)
 the absolute Frobenius endomorphism of $S$ (resp., $X$).
 The {\bf Frobenius twist} of $X$ over $S$ is, by definition,
 the base-change
 \begin{align}
 X^{(1)} \left(:= S \times_{F_S, S} X \right)
 \end{align}
 of $X$ along 
   $F_S$.
 Denote by $f^{(1)} : X^{(1)} \rightarrow S$ the structure morphism of $X^{(1)}$.
 The {\bf relative Frobenius morphism} of $X$ over $S$ is the unique 
 morphism $F_{X/S} : X \rightarrow X^{(1)}$ over $S$ that fits into the following commutative diagram:
 \begin{align} \label{QS3360}
\vcenter{\xymatrix@C=36pt@R=36pt{
X \ar@/^10pt/[rrrrd]^{F_X}\ar@/_10pt/[ddrr]_{f} \ar[rrd]_{  \ \ \ \  \  \   F_{X/S}} & & &   &   \\
& & X^{(1)}  \ar[rr]_{F_S \times \mr{id}_X} \ar[d]^-{f^{(1)}}  \ar@{}[rrd]|{\Box}  &  &  X \ar[d]^-{f} \\
&  &S \ar[rr]_{F_S} & &  S.
}}
\end{align}

A {\bf log scheme} means, in the present paper, a scheme equipped with a logarithmic structure in the sense of Fontaine-Illusie (cf. ~\cite[(1.2)]{KKa}).
For a log scheme (resp., a morphism of log schemes) indicated, say, by  ``$X^\mr{log}$" (resp., ``$f^\mr{log}$"), we shall use the notation ``$X$" (resp., ``$f$") to denote the underlying scheme (resp., the underlying morphism of schemes).
Recall that a log scheme is called {\bf fs} (i.e., {\bf fine and saturated}, in the sense of  ~\cite[\S\,1.3]{Ill}) if it admits an \'{e}tale local  chart modeled on a finite (= integral and finitely generated) and saturated monoid.

Given a morphism of log schemes $f^\mr{log} : X^\mr{log} \rightarrow S^\mr{log}$, we shall write $\mcT_{X^\mr{log}/S^\mr{log}}$ for the sheaf of logarithmic derivations on $X^\mr{log}$ relative to $S^\mr{log}$ and $\Omega_{X^\mr{log}/S^\mr{log}}$ for its dual, i.e., the sheaf of logarithmic $1$-forms on $X^\mr{log}$ relative to $S^\mr{log}$ (cf. ~\cite[(1.7)]{KKa}).
 Note that various notions and discussions in logarithmic geometry of schemes can be naturally generalized to  Deligne-Mumford stacks.
For instance,  one can define fs log stacks and log curves (in the sense of \S\,\ref{SS2}) over such stacks.

Denote by $\mbG_m$ the multiplicative group over $k$.
For a positive integer $n$,
we shall write $\mr{GL}_n$ (resp., $\mr{PGL}_n$) for the general (resp., projective) linear group of rank $n$ over $k$.

\LSP
\section{Opers and Miura opers on elliptic curves} \label{S2}
\LSP

This section deals with  the basics of  opers and Miura opers defined on stable elliptic curves in terms of logarithmic geometry.
We discuss  the representability of their  moduli spaces (cf. Corollaries \ref{Cor9} and \ref{Cor33}). 
These results  are obtained by applying, almost verbatim, the arguments developed  for hyperbolic curves.
However,  a key distinction lies in the presence of nowhere-vanishing tangent fields  on elliptic curves, which  allows for a significantly simpler description of both  (generic Miura) opers and their moduli space.
Although we assume that the base field $k$ has positive characteristic, this assumption is not essential at least in this section.
For discussions on opers and Miura opers in characteristic zero (without the logarithmic setting),
 we refer the reader to ~\cite{BeDr1}, ~\cite{Fre1}, ~\cite{Fre2}, ~\cite{Fre3}, and related works,  where these objects are studied  in the context of the geometric Langlands program.

\LSP
\subsection{Stable elliptic curves and their moduli stack} \label{SS2}

Let $S$ be a scheme.
By a {\bf curve} over $S$, we mean a proper flat scheme $f : X \rightarrow S$ over $S$ all of whose geometric fibers are connected, reduced, and of dimension $1$.
An {\bf $r$-pointed curve} over $S$ (where $r \in \mbZ_{> 0}$) is
a collection of data
\begin{align}
\msX := (f: X \rightarrow S, \{ \sigma_i \}_{i=1}^r)
\end{align}
consisting of 
  a curve $f : X \rightarrow S$ and an ordered set of $r$ sections $\sigma_i : S \rightarrow X$ ($i=1, \cdots, r$) of $f$ whose images are disjoint, i.e., $\mr{Im}(\sigma_i) \cap \mr{Im}(\sigma_j) = \emptyset$ if $i \neq j$. 
An $r$-pointed curve $\msX := (f : X \rightarrow S, \{ \sigma_i \}_i)$ is  called {\bf of (arithmetic) genus $g \in \mbZ_{\geq 0}$} if, for any $k$-rational  point $s : \mr{Spec}(k) \rightarrow S$ of $S$, the fiber $X_s := f^{-1} (s)$ satisfies the equality $\mr{dim}_{k} (H^1 (X_s, \mcO_{X_s}))= g$.

Suppose that, for each $\Box \in \{\circ, \bullet \}$,   we are given a scheme $S_\Box$ and an $r$-pointed curve $\msX_\Box := (f_\Box : X_\Box \rightarrow S_\Box, \{ \sigma_{\Box i} \}_{i=1}^r)$ over $S_\Box$.
Then, a {\bf morphism of $r$-pointed curves} from $\msX_\circ$ to $\msX_\bullet$ is defined as a pair of morphisms between schemes $(\phi : S_\circ \rightarrow S_\bullet, \Phi : X_\circ \rightarrow X_\bullet)$ such that the following squares form cartesian commutative diagrams:
\begin{align}
\hspace{20mm}
\vcenter{\xymatrix@C=46pt@R=36pt{
 X_\circ \ar[r]^-{\Phi} \ar[d]_-{f_\circ} & X_\bullet \ar[d]^-{f_\bullet} \\
 S_\circ \ar[r]_-{\phi} & S_\bullet,
 }}
 \hspace{15mm}
 \vcenter{\xymatrix@C=46pt@R=36pt{
 X_\circ \ar[r]^-{\Phi}  & X_\bullet  \\
 S_\circ \ar[r]_-{\phi} \ar[u]^-{\sigma_{\circ i}} & S_\bullet \ar[u]_-{\sigma_{\bullet i}}
 }} \ \ \ (i=1, \cdots, r).
\end{align}

For a pair of nonnegative integers $(g, r)$ with $2g-2 +r > 0$,
an {\bf $r$-pointed stable curve of genus $g$} over $S$ is defined as an $r$-pointed curve $\msX := (f : X \rightarrow S, \{\sigma_i \})$ of genus $g$ such that, for each geometric  point $s : \mr{Spec} (k') \rightarrow S$ of $S$, any $k'$-rational point of $X_s$ is either a smooth point or a nodal point (i.e., a point $x$ such that the complete local ring $\widehat{\mcO}_{X_s, x}$ is isomorphic to $k' [\![u, v ]\!]/(uv)$).
Then, $r$-pointed stable curves of genus $g$ over $k$-schemes  and $k$-morphisms between them together  form a category
$\overline{\mcM}_{g, r}$ fibered in groupoids over the category of $k$-schemes $\mcS ch_{/k}$ (cf. ~\cite[Definition 1.1]{Knu}).
It is well-known that $\overline{\mcM}_{g, r}$ can be represented by a connected, proper, and smooth Deligne-Mumford stack over $k$ of dimension $3g-3+r$ (cf. ~\cite[Corollary 2.6 and Theorem 2.7]{Knu}, ~\cite[\S\,5]{DeMu}).

For convenience of our discussion,
we shall write
\begin{align}
\overline{\mcM}_{\mr{ell}} := \overline{\mcM}_{1, 1},
\end{align}
and a pointed curve classified by this stack will be referred to as a {\bf stable elliptic curve}.
Note that $\overline{\mcM}_\mr{ell}$ has a dense open substack $\mcM_{\mr{ell}}$ classifying smooth curves (which are called {\bf smooth elliptic curves}, or simply, {\bf elliptic curves}).
Also, write $f_\mr{ell} : \mcC_{\mr{ell}} \rightarrow \overline{\mcM}_\mr{ell}$ for
the universal family of curves over $\overline{\mcM}_\mr{ell}$.

The divisor at infinity $\partial \overline{\mcM}_\mr{ell} := \overline{\mcM}_\mr{ell} \setminus \mcM_\mr{ell}$ determines a log structure on $\overline{\mcM}_{\mr{ell}}$, and
the pull-back of $\partial \overline{\mcM}_\mr{ell}$ via  $f_\mr{ell}$ determines a log structure on $\mcC_\mr{ell}$.
We denote the resulting fs log stacks by $\overline{\mcM}_\mr{ell}^\mr{log}$ and $\mcC_\mr{ell}^\mr{log}$, respectively.
(Note that the log structure of $\mcC^\mr{log}_\mr{ell}$ mentioned here is different from that of ~\cite{FKa}.
In fact, the log structure defined in that reference arises from the divisor which is the union of the pull-back of $\partial\overline{\mcM}_{\mr{ell}}$ {\it and  the image of the universal  marked point $\overline{\mcM}_\mr{ell} \rightarrow \mcC_\mr{ell}$}.)
In particular, $f_\mr{ell}$ extends naturally to a morphism of  fs log stacks $f_\mr{ell}^\mr{log} : \mcC_\mr{ell}^\mr{log} \rightarrow \overline{\mcM}_\mr{ell}^\mr{log}$, forming a family of log curves.
Here, a {\bf log curve}  means a log smooth integral morphism of fs log schemes (or more generally, of fs log stacks)  $f^\mr{log} : X^\mr{log} \rightarrow S^\mr{log}$  such that 
the underlying morphism of schemes $f : X \rightarrow S$ is proper and 
all of its geometric fibers are   connected and reduced $1$-dimensional schemes.
(Both   $\mcT_{X^\mr{log}/S^\mr{log}}$ and  $\Omega_{X^\mr{log}/S^\mr{log}}$ associated to such a log curve  are line bundles,
 and the underlying morphism $f : X \rightarrow S$ is flat;  see ~\cite[Corollary 4.5]{KKa}).

Next, let $S$ be a $k$-scheme and $\msX := (f : X \rightarrow S, \sigma : S \rightarrow X)$ be a stable elliptic curve  over $S$.
It corresponds to a morphism $\phi : S \rightarrow \overline{\mcM}_\mr{ell}$ with $X \cong S \times_{\phi, \overline{\mcM}_\mr{ell}} \mcC_\mr{ell}$.
Pulling-back the log structures of $\overline{\mcM}_\mr{ell}^\mr{log}$ and $\mcC_{\mr{ell}}^\mr{log}$ yields
 log structures on $S$ and $X$, respectively;
 we denote the resulting log schemes by $S^\mr{log}$ and $X^\mr{log}$, respectively.
 Both $S^\mr{log}$ and $X^\mr{log}$ are fs log schemes, and the morphism $f: X \rightarrow S$ extends to a log curve $f^\mr{log} : X^\mr{log} \rightarrow S^\mr{log}$ over $S^\mr{log}$.
(When considering an {\it $S^\mr{log}$-connection} in the sense of the subsequent discussion, we are essentially dealing with a connection over  {\it unpointed} genus-$1$ curves, as   it does not permit  a  logarithmic   pole at  the unique marked point.)
The canonical  bundle $\omega_{X/S}$ of $X/S$ is naturally isomorphic to $\Omega_{X^\mr{log}/S^\mr{log}}$.
Note that if $\msX$ is smooth, then we have $S^\mr{log}  = S$ and $X^\mr{log} = X$.

\LSP
\subsection{Logarithmic connections on principal bundles} \label{SS3}

Let $S^\mr{log}$ be an fs log scheme over $k$ and $f^\mr{log} : X^\mr{log} \rightarrow S^\mr{log}$ a log curve over $S^\mr{log}$.
For simplicity, we write $\mcT := \mcT_{X^\mr{log}/S^\mr{log}}$ and $\Omega := \Omega_{X^\mr{log}/S^\mr{log}}$.
Also, let $G$ be a connected smooth affine algebraic group over $k$ with Lie algebra $\mfg$ and 
$\pi : \mcE \rightarrow X$  a $G$-bundle on $X$.

Given a $k$-vector space $\mfh$ equipped with a $G$-action, we write $\mfh_G$ for the vector bundle on $X$ associated to the affine space $\mcE \times^G \mfh \left(= (\mcE \times_k \mfh)/G \right)$ over $X$.
In particular, by considering the $k$-vector space $\mfg$ equipped with the $G$-action given by the adjoint representation $\mr{Ad}_G : G \rightarrow \mr{GL} (\mfg)$, we obtain the {\bf  adjoint bundle} $\mfg_\mcE$ associated to $\mcE$.

The pull-back of the log structure on $X^\mr{log}$ via $\pi$ defines a log structure on $\mcE$; we denote the resulting log scheme by $\mcE^\mr{log}$.
The $G$-action on $\mcE$ induces a $G$-action on the direct image $\pi_* (\mcT_{\mcE^\mr{log}/S^\mr{log}})$ 
of $\mcT_{\mcE^\mr{log}/S^\mr{log}}$, and
 we obtain the subsheaf of $G$-invariant sections  $\widetilde{\mcT}_{\mcE^\mr{log}/S^\mr{log}} := \pi_* (\mcT_{\mcE^\mr{log}/S^\mr{log}})^G$.
Since the adjoint bundle $\mfg_\mcE$ coincides with the sheaf of $G$-invariant sections in  $\pi_* (\mcT_{\mcE^\mr{log}/X^\mr{log}}) \left(= \pi_* (\mcT_{\mcE/X}) \right)$,  the differential  of $\pi$ yields 
a short exact sequence of $\mcO_X$-modules
\begin{align} \label{Eq1}
0 \rightarrow \mfg_\mcE \rightarrow \widetilde{\mcT}_{\mcE^\mr{log}/S^\mr{log}} \xrightarrow{d_\mcE} \mcT \rightarrow 0. 
\end{align}
We shall regard $\mfg_\mcE$ as an $\mcO_X$-submodule of $\widetilde{\mcT}_{\mcE^\mr{log}/S^\mr{log}}$ via the second arrow $\mfg_\mcE \hookrightarrow \widetilde{\mcT}_{\mcE^\mr{log}/S^\mr{log}}$.

Recall from ~\cite[Definition 1.17]{Wak5} that an {\bf $S^\mr{log}$-connection} on $\mcE$ is a split injection of  \eqref{Eq1}, i.e., an $\mcO_X$-linear morphism
\begin{align}
\nabla : \mcT  \rightarrow \widetilde{\mcT}_{\mcE^\mr{log}/S^\mr{log}}
\end{align}
with  $d_\mcE \circ \nabla = \mr{id}_{\mcT}$.
Since $\mcT$ is a line bundle,  each  $S^\mr{log}$-connection  is automatically flat, i.e.,  its curvature vanishes identically (cf. ~\cite[Definition 1.23]{Wak5} for the definition of the curvature of a logarithmic connection).
By a {\bf flat $G$-bundle} on $X^\mr{log}/S^\mr{log}$, we mean a pair $(\mcE, \nabla)$ consisting of a $G$-bundle $\mcE$ on $X$ and an $S^\mr{log}$-connection $\nabla$ on $\mcE$.

Let $\nabla$ be  an $S^\mr{log}$-connection  on $\mcE$ and 
 $\eta : \mcE \xrightarrow{\sim} \mcE$  a $G$-equivariant  automorphism of $\mcE$ over $X$.
 The differential of $\eta$ gives a $G$-equivariant automorphism  of $
\pi_*(\mcT_{\mcE^\mr{log}/S^\mr{log}})$, and hence restricts to an automorphism $d\eta^G$ of  $\widetilde{\mcT}_{\mcE^\mr{log}/S^\mr{log}}$.
Then,  the composite
\begin{align}
\eta_* (\nabla) := d\eta^G \circ \nabla : \mcT \rightarrow \widetilde{\mcT}_{\mcE^\mr{log}/S^\mr{log}}
\end{align}
defines  an $S^\mr{log}$-connection, which is  called  the {\bf gauge transformation}
 of $\nabla$ by $\eta$ (cf. ~\cite[Definition 1.20]{Wak5}).
Two $S^\mr{log}$-connections on $\mcE$ are called  {\bf gauge equivalent} if they coincide with each other after applying a gauge transformation.
The binary relation of ``{\it being gauge equivalent}" specifies  an equivalence relation in the set of $S^\mr{log}$-connections.

\LSP
\subsection{Connections on the trivial principal bundle} \label{SS42}

The {\bf trivial $G$-bundle} is the product $\mcE_{G, \mr{triv}} := X \times_k G$, where $G$ acts only on the second factor, considered as  a scheme over $X$ by the first projection $X \times_k G \rightarrow X$.
For an element
$h$ of $G (S)$,
the left-translation by $h$ determines an automorphism $L_h : \mcE_{G, \mr{triv}} \xrightarrow{\sim} \mcE_{G, \mr{triv}}$ (cf. ~\cite[\S\,1.1.2]{Wak5}).
Conversely, since  the equality $G (S) = G (X)$ holds because of the affinity of $G$, any automorphism of $\mcE_{G, \mr{triv}}$ coincides with $L_h$ for a unique $h \in G (S)$.

Note that there exists a canonical decomposition
\begin{align} \label{Eq61}
\widetilde{\mcT}_{\mcE_{G, \mr{triv}}^\mr{log}/S^\mr{log}} = \mcT \oplus \mcO_X \otimes_k \mfg
\end{align}
(cf. ~\cite[Example 1.12]{Wak5}), where the  direct summand $\mcO_X \otimes_k \mfg$ in the right-hand side is identified with the adjoint bundle $\mfg_{\mcE_{G, \mr{triv}}}$.
In particular, 
if  $\mr{pr}_2$ denotes the projection onto the second factor $\mcT \oplus \mcO_X \otimes_k \mfg \twoheadrightarrow \mcO_X \otimes_k \mfg$, then
each $S^\mr{log}$-connection $\nabla$ on $\mcE_{G, \mr{triv}}$ can be determined by the associated morphism 
\begin{align} \label{Eq65}
\nabla^\mfg := \mr{pr}_2 \circ \nabla : \mcT \rightarrow \mcO_X \otimes_k \mfg.
\end{align}

Next, suppose that  we are given a generator $\partial$ of $\mcT$, i.e., $\mcT = \mcO_X \partial$.
Given an element $v$ of $\mfg (S) \left(= H^0 (X, \mcO_X \otimes_k \mfg) \right)$, we obtain  an $S^\mr{log}$-connection 
\begin{align} \label{Eq60}
\nabla_{\partial, v} : \mcT \rightarrow  \mcT \oplus \mcO_X \otimes_k \mfg \left(\stackrel{\eqref{Eq61}}{=} \widetilde{\mcT}_{\mcE^\mr{log}_{G, \mr{triv}}/S^\mr{log}} \right)
\end{align}
 on $\mcE_{G, \mr{triv}}$ given by $a \cdot \partial \mapsto (a \cdot \partial, a \cdot  v)$ for any $a \in \mcO_X$.
 Conversely, each $S^\mr{log}$-connection on $\mcE_{\mr{triv},G}$ coincides with $\nabla_{\partial, v}$ for a unique $v \in \mfg (S)$.

\SSP
\bpr  \label{Prop15}
Let $\partial$ be as above and 
$h$ an element of $G (S)$.
Then, the gauge transformation $L_{h *}(\nabla_{\partial, v})$ of $\nabla_{\partial, v}$ by $L_{h}$ satisfies the equality 
\begin{align} \label{Eq62}
L_{h *}(\nabla_{\partial, v})^\mfg = \mr{Ad}_G (h) \circ \nabla^\mfg
\end{align}
 of morphisms $\mcT \rightarrow \mcO_X \otimes_k \mfg$.
 In particular, 
 the bijection 
 \begin{align}
 \mfg (S) \xrightarrow{\sim}
 \left(\begin{matrix} \text{the set of} \\ \text{$S^\mr{log}$-connections on $\mcE_{G, \mr{triv}}$} \end{matrix} \right)
 \end{align}
 given by  $v \mapsto \nabla_{\partial, v}$ induces  a bijection 
 \begin{align}
 \mfg (S)/G(S) \xrightarrow{\sim} \left(\begin{matrix} \text{the set of gauge equivalence classes} \\ \text{of $S^\mr{log}$-connections on $\mcE_{G, \mr{triv}}$} \end{matrix} \right),
 \end{align}
 where the domain of this bijection is the quotient set of $\mfg (S)$ with respect to the adjoint action of $G (S)$.
\epr
\begin{proof}
 Let  $\Theta_G \in H^0 (G, \Omega_{G/k} \otimes_k \mfg)$ be  the Maurer-Cartan form on $G$ (cf. ~\cite[Definition 1.1]{Wak5}).
The pull-back $(h^{-1})^*(\Theta_G)$ of $\Theta$ via $h^{-1} : X \rightarrow G$ defines an element of $H^0 (X, \Omega \otimes_k \mfg) \left(= \mr{Hom}_{\mcO_X} (\mcT, \mcO_X \otimes_k \mfg) \right)$   via the natural morphism $(h^{-1})^* (\Omega_{G/k}) \rightarrow \Omega$.
Since $h^{-1}$ factors through $f : X \rightarrow S$, we have 
 $(h^{-1})^* (\Theta_G) = 0$.
Hence, 
it follows from ~\cite[Proposition 1.22]{Wak5}
that
 \begin{align}
L_{\eta *}(\nabla_{\partial, v})^\mfg = (h^{-1})^* (\Theta_G) + \mr{Ad}_G (h) \circ \nabla^\mfg =  \mr{Ad}_G (h) \circ \nabla^\mfg,
\end{align}
which completes the proof of the first assertion.
Moreover, the second assertion follows from the first assertion together with the respective  final sentences in the first and third  paragraphs of this subsection.
\end{proof}

\LSP
\subsection{Opers on a log curve} \label{SS4}

Hereinafter, suppose that $G$ is 
a simple algebraic $k$-group of adjoint type satisfying the inequality  $2h_G < p$,
  where $h_G$ denotes the Coxeter number of $G$.
Choose 
a maximal torus $T$ of $G$ and a Borel subgroup $B$ of $G$ containing $T$.
Also, write $\mfb$ and $\mft$ for the Lie algebras of $B$ and $T$, respectively (hence $\mft \subseteq \mfb \subseteq \mfg$).  

Let $\Gamma \subseteq \mr{Hom}(T, \mbG_m)$ (:= the additive group of all characters of $T$) denote the set of simple roots in $B$ with respect to $T$.
For each $\beta \in \mr{Hom}(T, \mbG_m)$, we set
\begin{align}
\mfg^\beta := \left\{ x \in \mfg \, | \, \mr{Ad}_G (t) (x) = \beta (t) \cdot x  \ \text{for all $t \in T$}\right\}.
\end{align}
There exists a unique Lie algebra grading
$\mfg = \bigoplus_{j \in \mbZ} \mfg_j$ (i.e., $[\mfg_{j_1}, \mfg_{j_2}] \subseteq \mfg_{j_1 + j_2}$ for $j_1, j_2 \in \mbZ$) satisfying the following conditions:
\begin{itemize}
\item
$\mfg_{j} = \mfg_{-j} = 0$ for all $j> \mr{rk}(\mfg)$,  where $\mr{rk}(\mfg)$ denotes the rank of $\mfg$;
\item
$\mfg_0 = \mft$, $\mfg_1 = \bigoplus_{\alpha \in \Gamma} \mfg^{\alpha}$, and $\mfg_{-1} = \bigoplus_{\alpha \in \Gamma} \mfg^{-\alpha}$.
\end{itemize}
The associated decreasing filtration $\{ \mfg^j \}_{j \in \mbZ}$  on $\mfg$ is defined by $\mfg^j := \bigoplus_{\ell \geq j} \mfg_\ell$,
which is closed under the adjoint action of $B$ and satisfies the following conditions:
\begin{itemize}
\item
$\mfg^{-j} = \mfg$ and $\mfg^j = 0$ for all   $j > \mr{rk}(\mfg)$;
\item
$\mfg^0 = \mfb$ and $[\mfg^{j_1}, \mfg^{j_2}] \subseteq \mfg^{j_1 + j_2}$ for $j_1, j_2 \in \mbZ$.
\end{itemize}

Now, let $\pi_B : \mcE_B \rightarrow X$ be a $B$-bundle on $X$, and
denote by $\pi_G : \left(\mcE_B \times^B G =: \right) \mcE_G \rightarrow X$ the $G$-bundle obtained via change of structure group by the inclusion $B \hookrightarrow G$. 
Hence, $\mcE_B$ specifies a $B$-reduction of $\mcE_G$, and
the  differential of the natural morphism $\mcE_B \rightarrow \mcE_G$ yields an identification  $\mfg_{\mcE_B} =  \mfg_{\mcE_G}$, as well as   
 the following injective morphism of short exact sequences:
\begin{align} \label{Ed37}
\vcenter{\xymatrix@C=46pt@R=36pt{
 0  \ar[r] & \mfb_{\mcE_B} \ar[r] \ar[d]^-{\iota_{\mfg /\mfb}} & \widetilde{\mcT}_{\mcE_B^\mr{log}/S^\mr{log}}\ar[d]^-{\widetilde{\iota}_{\mfg/\mfb}} \ar[r]^-{d_{\mcE_B}}  & \mcT \ar[r] \ar[d]^-{\mr{id}}_-{\wr}& 0\\
0 \ar[r] & \mfg_{\mcE_G} \ar[r] & \widetilde{\mcT}_{\mcE_G^\mr{log}/S^\mr{log}} \ar[r]_-{d_{\mcE_G}} & \mcT  \ar[r] & 0. 
 }}
\end{align}
Each $\mfg^j \left(\subseteq \mfg \right)$ is closed under the adjoint $B$-action, so induces a subbundle $\mfg^j_{\mcE_B}$ of $\mfg_{\mcE_B} \left(=  \mfg_{\mcE_G} \right)$.
The collection $\{ \mfg^j_{\mcE_B} \}_{j \in \mbZ}$ forms a decreasing filtration on $\mfg_{\mcE_B}$.
Moreover, 
$\widetilde{\mcT}_{\mcE_G^\mr{log}/S^\mr{log}}$ admits 
a decreasing filtration $\{ \widetilde{\mcT}^j_{\mcE_G^\mr{log}/S^\mr{log}} \}_{j \in \mbZ}$ defined by 
\begin{align}
\widetilde{\mcT}^j_{\mcE_G^\mr{log}/S^\mr{log}} := \mr{Im}(\widetilde{\iota}_{\mfg/\mfb}) + \mfg^j_{\mcE_B} \left(\subseteq \widetilde{\mcT}_{\mcE_G^\mr{log}/S^\mr{log}} \right)
\end{align}
for each $j \in \mbZ$.
The equality $\widetilde{\mcT}^j_{\mcE^\mr{log}/S^\mr{log}} = \widetilde{\mcT}^0_{\mcE^\mr{log}/S^\mr{log}}$ holds when  $j \in \mbZ_{\geq 0}$, and the  inclusion $\mfg_{\mcE_B}\hookrightarrow \widetilde{\mcT}_{\mcE_G^\mr{log}/S^\mr{log}}$ induces an isomorphism
\begin{align}
\mfg^{j-1}_{\mcE_B}/\mfg^j_{\mcE_B} \xrightarrow{\sim} \widetilde{\mcT}^{j-1}_{\mcE_G^\mr{log}/S^\mr{log}}/\widetilde{\mcT}^j_{\mcE_G^\mr{log}/S^\mr{log}}
\end{align}
for every $j \leq 0$.
Note that each $\mfg^{-\alpha}$ ($\alpha \in \Gamma$) is closed under the $B$-action defined  by the composite $B \twoheadrightarrow \left(B/[B, B] =\right) T \xrightarrow{\mr{Ad}_G} \mr{Aut} (\mfg^{-\alpha})$.
Hence, the decomposition $\left(\mfg_{-1} = \right) \mfg^{-1}/\mfg^0 = \bigoplus_{\alpha \in \Gamma} \mfg^{-\alpha}$ gives rise to a canonical decomposition 
\begin{align} \label{Eq2}
\widetilde{\mcT}^{-1}_{\mcE_G^\mr{log}/S^\mr{log}}/\widetilde{\mcT}^0_{\mcE_G^\mr{log}/S^\mr{log}} = \bigoplus_{\alpha \in \Gamma} \mfg^{-\alpha}_{\mcE_B}.
\end{align}

\SSP
\bde[cf. ~\cite{Wak5}, Definition 2.1, for the logarithmic setting] \label{Def1}
\begin{itemize}
\item[(i)]
Let us consider a pair 
\begin{align}
\msE^\spadesuit := (\mcE_B, \nabla)
\end{align}
consisting of a $B$-bundle $\mcE_B$ on $X$ and an $S^\mr{log}$-connection $\nabla$ on the  $G$-bundle $\mcE_G$ $\left(:= \mcE_B \times^B G \right)$ associated to  $\mcE_B$.
Then, we shall say that $\msE^\spadesuit$ is a {\bf $G$-oper} on $X^\mr{log}/S^\mr{log}$ if $\nabla (\mcT) \subseteq \widetilde{\mcT}^{-1}_{\mcE_G^\mr{log}/S^\mr{log}}$ and,  for any $\alpha \in \Gamma$, the composite
\begin{align}
\mcT \xrightarrow{\nabla} \widetilde{\mcT}^{-1}_{\mcE_G^\mr{log}/S^\mr{log}} \twoheadrightarrow \widetilde{\mcT}^{-1}_{\mcE_G^\mr{log}/S^\mr{log}}/\widetilde{\mcT}^0_{\mcE_G^\mr{log}/S^\mr{log}} \twoheadrightarrow \mfg_{\mcE_B}^{-\alpha}
\end{align}
is an isomorphism, where the third arrow denotes the natural projection with respect to the decomposition  \eqref{Eq2}.
If $X^\mr{log}/S^\mr{log}$ arises from a stable elliptic curve $\msX$  (as explained  in the final paragraph of \S\,\ref{SS2}), then we refer to any $G$-oper on $X^\mr{log}/S^\mr{log}$ as a {\bf $G$-oper on $\msX_\mr{ell}$}.
\item[(ii)]
Let $\msE_\circ^\spadesuit := (\mcE_{\circ B}, \nabla_\circ)$, $\msE_\bullet^\spadesuit := (\mcE_{\bullet B}, \nabla_\bullet)$ be $G$-opers on $X^\mr{log}/S^\mr{log}$.
An {\bf isomorphism of $G$-opers} from $\msE^\spadesuit_\circ$ to $\msE_\bullet^\spadesuit$ is defined as an isomorphism of $B$-bundles $\mcE_{\circ B} \xrightarrow{\sim} \mcE_{\bullet B}$ such that $\nabla_\circ$ commutes with $\nabla_\bullet$ via 
the associated  isomorphism of $G$-bundles $\mcE_{\circ G} \xrightarrow{\sim} \mcE_{\bullet G}$. 
\end{itemize} 

\ede

For a morphism $(\phi, \Phi) : \msX   \rightarrow \msX'$ in $\overline{\mcM}_\mr{ell}$,
the {\bf base-change} $\phi^*(\msE^\spadesuit)$ of each $G$-oper $\msE^\spadesuit$ on $\msX'$ along  this morphism can be formulated  (cf. ~\cite[\S\,2.1.5]{Wak5}).
Thus, we obtain a fibered category
\begin{align} \label{Eq100}
\mcO p_G
\end{align}
 over  $\mcS ch_{/k}$ defined as follows:
\begin{itemize}
\item
The objects are the pairs $(\msX, \msE^\spadesuit)$, where $\msX$ is 
a stable elliptic curve
 over a $k$-scheme $S$ and $\msE^\spadesuit$ is a $G$-oper on $\msX_\mr{ell}$;
\item
The morphisms from $(\msX, \msE^\spadesuit)$ to $(\msX', \msE'^\spadesuit)$ are morphisms of pointed curves $(\phi, \Phi) : \msX \rightarrow \msX'$   with $\msE^\spadesuit \cong \phi^* (\msE'^\spadesuit)$;
\item
The projection $\mcO p_{G} \rightarrow \mcS ch_{/k}$ is given by assigning, to each pair $(\msX, \msE^\spadesuit)$ as above, the base $k$-scheme $S$ of $\msX$.
\end{itemize}
Since  any $G$-oper does not have nontrivial automorphisms (cf. ~\cite[Proposition 2.9]{Wak5}), 
the functor 
  $\mcO p_G \rightarrow \mcS ch_{/k}$ turns out to be {\it  fibered in  equivalence relations},  in the sense of ~\cite[Definition 3.9]{FGA},
i.e., specifies a set-valued  sheaf on $\mcS ch_{/k}$ (with respect to the big \'{e}tale topology).
Forgetting the $G$-oper structures  yields a functor  $\mcO p_G \rightarrow \overline{\mcM}_{\mr{ell}}$.

\LSP
\subsection{Adjoint quotient of the Lie algebra} \label{SS90}

To proceed our discussion, we shall fix  a collection of data $\{ x_\alpha \}_{\alpha \in \Gamma}$, where each $x_\alpha$ denotes a generator of $\mfg^\alpha$.
(Two such collections are conjugate by an element of $T (k)$, and particularly
 the results obtained in the subsequent discussion are essentially independent of the choice of $\{ x_\alpha \}_\alpha$.)
Let us set
\begin{align}
q_1 := \sum_{\alpha \in \Gamma} x_\alpha, \hspace{10mm}
\check{\rho}:= \sum_{\alpha \in \Gamma} \check{\omega}_\alpha,
\end{align}
where $\check{\omega}_\alpha$ (for each $\alpha \in \Gamma$)
denotes the fundamental coweight of $\alpha$, considered as an element of $\mft$ via differentiation.
Then, there exists a unique collection $\{ y_\alpha \}_{\alpha \in \Gamma}$, where $y_\alpha$ is a generator of $\mfg^{-\alpha}$, such that if we write 
\begin{align}
q_{-1} := \sum_{\alpha \in \Gamma} y_\alpha \in \mfg_{-1},
\end{align}
then the set $\{ q_{-1}, 2\check{\rho}, q_1 \}$ forms an $\mfs \mfl_2$-triple.

Denote by $W$ the  Weyl group of $(G, T)$, and consider the GIT quotient $\mfg /\!\!/G$ (resp., $\mft /\!\!/W$) of $\mfg$ (resp., $\mft$) by the adjoint action of $G$ (resp., $W$).
That is to say,
if
$S_k (\mfg^\vee)$ (resp., $S_k (\mft^\vee)$) denotes the symmetric algebra on $\mfg$ (resp., $\mft$) over $k$, regarded as   its coordinate ring, then
 $\mfg /\!\!/G$ (resp., $\mft /\!\!/W$)  is defined as the spectrum of the  subring of polynomial invariants $S_k (\mfg^\vee)^G$ (resp., $S_k (\mft^\vee)^W$).
Let us use the notation $\mfc$ to denote the $k$-scheme $\mft/\!\!/W$.
This $k$-scheme carries  a $\mbG_m$-action arising  from the homotheties on $\mfg$.

A Chevalley's theorem asserts (cf. ~\cite[Theorem 1.1.1]{Ngo}, ~\cite[Chap.\,VI, Theorem 8.2]{KiWe}) that the $k$-morphism $\mfc \rightarrow \mfg/\!\!/G$ induced by the natural $k$-algebra morphism $S_k (\mfg^\vee)^G \rightarrow S_k (\mft^\vee)^W$ is an isomorphism because we have assumed
the inequality $2h_G < p$.
Also, in $S_k (\mfg^\vee)^G \left(\cong S_k (\mft^\vee)^W \right)$,
there are $\mr{rk}(\mfg)$ algebraically independent (relative to $k$) homogenous elements, which generates $S_k (\mfg^\vee)^G$.   
Denote by 
\begin{align}
\chi : \mfg \rightarrow \mfc
\end{align}
the morphism of $k$-schemes defined as the composite of the  quotient $\mfg \twoheadrightarrow \mfg/\!\!/G$ and the inverse of the isomorphism $\mfc \xrightarrow{\sim} \mfg/\!\!/G$ resulting from Chevalley's theorem.

Finally, we shall write $\mfg^{\mr{ad}(q_1)} := \{ x \in \mfg \, | \, [q_1, x] = 0\}$ and write
\begin{align}
\kappa : \mfg^{\mr{ad}(q_1)} \hookrightarrow \mfg \xrightarrow{\chi} \mfc
\end{align}
for  the composite of $\chi$ and the injection  $\mfg^{\mr{ad}(q_1)}\hookrightarrow \mfg$ given by $v \mapsto q_{-1} + v$ for any $v \in \mfg^{\mr{ad}(q_1)}$.
Since $2h_G < p$,
$\kappa$ is verified  to be an isomorphism of $k$-schemes (cf. ~\cite[Lemma 1.2.1]{Ngo}).
In particular,  the equalities   $\mr{dim}(\mfc) = \mr{rk}(\mfg) = \mr{dim}(\mfg^{\mr{ad}(q_1)})$ hold.

\LSP
\subsection{Representability of the moduli stack of opers} \label{SS25}

We shall set 
\begin{align} \label{Eq25}
\mcN := T^\vee \overline{\mcM}_{\mr{ell}} \setminus \mr{Im} (0_{T^\vee\overline{\mcM}_{\mr{ell}}}),
\end{align}
where $T^\vee \overline{\mcM}_{\mr{ell}}$ denotes the total space of the cotangent bundle of $\overline{\mcM}_{\mr{ell}}$ over $k$ and $0_{T^\vee\overline{\mcM}_{\mr{ell}}}$ denotes the zero section $\overline{\mcM}_{\mr{ell}} \rightarrow T^\vee\overline{\mcM}_{\mr{ell}}$.
This stack defines a $\mbG_m$-bundle on $\overline{\mcM}_\mr{ell}$  and classifies   pairs $(\msX, \delta)$ consisting of a stable elliptic curve $\msX := (X/S, \sigma)$ and 
a global generator $\delta$ of $\Omega \left(=\Omega_{X^\mr{log}/S^\mr{log}}\right)$, i.e., $\Omega = \mcO_X\delta$. 
For example, 
if $\msX$ is an elliptic curve given by the usual Weierstrass equation
\begin{align}
y^2 + a_1 xy + a_3 y^2 = x^3 + a_2 x^2 + a_4 x + a_6 \ \ (a_i \in k),
\end{align}
then the  {\it invariant differential} 
$\frac{dx}{2y + a_1 x + a_3} \left(= \frac{dy}{3x^2 + 2a_2 x + a_4 -a_1 y} \right)$ determines  an object in $\mcN$.

Let us take a $k$-scheme $S$ and  an object $(\msX, \delta)$ of $\mcN$ over $S$.
The dual $\delta^\vee$ of $\delta$ determines a generator of $\mcT$.
Then, for each $v \in \mfb (S)$,
 the  pair
\begin{align} \label{Eq89}
\msE^\spadesuit_{\delta, v} := (\mcE_{B, \mr{triv}}, \nabla_{\delta^\vee, q_{-1} + v})
\end{align}
 forms a $G$-oper, where the trivial $B$-bundle $\mcE_{B, \mr{triv}}$ is regarded as a $B$-reduction of $\mcE_{G, \mr{triv}}$.
The  assignment $(\msX, \delta, \rho) \mapsto (\msX, \delta, \msE^\spadesuit_{\delta, \kappa^{-1}(\rho)})$ for each $\rho \in \mfc (S)$ is functorial with respect to $S$, and hence determines  a functor
\begin{align} \label{Eq31}
\mcN \times_k \mfc \rightarrow \mcN \times_{\overline{\mcM}_\mr{ell}} \mcO p_G
\end{align}
over $\mcN$.
Regarding this functor, we can prove the following assertion.

\SSP
\bpr \label{Prop14}
The functor  \eqref{Eq31} is an equivalence of categories.
That is to say,
if we are given   a stable elliptic curve $\msX$ over $S \in \mr{ob}(\mcS ch_{/k})$  equipped with a choice of a global  generator $\delta$ of $\Omega$, then
 for each $G$-oper $\msE^\spadesuit$ on $\msX_\mr{ell}$ there exists a unique element
\begin{align} \label{Eq93}
\rho_{\delta, \msE^\spadesuit} \in \mfc (S)
\end{align}
that is mapped to $\msE^\spadesuit$ via \eqref{Eq31}.
\epr
\begin{proof}
We shall construct an inverse of \eqref{Eq31}.
Let $S$  be a $k$-scheme and $(\msX, \delta, \msE^\spadesuit)$ a collection of data classified by  the fiber of $\mcN \times_{\overline{\mcM}_\mr{ell}} \mcO p_G$ over $S$.
The associated log curve $X^\mr{log}/S^\mr{log}$ gives   a specific  $B$-bundle
 ${^\dagger}\mcE_B$ on $X$ constructed in the manner of ~\cite[Eq.\,(210)]{Wak5}.
 According to ~\cite[Proposition 2.19]{Wak5},
  there exists uniquely a pair $({^\dagger}\mcE^\spadesuit, \mr{nor}_{\msE^\spadesuit})$ consisting of a $G$-oper ${^\dagger}\mcE^\spadesuit := ({^\dagger}\mcE_B, {^\dagger}\nabla)$ on $\msX_\mr{ell}$ that is $\{ x_\alpha \}_\alpha$-normal (cf. ~\cite[Definition 2.14]{Wak5})  and an isomorphism of $G$-opers $\mr{nor}_{\msE^\spadesuit} : {^\dagger}\msE^\spadesuit \xrightarrow{\sim} \msE^\spadesuit$.
Since the pair $(X, \delta^\vee)$ defines a log chart on $X^\mr{log}/S^\mr{log}$ in the sense of ~\cite[Definition 1.42]{Wak5},
it gives a trivialization $\mcE_{B, \mr{triv}} \xrightarrow{\sim}{^\dagger}\mcE_B$ of ${^\dagger}\mcE_B$.
Under the identification $\mcE_{B, \mr{triv}} ={^\dagger}\mcE_B$ given by  this trivialization, ${^\dagger}\nabla$ can be expressed as ${^\dagger}\nabla = \nabla_{\delta^\vee, q_{-1} +v}$ for a unique $v \in \mfg^{\mr{ad}(q_1)} (S)$.
 Thus, we obtain a well-defined element $\rho_{\delta, \msE^\spadesuit} := \kappa (v)$ of $\mfc (S)$.
 The resulting assignment $\msE^\spadesuit \mapsto \rho_{\delta, \msE^\spadesuit}$ is functorial with respect to $S$, and hence determines  a functor $\mcN \times_{\overline{\mcM}_\mr{ell}} \mcO p_G \rightarrow \mcN \times_k \mfc$
over $\mcN$.
This functor  is verified to define  an inverse of \eqref{Eq31},  which means that \eqref{Eq31} is an equivalence of categories.
\end{proof}
\SSP

The following assertion, i.e., the elliptic version of ~\cite[Theorem 2.39]{Wak5},  is  a direct consequence of the above proposition using  the descent property of the fpqc morphism $\mcN \rightarrow \overline{\mcM}_\mr{ell}$.

\SSP
\bco \label{Cor9}
The category $\mcO p_G$ can be represented by a nonempty,  connected, and  smooth Deligne-Mumford stack over $k$ and each geometric fiber of  the projection $\mcO p_G \rightarrow \overline{\mcM}_\mr{ell}$ 
is isomorphic to the $k$-scheme $\mfc$.
\eco

\LSP
\subsection{Miura opers on stable elliptic curves} \label{SS1331}

Next, we  recall the notion of  a (generic) Miura oper on a log curve (cf. ~\cite[Definition 3.2.1]{Wak7}).
Let $S^\mr{log}$ be an fs log scheme over $k$ and $f^\mr{log} : X^\mr{log} \rightarrow S^\mr{log}$ a log curve over $S^\mr{log}$.

\bde[cf. ~\cite{Wak7}, Definition 3.2.1, for the logarithmic setting] \label{Def}
A {\bf Miura $G$-oper} on $X^\mr{log}/S^\mr{log}$ is a triple
\begin{align}
\widehat{\mcE}^\spadesuit := (\mcE_B, \mcE'_B, \nabla)
\end{align}
 such that 
 $(\mcE_B, \nabla)$
 is a $G$-oper  on $X^\mr{log}/S^\mr{log}$ and $\mcE'_B$ is another  $B$-reduction of $\mcE_G := \mcE_B \times^B G$   horizontal with respect to $\nabla$, i.e., 
 the $S^\mr{log}$-connection $\nabla : \mcT \rightarrow \widetilde{\mcT}_{\mcE_G^\mr{log}/S^\mr{log}}$
 (where $\mcE_G := \mcE_B \times^B G$ and $\mcT := \mcT_{X^\mr{log}/S^\mr{log}}$)
  factors through
 the injection $\widetilde{\mcT}_{{\mcE'}_{B}^{\mr{log}}} \hookrightarrow \widetilde{\mcT}_{\mcE_G^\mr{log}/S^\mr{log}}$ induced from the inclusion $\mcE'_B\hookrightarrow \mcE_G$.
 If $X^\mr{log}/S^\mr{log}$ arises from a stable elliptic curve $\msX$, then we refer to any Miura $G$-oper on $X^\mr{log}/S^\mr{log}$ as a {\bf Miura $G$-oper on $\msX_\mr{ell}$}.
 Moreover,
the definition of an isomorphism between Miura $G$-opers can be formulated in a natural fashion.
\ede

Let $\widehat{\msE}^\spadesuit := (\mcE_B, \mcE'_B, \nabla)$ be a Miura $G$-oper on  $X^\mr{log}/S^\mr{log}$.
By twisting the flag variety $G/B$ by the $B$-bundle  $\mcE_B$,
 we obtain a proper scheme 
$(G/B)_{\mcE_B} := \mcE_B \times^B (G/B)$ over $X$.
Here,  denote by $w_0$ the longest element of $W$.
The Bruhat decomposition $G = \coprod_{w \in W} B  w  B$ gives rise to
 a decomposition
\begin{align}
(G/B)_{\mcE_B} = \coprod_{w \in W} S_{\mcE_B, w}
\end{align}
of $(G/B)_{\mcE_B}$, where
each $S_{\mcE_B, w}$ denotes  the $\mcE_B$-twist of $B w_0 w B$, i.e., $S_{\mcE_B, w} := \mcE_B \times^B (B w_0 w B)$.
The $B$-reduction $\mcE'_B$ determines a section 
$\sigma_{\mcE_B, \mcE'_B} :  X \rightarrow (G/B)_{\mcE_B}$
  of the natural projection $(G/B)_{\mcE_B} \rightarrow X$.
Given an element $w$ of $W$ and a point $x$ of $X$, 
we  say that $\mcE_B$ and $\mcE'_B$ are in {\bf relative position $w$ at $x$} if $\sigma_{\mcE_B, \mcE'_B} (x)$ belongs to $S_{\mcE_B, w}$.
In particular, if $\sigma_{\mcE_B, \mcE'_B} (x)$ belongs to $S_{\mcE_B, 1}$, then we say that $\mcE_{B}$ and $\mcE'_B$ are in {\bf generic position at $x$}.

\SSP
\bde[cf. ~\cite{Wak7}, Definition 3.3.1] \label{Def67}
We shall say that $\widehat{\msE}^\spadesuit$ is {\bf generic} $\mcE_B$ and $\mcE'_B$ are in generic position at every  point of $X$. 
\ede
\SSP

Just as in the case of $\mcO p_G$, we obtain the category
\begin{align} \label{Eq101}
\mcM \mcO p_G
\end{align}
fibered in equivalence relations over $\mcS ch_{/k}$ classifying pairs $(\msX, \widehat{\msE}^\spadesuit)$ consisting a stable elliptic curve $\msX$ and a generic Miura $G$-oper $\widehat{\msE}^\spadesuit$ on $\msX_\mr{ell}$.
Forgetting the data of generic Miura $G$-opers yields a functor  $\mcM \mcO p_G \rightarrow \overline{\mcM}_\mr{ell}$, and 
the assignment $(\mcE_B, \mcE'_B, \nabla) \mapsto (\mcE_B, \nabla)$ determines a functor
\begin{align} \label{Eq102}
\mcM \mcO p_G \rightarrow \mcO p_G
\end{align}
over $\overline{\mcM}_\mr{ell}$.

Let us take a  $k$-scheme $S$ and an object $(\msX, \delta)$ of $\mcN$ over $S$.
Denote by $B^{-}$ the opposite Borel subgroup of $B$ relative to $T$, which determines a $B$-reduction $\mcE_{B^{-1}, \mr{triv}}$ of $\mcE_{G, \mr{triv}}$.
Then, the triple
\begin{align} \label{Eq105}
\widehat{\msE}_{\delta, v}^\spadesuit := (\mcE_{B, \mr{triv}}, \mcE_{B^{-1}, \mr{triv}}, \nabla_{\delta^\vee, q_{-1} + v})
\end{align}
defined for each $v \in \mft (S)$  forms a generic Miura $G$-oper, and 
the resulting assignment $(\msX, \delta, v) \mapsto (\msX, \delta, \widehat{\msE}^\spadesuit_{\delta, v})$ determines a functor
\begin{align} \label{Eq104}
\mcN \times_k \mft \rightarrow \mcN \times_{\overline{\mcM}_\mr{ell}} \mcM \mcO p_G
\end{align}
over $\mcN$, making the following square diagram commute:
\begin{align} \label{Eq114}
\vcenter{\xymatrix@C=46pt@R=36pt{
 \mcN \times_k \mft \ar[r]_-{}^-{\eqref{Eq104}} \ar[d]_-{\mr{id}_\mcN \times \chi |_\mft} & 
 \mcN \times_{\overline{\mcM}_\mr{ell}} \mcM \mcO p_G
 \ar[d]^-{\mr{id}_\mcN \times \eqref{Eq102}} \\
 \mcN \times_k \mfc  \ar[r]^-{\sim}_-{\eqref{Eq31}} & \mcN \times_{\overline{\mcM}_\mr{ell}}  \mcO p_G.
 }}
\end{align}

Similarly to Proposition \ref{Prop14}, we obtain the following assertion.

\SSP
\bpr \label{Prop136}
The functor \eqref{Eq104} is an equivalence of categories.
That is to say, 
if we are given a stable elliptic curve $\msX$ over $S \in \mr{ob} (\mcS ch_{/k})$ equipped with a choice of a global  generator $\delta$ of $\Omega$,
then for 
each generic Miura $G$-oper $\widehat{\msE}^\spadesuit$ on   $\msX_\mr{ell}$
 there exists a unique element
\begin{align} \label{Eq139}
\mu_{\delta, \widehat{\msE}^\spadesuit} \in \mft (S)
\end{align}
that is mapped to $\widehat{\msE}^\spadesuit$ via \eqref{Eq104}.
\epr
\begin{proof}
Let $\widehat{\msE}^\spadesuit := (\mcE_B, \mcE'_B, \nabla)$ be a generic Miura $G$-oper on $\msX_\mr{ell}$.
As discussed in the proof of Proposition \ref{Prop14},
the choice of $\delta$ allows us to  identify $\mcE_B$  and $\mcE_G \left(:= \mcE_B \times^B G \right)$ with $\mcE_{B, \mr{triv}}$ and $\mcE_{G, \mr{triv}}$, respectively.
Also, after possibly transposing $\mcE'_B$ via 
the left-translation $L_h$ by some $h \in B(S)$,
we may assume that $\mcE'_B = \mcE_{B^{-}, \mr{triv}}$ and $\nabla = \nabla_{\delta^\vee, q_{-1} + \mu_{\delta, \widehat{\msE}^\spadesuit}}$ for some $\mu_{\delta, \widehat{\msE}^\spadesuit} \in \mft (S)$ (cf. Proposition \ref{Prop15}).
Note that an automorphism of $\mcE_{G, \mr{triv}}$ preserving two $B$-reductions $\mcE_{G, \mr{triv}}$ and  $\mcE_{G, \mr{triv}}$ arises from the  left-translation by 
an $S$-rational point of $B \cap B^{-} = T$.
Also,   for  
two elements $v, v' \in \mft (S)$,
$\nabla_{\delta^\vee, q_{-1}+ v}$ coincides with 
$\nabla_{\delta^\vee, q_{-1}+ v}$ up to  a gauge transformation by an element of $T$  if and only if  the equality $v = v'$ holds.
Hence, the element $\mu_{\delta, \widehat{\msE}^\spadesuit}$ is uniquely determined and depends only on $\widehat{\msE}^\spadesuit$.
It follows that the assignment $\widehat{\msE}^\spadesuit \mapsto \mu_{\delta, \widehat{\msE}^\spadesuit}$ is well-defined and gives a functor  
$\mcN \times_{\overline{\mcM}_\mr{ell}} \mcM \mcO p_G \rightarrow \mcN \times_k \mft$.
This functor turns out to be an inverse of \eqref{Eq104}, so \eqref{Eq104} defines  an equivalence of categories.
\end{proof}
\SSP

The above proposition together with the descent property of $\mcN \rightarrow \overline{\mcM}_\mr{ell}$ implies the following assertion.

\SSP
\bco \label{Cor33}
The category $\mcM \mcO p_G$ can be represented by a nonempty, connected, and smooth Deligne-Mumford stack over $k$ and each geometric fiber of the projection $\mcM \mcO p_{G} \rightarrow \overline{\mcM}_\mr{ell}$ is isomorphic to the $k$-scheme $\mft$.
Moreover, the morphism $\mcM \mcO p_G \rightarrow \mcO p_G$ is finite, faithfully flat, and of degree $\sharp (W)$.
\eco

\LSP
\section{Moduli space of dormant opers on elliptic curves} \label{S44}

In this section, we review  dormant  (Miura) opers  and examine the geometric structure of 
their moduli stacks.
The proofs of the various facts asserted  in Theorems  \ref{Th4992}-\ref{ThF} are provided.
In particular, we establish
the connectedness of  the moduli stacks, meaning  that any two elliptic curves equipped with  dormant (generic Miura) opers
can be continuously deformed into each other.
The proof proceeds by explicitly describing such structures in both the ordinary and supersingular cases for the underlying curve.

For  previous studies  on dormant (Miura) opers on pointed stable curves in positive characteristic,
we refer the reader to ~\cite{Moc2}, ~\cite{Wak2}, ~\cite{Wak7}, ~\cite{Wak5},  ~\cite{Wak8}, ~\cite{Wak13}, and ~\cite{Wak9}, and related works.

\LSP
\subsection{$p$-curvature of a logarithmic connection} \label{SS7}

Let $S$ be a $k$-scheme and $\msX := (f : X \rightarrow S, \sigma)$ a stable elliptic curve over $S$;
for simplicity, we write $\mcT := \mcT_{X^\mr{log}/S^\mr{log}}$ and $\Omega := \Omega_{X^\mr{log}/S^\mr{log}}$.
Also,
let  us fix  a flat $G$-bundle $(\mcE, \nabla)$ on $X^\mr{log}/S^\mr{log}$.

Recall here that, for a morphism of fs log scheme $f : Y^\mr{log} \rightarrow T^\mr{log}$,
 the logarithmic tangent bundle $\mcT_{Y^\mr{log}/T^\mr{log}}$ has a structure of restricted $\mcO_T$-algebra, in the sense of ~\cite[Definition 3.4]{Wak5}.
In particular, it is equipped with a $p$-th power operation $(-)^{[p]} : \mcT_{Y^\mr{log}/T^\mr{log}} \rightarrow \mcT_{Y^\mr{log}/T^\mr{log}}$ (cf. ~\cite[Proposition 1.2.1]{Ogu1}).
If both the log structures of $Y^\mr{log}$ and $T^\mr{log}$ are  trivial, then 
the assignment $(-)^{[p]}$ can be realized as taking the $p$-th iterations of  locally defined derivations  $\mcO_Y \rightarrow \mcO_Y$ corresponding to sections in $\mcT_{Y/T}$.
Let us  apply the formation of this operation to the cases where $Y^\mr{log}/T^\mr{log}$ is taken to be $X^\mr{log}/S^\mr{log}$ and $\mcE^\mr{log}/S^\mr{log}$.
Hence,   we obtain $p$-th power operations on $\mcT$ and $\widetilde{\mcT}_{\mcE^\mr{log}/S^\mr{log}}$.
According to ~\cite[\S\,3.3]{Wak5},
there exists an $\mcO_X$-linear morphism
\begin{align} \label{Eq3990}
\psi (\nabla) : F^*_{X} (\mcT) \rightarrow \mfg_\mcE \left(\subseteq \widetilde{\mcT}_{\mcE^\mr{log}/S^\mr{log}}\right)
\end{align}
determined by the condition that $F^{-1}_X (\partial)$ is sent to $\nabla (\partial)^{[p]}-\nabla (\partial^{[p]})$
for any local section $\partial \in \mcT$.
We refer to $\psi (\nabla)$ as the {\bf $p$-curvature} of $\nabla$.
See ~\cite[\S\,1.2]{Ogu1} for the $p$-curvature of a logarithmic connection defined on an $\mcO_X$-module.

\LSP
\subsection{Hasse invariant} \label{SS35}

Recall that $\msX$ is called {\bf ordinary} (resp., {\bf supersingular}) if the  $\mcO_S$-linear morphism
\begin{align} \label{Eq90}
\varphi_{\msX} : F_S^* (\mbR^1 f_*(\mcO_X)) \rightarrow \mbR^1 f_* (\mcO_X)
\end{align}
induced by the Frobenius morphism $F_X$ is an isomorphism (resp., vanishes identically).
The locus 
\begin{align} \label{Eq91}
\overline{\mcM}_{\mr{ell}}^\mr{ord} \ \left(\text{resp.,} \ \overline{\mcM}_{\mr{ell}}^\mr{ss}  \right)
\end{align}
in $\overline{\mcM}_{\mr{ell}}$ classifying {\it ordinary} (resp., {\it supersinguler}) stable elliptic curves forms a dense  open (resp., a proper closed) substack.
In particular, $\overline{\mcM}_{\mr{ell}}^\mr{ord}$ contains $\partial \overline{\mcM}_\mr{ell}$.

Suppose that we are given 
a global generator $\delta$ of $\Omega$.
Since its dual $\delta^\vee$ generates $\mcT$,
there exists a unique element 
$H (\msX, \delta) \in H^0 (S, \mcO_S) \left(= H^0 (X, \mcO_X)\right)$ satisfying
\begin{align}
(\delta^\vee)^{[p]} = H (\msX, \delta) \cdot \delta^\vee.
\end{align}
We refer to this element  as  the {\bf Hasse invariant}  of $(\msX, \delta)$.
For each $\lambda \in H^0 (S, \mcO_S^\times)$, the equality
$H (\msX, \lambda \cdot \delta) = \lambda^{1-p} \cdot H (\msX, \delta)$ holds
 (cf. ~\cite[\S\,6.2.11]{Hid} for the non-logarithmic version).
 
 Note that $\msX$ is ordinary (resp., supersingular) if and only if $H (\msX, \delta)$ lies in $H^0 (S, \mcO^\times_S)$ (resp., $H (\msX, \delta)$ vanishes).
Indeed,
if $C : F_{X/S*} (\Omega) \rightarrow \Omega_{X^{(1)\mr{log}}/S^\mr{log}}$ denotes  the Cartier operator  (i.e., the map ``$C$" resulting from  ~\cite[Proposition 1.2.4]{Ogu1} in the case where ``$(E, \nabla)$" is taken to be the trivial one $(\mcO_X, d)$),
then it is well-known that
its direct image
\begin{align}
f^{(1)}_*(C) :  f_* (\Omega) \left(= f_*^{(1)} (F_{X/S*} (\Omega)) \right) \rightarrow \left(f^{(1)}_* (\Omega_{X^{(1)\mr{log}}/S^\mr{log}})  \cong \right)F_S^* (f_*(\Omega))
\end{align}
via $f^{(1)} : X^{(1)} \rightarrow S$ is dual to $\varphi_{\msX}$ via Grothendieck-Serre duality.
 On the other hand,
$f^{(1)}_* (C)$ is  dual to the $p$-th power operation $(-)^{[p]}$ because
\begin{align}
\langle  C (\delta), F^{-1}_S (\delta^\vee)\rangle = \langle \delta, (\delta^\vee)^{[p]} \rangle - (\delta^{\vee})^{p-1} (\langle \delta, \delta^\vee \rangle) = \langle \delta, (\delta^\vee)^{[p]} \rangle,
\end{align}
where the first equality follows from  ~\cite[Eq.\,(1.2.4.5)]{Ogu1}.
It follows that $(-)^{[p]}$ coincides with  $\varphi_{\msX}$ under the canonical identification $\mbR^1f_* (\mcO_X) = 
f_*(\mcT)$, which implies the requited equivalence.

For 
each element   $v \in \mfg (S)$,
the $p$-curvature of the $S^\mr{log}$-connection $\nabla_{\delta^\vee, v}$ on  $\mcE_{G, \mr{triv}}$ satisfies
\begin{align} \label{Eq100}
\psi (\nabla_{\delta^\vee, v}) (\delta^\vee) = \nabla (\delta^\vee)^{[p]} - \nabla ((\delta^\vee)^{[p]}) =  v^{[p]} - H (\msX, \delta) \cdot v \in \mfg (S).
\end{align}

\LSP
\subsection{Hitchin-Mochizuki morphism and $p$-nilpotent opers} \label{SS53}

In the rest of this  paper, suppose that $G$, $B$, and $T$  are all defined over $\mbF_p$ (hence $\mfg$, $\mfb$, $\mft$, and $\mfc$ can be  defined over $\mbF_p$).

Denote by $\A$ (resp., $\A^\nabla$) the set-valued sheaf on
$\overline{\mcM}_{\mr{ell}}$ which assigns, to any object $\msX \in  \mr{ob}(\overline{\mcM}_{\mr{ell}})$ (where $\msX := (X/S, \sigma)$) over a $k$-scheme $S$, the set 
of morphisms $\varsigma : X \rightarrow (F_X^* (\Omega))^\times \times^{\mbG_m} \mfc$ (resp., $X \rightarrow \Omega^\times \times^{\mbG_m} \mfc$) over $X$.
Here, for each line bundle $\mcL$, we denote by $\mcL^\times$ the $\mbG_m$-bundle corresponding to $\mcL$.
Both $\A$ and $\A^\nabla$ can  be represented by a relative scheme over $\overline{\mcM}_\mr{ell}$ of finite type.
Also, the pull-back via  $F_{X/S}$ induces an $\overline{\mcM}_\mr{ell}$-morphism $F^{\nabla \Rightarrow \emptyset} : \A^\nabla \rightarrow \A$, which is verified to be a closed immersion because the proof of  ~\cite[Proposition 3.23]{Wak5}  is available even in our non-hyperbolic situation.

Let $(\msX, \msE^\spadesuit)$, where $\msX := (X/S, \sigma)$ and $\msE^\spadesuit := (\mcE_B, \nabla)$, be an object of $\mcO p_G$ over a $k$-scheme $S$.
The pair $(\mcE_G, \psi (\nabla))$  specifies a global section of the quotient stack
$[((F_X^*(\Omega))^\times \times^{\mbG_m} \mfg)/G]$ over $X$.
The composite of this section with the morphism
\begin{align}
[((F_X^*(\Omega))^\times \times^{\mbG_m} \mfg)/G] \rightarrow (F_X^* (\Omega))^\times \times^{\mbG_m} \mfc
\end{align}
induced by $\chi : \mfg \rightarrow \mfc$ defines an $X$-rational point $\overline{\psi (\nabla)}$ of $(F_X^* (\Omega))^\times \times^{\mbG_m} \mfc$, or equivalently, an $S$-rational point of $\A$. 
The resulting assignment $\msE^\spadesuit \mapsto \overline{\psi (\nabla)}$ determines an $\overline{\mcM}_\mr{ell}$-morphism 
\begin{align}
\mr{HM} : \mcO p_G \rightarrow \A.
\end{align}
One can  apply 
 the proof of ~\cite[Proposition 3.2]{LasPa} (or ~\cite[Proposition-Definition 3.24]{Wak5}) to our situation, and we see that
  $\mr{HM}$ factors through $F^{\nabla \Rightarrow \emptyset}$.
Thus, $\mr{HM}$ restricts to  a morphism
\begin{align}
\mr{HM}^\nabla : \mcO p_G \rightarrow \A^\nabla
\end{align}
over $\overline{\mcM}_\mr{ell}$, referred to as the {\bf Hitchin-Mochizuki morphism}.

Next, let us take an object $(\msX, \delta, \varsigma)$ of $\mcN \times_{\overline{\mcM}_\mr{ell}} \A^\nabla$, where $\msX$ is as above.
The trivialization $\Omega = \mcO_X$ of $\Omega$ determined by $\delta$ specifies  an identification $\Omega^\times \times^{\mbG_m} \mfc = X \times_k \mfc$.
Under this identification, $\varsigma$ corresponds to an element $\varsigma' \in \mfc (X) \left(= \mfc (S) \right)$.
The assignment $(\msX, \delta, \varsigma) \mapsto (\msX, \delta, \varsigma')$ gives  
an isomorphism
\begin{align} \label{Eq77}
\mcN \times_{\overline{\mcM}_\mr{ell}}\A^\nabla \xrightarrow{\sim} \mcN \times_k \mfc
\end{align}
over $\mcN$.
Using this isomorphism,
we obtain a composite endomorphism
\begin{align}
\xi : \mcN \times_k \mfc \xrightarrow{\eqref{Eq31}}  \mcN \times_{\overline{\mcM}_\mr{ell}} \mcO p_G \xrightarrow{\mr{id}_\mcN \times \mr{HM}^\nabla} \mcN \times_{\overline{\mcM}_\mr{ell}} \A^\nabla \xrightarrow{\eqref{Eq77}} \mcN \times_k \mfc
\end{align}
 of $\mcN \times_k \mfc$ over $\mcN$.
According to \eqref{Eq100}, this morphism  sends $(\msX, \delta, \rho)$ to $(\msX, \delta, \gamma_{H (\msX, \delta)} (\rho))$, where 
for each $a \in H^0 (S, \mcO_S)$ and $\rho \in \mfc (S)$, we set
\begin{align} \label{Eq78}
\gamma_{a} (\rho) :=\chi ((q_{-1} +\kappa^{-1}(\rho))^{[p]} - a  \cdot (q_{-1}+\kappa^{-1}(\rho))).
\end{align}

\SSP
\bpr \label{Cor17}
The morphism $\xi$  is finite, faithfully flat, and  of degree $p^{\mr{rk}(\mfg)}$. 
Moreover, if we write $\mcN^{\mr{ss}} := \mcN \times_{\overline{\mcM}_\mr{ell}} \overline{\mcM}_\mr{ell}^\mr{ss}$, then
the restriction $\xi^\mr{ss} : \mcN^\mr{ss} \times_k \mfc \rightarrow \mcN^\mr{ss} \times_k \mfc$ of $\xi$ over $\mcN^\mr{ss} \left(\subseteq \mcN \right)$ coincides with $\mr{id}_{\mcN^\mr{ss}} \times F_{\mfc/k}$.
\epr
\begin{proof}
First, let us  prove the finiteness of $\xi$.
We shall set $R := S_k [\mft^\vee]^W \left(= S_k [\mfg^\vee]^G \right)$.
The $\mcO_\mcN$-algebra corresponding to $\mcN \times_k \mfc$ is given by $\mcO_\mcN \otimes_k R$.
The  $\mbG_m$-action on $\mfc$ opposite to the usual one 
 induces a $\mbZ_{\geq 0}$-grading 
$R = \bigoplus_{j \in \mbZ_{\geq 0}}  R_j$.
By letting $R^j := \bigoplus_{j' \leq j} R_{j'}$,
we obtain an increasing filtration $\{ \mcO_\mcN \otimes_k R^j \}_{j \in \mbZ_{\geq 0}}$ on $\mcO_\mcN \otimes_k R$.
It follows from \eqref{Eq78} that
the $\mcO_\mcN$-algebra endomorphism 
$\xi^* : \mcO_\mcN \otimes_k R \rightarrow \mcO_\mcN \otimes_k R$ corresponding to $\xi$ satisfies 
$\xi^* (\mcO_\mcN \otimes_k R^j) \subseteq \mcO_\mcN \otimes_k R^{p\cdot j}$ for every $j \in \mbZ_{\geq 0}$.
This endomorphism  induces, for each $j$, a morphism $\mr{gr}^j (\xi^*) : \mr{gr}^j (\mcO_\mcN \otimes_k R) \rightarrow \mr{gr}^{p \cdot j} (\mcO_\mcN \otimes_k R)$, 
where $\mr{gr}^{j'} (\mcO_\mcN \otimes_k R) := (\mcO_\mcN \otimes_k R^{j'})/(\mcO_\mcN \otimes_k R^{j'-1}) \left(= \mcO_\mcN \otimes_k R_{j'} \right)$.
The direct sum of $\mr{gr}^{j} (\xi^*)$'s defines an $\mcO_\mcN$-algebra endomorphism 
\begin{align}
\bigoplus_{j}\mr{gr}^j (\xi^*) : \mcO_\mcN \otimes_k R \rightarrow \mcO_\mcN \otimes_k R
\end{align}
of  $\mcO_\mcN \otimes_k R \left(= \bigoplus_{j \in \mbZ_{\geq 0}} \mcO_\mcN \otimes_k R_j \right)$.
The explicit description of $\xi$ displayed in  \eqref{Eq78}
 shows that 
$\bigoplus_{j}\mr{gr}^j (\xi^*)$ coincides with the $\mcO_\mcN$-algebra endomorphism $\xi_0^*$ corresponding to the endomorphism $\xi_0 : \mcN \times_k \mfc \rightarrow \mcN \times_k \mfc$ given by $(\msX, \delta, \rho) \mapsto (\msX, \delta, \gamma_0 (\rho))$.
 Here, observe  
 the following sequence of equalities   for $\rho \in \mfc (S)$:
\begin{align}  \label{Eq105}
\gamma_{0} (\rho) & =\chi ((q_{-1} +\kappa^{-1}(\rho))^{[p]} - 0  \cdot (q_{-1}+\kappa^{-1}(\rho))) \\
&=  \chi ((q_{-1} + \kappa^{-1} (\rho))^{[p]}) \notag \\
&= \chi (q_{-1}+ \kappa^{-1}(\rho))^{[p]} \notag \\
& = F_{\mfc/k} (\rho),
\end{align}
where the last equality follows from ~\cite[Proposition 3.12]{Wak5}.
This implies  $\xi_0 = \mr{id}_\mcN \times F_{\mfc/k}$, and 
$\mcO_\mcN \otimes_k R$ turns out to be  a finite 
$\mr{Im}(\xi_0^*) \left(= \mr{Im} (\bigoplus_{j}\mr{gr}^j (\xi^*))\right)$-module.
It follows that $\mcO_\mcN \otimes_k R$ is finitely generated as an  $\mr{Im}(\xi^*)$-module.
This concludes the finiteness of $\xi$.

Also, since  the equality $\gamma_{H (\msX, \delta)} = \gamma_0$ holds  for any $(\msX, \delta) \in \mr{ob}(\mcN^\mr{ss})$, 
the sequence \eqref{Eq105} also shows $\xi^\mr{ss} = \mr{id}_{\mcN^\mr{ss}} \times F_{\mfc/k}$, which completes the proof of the latter assertion.

Next, let us consider the surjectivity of $\xi$.
Since $\mcN \times_k \mfc$ is irreducible,
the image $\mr{Im} (\xi)$ of $\xi$ is  irreducible.
The equality $\xi^\mr{ss} = \mr{id}_{\mcN^\mr{ss}} \times F_{\mfc/k}$ just obtained  and the surjectivity of $F_{\mfc /k}$ together  imply that the image of $\xi$ contains $\mcN^\mr{ss} \times_k \mfc \left(\subseteq \mcN \times_k \mfc \right)$, which is of codimension
$1$.
Hence, $\mr{Im} (\xi)$ must be the entire space $\mcN \times_k \mfc$ because each fiber of the $1$-st projection $\mcN \times_k \mfc \rightarrow \mcN$   contains at least  one  point in   $\mr{Im} (\xi)$.
That is to say, $\xi$ is surjective.

The flatness of $\xi$ follows from ~\cite[Chap.\,III, Exercise 10.9]{Har} together with the fact that both the domain and codomain of $\xi$ are smooth over $k$ and each fiber of $\xi$ is a nonempty $0$-dimensional scheme.

Finally, the degree computation  of $\xi$ can be immediately verified since  $\xi^\mr{ss} \left( = \mr{id}_{\mcN^\mr{ss}} \times F_{\mfc/k}\right)$ is of degree $p^{\mr{rk} (\mfg)} \left(= p^{\mr{dim}(\mfc)} \right)$.
Thus, we have finished the proof of this proposition.
\end{proof}
\SSP

By the descent property of the fpqc morphism $\mcN \rightarrow \overline{\mcM}_\mr{ell}$, the above proposition implies the following assertion.

\SSP
\bco \label{Cor13}
The Hitchin-Mochizuki morphism $\mr{HM}^\nabla$ is finite, faithfully flat, and  of degree $p^{\mr{rk}(\mfg)}$. 
Moreover, 
the fiber of $\mr{HM}^\nabla$ over 
any  $k$-rational point $a \in \overline{\mcM}_\mr{ell}^\mr{ss}$ coincides with $F_{\mfc/k}$ under certain identifications $\mcO p_G |_a = \A^\nabla |_a = \mfc$ (cf. \eqref{Eq31} and \eqref{Eq77}).
\eco
\SSP

Denote by $0_{\A^\nabla}$ the global section $\overline{\mcM}_\mr{ell} \rightarrow \A^\nabla$ determined by the zero section $X \rightarrow \Omega^\times \times^{\mbG_m} \mfc$.
The inverse image of $0_{\A^\nabla}$ via $\mr{HM}^\nabla$
defines a closed substack
\begin{align}
\mcO p_G^{^{p\text{-}\mr{nilp}}} := (\mr{HM}^\nabla)^{-1} (0_{\A^\nabla})
\end{align}
of $\mcO p_G$.
This stack classifies, by definition, $G$-opers $\msE^\spadesuit := (\mcE_B, \nabla)$ with $\overline{\psi (\nabla)} = 0$; such $G$-opers are called  {\bf $p$-nilpotent} (cf. ~\cite[Definition 3.15]{Wak5}).

\SSP
\bt \label{Th13}
The stack $\mcO p_G^{^{p\text{-}\mr{nilp}}}$ is finite,  faithfully flat, and  of degree $p^{\mr{rk}(\mfg)}$ over $\overline{\mcM}_\mr{ell}$.
Moreover, the fiber of the projection $\mcO p_G^{^{p\text{-}\mr{nilp}}} \rightarrow \overline{\mcM}_\mr{ell}$ over any $k$-rational point $a \in \overline{\mcM}_\mr{ell}^\mr{ss}$ is isomorphic to $\mr{Spec} (k[\epsilon_1, \cdots, \epsilon_{\mr{rk}(\mfg)}]/(\epsilon^p_1, \cdots, \epsilon^p_{\mr{rk}(\mfg)}))$.
\et
\begin{proof}
The assertion follows from  Corollary \ref{Cor13}.
\end{proof} 

\LSP
\subsection{Dormant opers and the finiteness of the moduli stack} \label{SS53}

\SSP
\bde[cf, ~\cite{Wak7}, Definition 3.8.1, ~\cite{Wak5}, Definition 3.15] \label{Def4}
We shall say that a $G$-oper $\msE^\spadesuit:= (\mcE_B, \nabla)$ (resp., a  Miura $G$-oper $\widehat{\msE}^\spadesuit := (\mcE_B, \mcE'_B, \nabla)$)
is {\bf dormant} if $\psi (\nabla)$ vanishes identically.
\ede
\SSP

The collection of dormant $G$-opers (resp., dormant generic Miura $G$-opers) 
 form a closed substack
\begin{align} \label{eqwq2}
\mcO p_{G}^{^\mr{Zzz...}}  \ \left(\text{resp.,} \ \mcM \mcO p_G^{^\mr{Zzz...}} \right)
\end{align}
of $\mcO p_G$ (resp., $\mcM \mcO p_G$), which has  the natural projection
\begin{align} \label{Eq40}
\Pi_G : \mcO p_G^{^\mr{Zzz...}} \rightarrow \overline{\mcM}_\mr{ell} \ \left(\text{resp.,} \ \widehat{\Pi}_G :   \mcM \mcO p_G^{^\mr{Zzz...}} \rightarrow \overline{\mcM}_\mr{ell} \right).
\end{align}
Since 
any dormant $G$-oper is $p$-nilpotent, 
$\mcO p^{^\mr{Zzz...}}_G$ (resp., $\mcM \mcO p^{^\mr{Zzz...}}_G$) can also be considered as a closed substack of $\mcO p_G^{^{p\text{-}\mr{nilp}}}$ (resp., $\mcO p_G^{^{p\text{-}\mr{nilp}}} \times_{\mcO p_G} \mcM \mcO p_G$).
The morphism \eqref{Eq102} restricts to a morphism
\begin{align} \label{EQ127}
\Xi_G : \mcM \mcO p_G^{^\mr{Zzz...}} \rightarrow \mcO p_G^{^\mr{Zzz...}}
\end{align}
over $\overline{\mcM}_\mr{ell}$, which is finite, faithfully flat, and of  degree $\sharp (W)$ because of the commutativity of  the square diagram \eqref{Eq114}.
Hence, Theorem \ref{Th13} implies the following assertion.

\SSP
\bt \label{Th12}
The  stack $\mcO p^{^\mr{Zzz...}}_G$ (resp., $\mcM \mcO p^{^\mr{Zzz...}}_G$) is finite over $\overline{\mcM}_\mr{ell}$, and hence, proper over $k$.
Moreover, the fiber of $\Pi_G$ (resp., $\widehat{\Pi}_G$) over any $k$-rational point  $a \in \overline{\mcM}_\mr{ell}^\mr{ss}$ consists exactly of a single point.
\et
\begin{proof}
The remaining part is   to prove the second assertion.
First, we shall consider the non-resp'd portion.
Let $\msX := (X, \sigma)$ be the supersingular elliptic curve classified by $a$ and choose a global generator $\delta$  of $\Omega_{X^\mr{log}/k}$.
Then,
$\msE^\spadesuit_{\delta, 0} := (\mcE_{B, \mr{triv}}, \nabla_{\delta^\vee, q_{-1}})$ forms a $G$-oper. 
Since $p > \mr{rk}(\mfg)$ and 
$q_{-1} \in \mfg_{-1}$,
we have $\mr{ad} ((q_{-1})^{[p]}) = \mr{ad}(q_{-1})^p = 0$ (cf. ~\cite[Definition 3.4]{Wak5} for the first equality), which implies  $(q_{-1})^{[p]} = 0$.
It follows that 
\begin{align}
\psi (\nabla_{\delta^\vee, q_{-1}}) \stackrel{\eqref{Eq100}}{=} (q_{-1})^{[p]} - 0 \cdot q_{-1} = 0.
\end{align}
That is to say, $\msE^\spadesuit_{\delta, 0}$ is dormant and specifies a point $\widetilde{a}$ of  the fiber 
$\Pi_G^{-1}(a)$.
On the other hand,  Theorem \ref{Th13} says that $\Pi_G^{-1}(a)$ contains at most one point because $\mcO p_G^{^\mr{Zzz...}}$ is a substack of  $\mcO p_G^{^{p\text{-}\mr{nilp}}}$.
Consequently, the equality  $\Pi_G^{-1}(a) = \{ \widetilde{a} \}$ holds, as desired.

Next, let us consider
the resp'd assertion.
Observe that the set-theoretic preimage of the zero element of $\mfc$ via  $\chi |_\mft : \mft \rightarrow \mfc$  consists  exactly of a single point.
Hence, by the commutativity of \eqref{Eq114},
the fiber of $\Xi_G$ over the point corresponding to $(\msX, \delta, \msE^\spadesuit_{\delta, 0})$ is set-theoretically a single point.
It follows that  the assertion  follows  from the non-resp'd portion and the equality $\widehat{\Pi}_G =   \Pi_G \circ \Xi_G$.
\end{proof}

\LSP
\subsection{Frobenius-invariant differentials} \label{SS56}

Denote by $f_{\mr{ell}}^{\mr{ord}} : \mcC_{\mr{ell}}^\mr{ord} \rightarrow \overline{\mcM}_\mr{ell}^\mr{ord}$ the universal family of stable elliptic curves over $\overline{\mcM}_\mr{ell}^\mr{ord}$.
By definition, the $\mcO_{\overline{\mcM}_\mr{ell}^\mr{ord}}$-linear morphism
\begin{align}
\varphi : F^*_{\overline{\mcM}_\mr{ell}} (\mbR^1 f^\mr{ord}_{\mr{ell}*} (\mcO_{\mcC_{\mr{ell}}^\mr{ord} })) \rightarrow  \mbR^1 f^\mr{ord}_{\mr{ell}*} (\mcO_{\mcC_{\mr{ell}}^\mr{ord} })
\end{align}
induced from the Frobenius morphism $F_{\mcC_\mr{ell}^\mr{ord}}$
is an isomorphism.
In other words, the pair $(\mbR^1 f^\mr{ord}_{\mr{ell}*} (\mcO_{\mcC_{\mr{ell}}^\mr{ord} }), \varphi)$ specifies a unit $F$-crystal on $\overline{\mcM}_\mr{ell}^\mr{ord}$.
By the result of ~\cite[Proposition 4.1.1]{Kat3} (generalized to Deligne-Mumford stacks),
this unit $F$-crystal corresponds to a locally free \'{e}tale sheaf of $\mbF_p$-modules;
the total space of this sheaf minus the zero section will be denoted by
\begin{align}
\mcN_\mr{inv}.
\end{align}

Let $\msX := (f: X \rightarrow S, \sigma)$ be  an ordinary  stable elliptic curve  over a $k$-scheme $S$.
As mentioned in \S\,\ref{SS35}, the morphism $\varphi_\msX$ (cf. \eqref{Eq90}) coincides with the $p$-th power operation on $\mcT := \mcT_{X^\mr{log}/S^\mr{log}}$ under the canonical identification $\mbR^1 f_* (\mcO_X) = f_* (\mcT)$.
This implies that
each section  of $\mcN_\mr{inv}$ over the $S$-rational point of $\overline{\mcM}_\mr{ell}^\mr{ord}$ classifying $\msX$  can be regarded as  a global  generator  $\partial$ of $\mcT$
  with $\partial^{[p]} = \partial$ (i.e, $H (\msX, \partial^\vee) =1$).
In particular,  the assignment $\partial \mapsto \partial^\vee$ defines a closed immersion 
$\mcN_\mr{inv} \hookrightarrow \mcN$, by which we consider $\mcN_\mr{inv}$ as a closed substack of $\mcN$.

The opposite of the natural $\mbG_m$-action on $\mcN$ restricts, via the inclusion $\mbF^\times_p \hookrightarrow \mbG_m$,  to an $\mbF_p^\times$-action on $\mcN_\mr{inv}$, which is given by
$(c, \partial) \mapsto c \cdot \partial$ for $c \in \mbF_p^\times$.
By this action, 
the natural projection $\mcN_\mr{inv} \rightarrow \overline{\mcM}_\mr{ell}^\mr{ord}$ becomes  
an \'{e}tale Galois covering with Galois group $\mbF_p^\times$.

Now, recall that an element $v$ of $\mft$ is called {\bf regular} if
its stabilizer $G_v$ under the adjoint action coincides with $T$.
The element $v$ is regular precisely when
 $\alpha (t) \neq 0$ for all  root $\alpha$ (cf. ~\cite[Chap.\,VI, Theorem 7.2]{KiWe}), and such elements
 form a nonempty open subscheme $\mft_\mr{reg}$ of $\mft$.
Since $\mft_\mr{reg}$ can be defined over $\mbF_p$, it makes sense to speak of the set of $\mbF_p$-rational points $\mft_\mr{reg}(\mbF_p)$; it carries  a natural $W$-action, so we obtain the quotient  $\mft_\mr{reg}(\mbF_p)/W$ with respect to this action, which is a finite subset of $\mfc (k)$.

\SSP
\bpr \label{Prop19}
Let $\msX$ be as above and $\partial$ a global generator of $\mcT$ with $\partial^{[p]} = \partial$, i.e., $(\msX, \partial) \in \mr{ob} (\mcN_\mr{inv})$.
Also,
let $\msE^\spadesuit$ (resp., $\widehat{\msE}^\spadesuit$) be a $G$-oper (resp., a generic Miura $G$-oper) on $\msX_\mr{ell}$.
Then, 
$\msE^\spadesuit$ (resp., $\widehat{\msE}^\spadesuit$) is dormant if and only if
$\rho_{\partial^\vee, \msE^\spadesuit} \in \mft_\mr{reg}(\mbF_p) /W$ (resp., $\mu_{\partial^\vee, \widehat{\msE}^\spadesuit} \in \mft_\mr{reg}(\mbF_p)$).
\epr
\begin{proof}
Since each $G$-oper extends to a generic Miura $G$-oper after the base-change over an fpqc  covering of $S$ (because of the commutativity of \eqref{Eq114}),
it suffices to consider the resp'd assertion.
For simplicity, we write $\mu := \mu_{\partial^\vee, \widehat{\msE}^\spadesuit} \in \mft (S)$.

First, we shall prove the ``if" part of the required equivalence.
Suppose that $\mu \in \mft_{\mr{reg}}(\mbF_p)$.
By restricting $\widehat{\msE}^\spadesuit$ to every geometric fiber of $f : X \rightarrow S$,
one  can assume, without loss of generality, that
$S = \mr{Spec}(k)$. 
Since the morphism $G/T \times_k \mft_\mr{reg} \rightarrow \mfg_\mr{reg} := \chi^{-1}(\mft_\mr{reg}/W) \left(\subseteq \mfg \right)$ given by $(\overline{h}, w) \mapsto \mr{Ad}_G (h) (w)$ (where $h \in G$) is an \'{e}tale  Galois covering (cf. ~\cite[Chap.\,VI, Theorem 9.1]{KiWe}),
it follows from Proposition \ref{Prop15} that $\nabla_{\partial, q_{-1}+ \mu}$ is transposed into the $S^\mr{log}$-connection $\nabla_{\partial, \mu}$ via the gauge transformation by some element of $G (k)$.
The $p$-curvature $\psi (\nabla_{\partial, \mu})$ of $\nabla_{\partial, \mu}$  satisfies
\begin{align}
\psi (\nabla_{\partial, \mu}) (\partial) \stackrel{\eqref{Eq100}}{=} \mu^{[p]} - H (\msX, \partial^\vee) \cdot \mu = F_{\mft/k}(\mu) - \mu = 0,
\end{align}
where the second equality follows from the fact that the $p$-th power operation on $\mft$ coincides with $F_{\mft/k} : \mft \rightarrow \mft$ (cf.  ~\cite[Example 4.4.10, (2)]{Spr}), and the last equality follows from $\mu \in \mft (\mbF_p)$.
Hence, $\nabla_{\partial, q_{-1}+ \mu}$ has vanishing $p$-curvature, i.e., $\widehat{\msE}^\spadesuit$ is dormant.

Next, we shall prove the ``only if" part.
Suppose that $\widehat{\msE}^\spadesuit$ is dormant.
Let us take an arbitrary  $k$-rational point $a$ of $S$, and denote by  $\overline{\mu}$ the restriction of $\mu$ to the fiber of $\msX$ over $a$.
The element $q_{-1} + \overline{\mu} \in \mfg$ admits a Jordan decomposition $q_{-1} + \overline{\mu} = \mu_s + \mu_n$ with $\mu_s$ semisimple and $\mu_n$ nilpotent.
Denote by $\mr{ad} :   \mfg \rightarrow \mr{End} (\mfg)$ the adjoint representation of $\mfg$, which is injective and compatible with the respective restricted structures, i.e., $p$-th power operations.
One can find an isomorphism of restricted Lie algebras $\alpha: \mr{End} (\mfg) \rightarrow \mfg \mfl_{\mr{dim}(\mfg)}$ which sends $\alpha (\mr{ad} (\mu_s + \mu_n)) \left(= \alpha (\mr{ad} (\mu_s)) + \alpha (\mr{ad}(\mu_n)) \right)$ to a Jordan normal form.
Here, $\alpha (\mr{ad}(\mu_s))$ is diagonal and every entry of $\alpha (\mr{ad}(\mu_n))$ except the superdiagonal is $0$.
Let us observe the following sequence of equalities:
\begin{align} \label{Eq80}
\alpha (\mr{ad} (\mu_s)) + \alpha (\mr{ad} (\mu_n)) & = \alpha (\mr{ad}(\mu_s + \mu_n)) = \alpha (\mr{ad} ((\mu_s + \mu_n)^{[p]})) \\
& = \alpha (\mr{ad}(\mu_s + \mu_n))^p = (\alpha (\mr{ad}(\mu_s)) + \alpha (\mr{ad}(\mu_n)))^p,
\end{align}
where the second equality follows from both  \eqref{Eq100} and the dormancy assumption  on $\widehat{\msE}^\spadesuit$.
An explicit computation of $(\alpha (\mr{ad}(\mu_s)) + \alpha (\mr{ad}(\mu_n)))^p$ shows (by \eqref{Eq80}) that $\alpha (\mr{ad}(\mu_n)) =0$ (i.e., $\mu_n =0$),
so $q_{-1} + \overline{\mu}$ is conjugate to an element $v \in \mft$.
Also,   $q_{-1} + \overline{\mu}$
is regular (in the sense of  the comment preceding ~\cite[Lemma 1.2.3]{Ngo}) because
it  is conjugate to
$q_{-1} + \rho_{\partial^\vee, \msE^\spadesuit}$ with $\delta_{\partial^\vee, \msE^\spadesuit} \in \mfg^{\mr{ad}(q_1)}$, where $\msE^\spadesuit$ denotes the underlying $G$-oper of $\widehat{\msE}^\spadesuit$ (cf. Propositions \ref{Prop15} and \ref{Prop14}).
Hence,   $v$ belongs to $\mft_\mr{reg}$, which implies  $q_{-1} + \overline{\mu} \in \mfg_\mr{reg}$.
By considering this fact for every $a \in S$,  we have $q_{-1} + \mu \in \mfg_{\mr{reg}} (S)$.
Since the natural morphism $G/T \times_k \mft_{\mr{reg}} \rightarrow \mfg_\mr{reg}$ is an \'{e}tale Galois covering as recalled  above, 
$q_{-1} + \mu$ is conjugate to  $\mu$ after possibly replacing $S$ with its \'{e}tale covering.
On the other hand,
 the $p$-th power operation  $(-)^{[p]}$ on $\mft$ coincides with $F_{\mft/k}$, so the equality $\mu = \mu^{[p]}$ obtained from 
 \eqref{Eq100} and the dormancy condition on $\widehat{\msE}^\spadesuit$ 
 implies 
$\mu \in \mft (\mbF_p) \cap \mft_\mr{reg} (k) = \mft_\mr{reg} (\mbF_p)$.
This completes the proof of this proposition.
\end{proof}

\LSP
\subsection{Global structure of the module stacks} \label{SS58}

We occasionally regard a  finite set $A$ as the disjoint union of copies of $\mr{Spec}(k)$ indexed by the elements of $A$.
In particular,  the product with $\mcN_\mr{inv}$ yields 
\begin{align} \label{Eq131}
\mcN_\mr{inv} \times A := 
 \coprod_{v \in A} \mcN_{\mr{inv}, v},
\end{align}
where  $\mcN_{\mr{inv}, v}:= \mcN_{\mr{inv}}$.
If, moreover, $A$ is equipped with an $\mbF_p^\times$-action, then it induces an $\mbF_p^\times$-action on
$\mcN_\mr{inv} \times A$ given by $(c, (\msX, \partial, v)) \mapsto (\msX, c\cdot \partial, c \cdot v)$ for $c \in \mbF_p^\times$. 

We here consider the case where  $A = \mft_\mr{reg}(\mbF_p)/W$ (resp., $A = \mft_\mr{reg}(\mbF_p)$) equipped with the $\mbF_p^\times$-action arising from the homotheties on $\mft$.
A pair   of an element $v$ of $\mft_\mr{reg}(\mbF_p)/W \subseteq \mfc (k)$
(resp., $\mft_\mr{reg}(\mbF_p) \subseteq \mft (k)$) and an object 
$(\msX, \partial) \in \mcN_{\mr{inv}, v}$ specify
a $G$-oper $\msE^\spadesuit_{\partial^\vee, v}$ (resp., a generic Miura $G$-oper $\widehat{\msE}^\spadesuit_{\partial^\vee, v}$) on $\msX$.
It follows from Proposition \ref{Prop19} that
the resulting assignment $(\msX, \partial) \mapsto (\msX, \partial,  \msE^\spadesuit_{\partial^\vee, v})$ (resp., $(\msX, \partial) \mapsto (\msX, \partial,  \widehat{\msE}^\spadesuit_{\partial^\vee, v})$) determines an isomorphism  of stacks
\begin{align} \label{Eq29}
\Theta_G : \mcN_\mr{inv} \times (\mft_\mr{reg} (\mbF_p)/W)\xrightarrow{\sim} \mcN_\mr{inv} \times_{\overline{\mcM}_\mr{ell}} \mcO p_G^{^\mr{Zzz...}} \\
  \left(\text{resp.,} \  \widehat{\Theta}_G : \mcN_\mr{inv} \times \mft_\mr{reg} (\mbF_p)\xrightarrow{\sim} \mcN_\mr{inv} \times_{\overline{\mcM}_\mr{ell}} \mcM \mcO p_G^{^\mr{Zzz...}} \right) \notag \\
\end{align}
over $\mcN_\mr{inv}$, fitting into the following square diagram obtained by restricting
 \eqref{Eq114}:
\begin{align} \label{Eq123}
\vcenter{\xymatrix@C=46pt@R=36pt{
 \mcN_\mr{inv} \times \mft_\mr{reg} (\mbF_p) \ar[r]_-{}^-{\widehat{\Theta}_G} \ar[d]_-{\mr{id}_{\mcN_\mr{inv}} \times \chi |_{\mft_\mr{reg} (\mbF_p)}} & 
 \mcN_\mr{inv} \times_{\overline{\mcM}_\mr{ell}} \mcM \mcO p_G^{^\mr{Zzz...}}
 \ar[d]^-{\mr{id}_\mcN \times \eqref{Eq102}} \\
 \mcN_\mr{inv} \times (\mft_\mr{reg} (\mbF_p)/W)  \ar[r]^-{}_-{\Theta_G} & \mcN_\mr{inv} \times_{\overline{\mcM}_\mr{ell}}  \mcO p_G^{^\mr{Zzz...}}.
 }}
\end{align}
The domains of $\Theta_G$ and $\widehat{\Theta}_G$ have  $\mbF_p^\times$-actions as  discussed above and, in  the codomains, $\mbF_p^\times$ acts only on the respective first factors.
Then, all the morphisms in \eqref{Eq123} preserve the $\mbF_p^\times$-actions.
By using this diagram with $\mbF_p^\times$-actions, we can prove  the following assertion.

\SSP
\bt \label{Th7}
\begin{itemize}
\item[(i)]
The projection $\Pi_G : \mcO p_{G}^{^\mr{Zzz...}} \!\rightarrow \overline{\mcM}_{\mr{ell}}$
restricts to a surjective, finite, and   \'{e}tale morphism
\begin{align} \label{Eq152}
\Pi_G^\mr{ord} : \overline{\mcM}_\mr{ell}^\mr{ord} \times_{\overline{\mcM}_\mr{ell}}\mcO p_{G}^{^\mr{Zzz...}} \!\rightarrow \overline{\mcM}_\mr{ell}^\mr{ord},
\end{align}
and its 
 degree coincides with  $\sharp (\mft_\mr{reg}(\mbF_p)/W) \left(=\sharp (\mft_\mr{reg}(\mbF_p))/\sharp (W)\right)$.
\item[(ii)]
The restriction
\begin{align}
\Xi^\mr{ord}_G  : \overline{\mcM}^\mr{ord}_\mr{ell} \times_{\overline{\mcM}_\mr{ell}}
\mcM \mcO p_G^{^\mr{Zzz...}} \rightarrow \overline{\mcM}^\mr{ord}_\mr{ell} \times_{\overline{\mcM}_\mr{ell}}\mcO p_G^{^\mr{Zzz...}}
\end{align}
of  the projection $\Xi_G : \mcM \mcO p_G^{^\mr{Zzz...}} \rightarrow \mcO p_G^{^\mr{Zzz...}}$  over $\overline{\mcM}^\mr{ord}_\mr{ell}$ is  an \'{e}tale Galois covering 
with Galois group $W$.
\end{itemize}
\et
\begin{proof}
Since  the lower horizontal arrow in \eqref{Eq123} is an isomorphism,
assertion (i) follows from the descent property of the \'{e}tale surjection $\mcN_\mr{inv} \twoheadrightarrow \overline{\mcM}_\mr{ell}^\mr{ord}$.

Next, let us consider assertion (ii).
Note that the natural $W$-action on the upper left-hand  corner  in  \eqref{Eq123} (i.e., the $W$-action coming from that on $\mft_\mr{reg}(\mbF_p)$)  commutes  with the $\mbF_p^\times$-action on it.
It follows that
this  $W$-action
  descends, via the composite Galois covering
  \begin{align}
 \mcN_\mr{inv} \times \mft_\mr{reg} (\mbF_p) \xrightarrow{\widehat{\Theta}_G} \mcN_\mr{inv} \times_{\overline{\mcM}_\mr{ell}} \mcM \mcO p_G^{^\mr{Zzz...}} \twoheadrightarrow \overline{\mcM}_\mr{ell}^\mr{ord} \times_{\overline{\mcM}_\mr{ell}} \mcM \mcO p_G^{^\mr{Zzz...}},
  \end{align}
    to a $W$-action on  $\overline{\mcM}_\mr{ell}^\mr{ord} \times_{\overline{\mcM}_\mr{ell}} \mcM \mcO p_G^{^\mr{Zzz...}}$ with $(\overline{\mcM}_\mr{ell}^\mr{ord} \times_{\overline{\mcM}_\mr{ell}} \mcM \mcO p_G^{^\mr{Zzz...}})/W \cong \overline{\mcM}_\mr{ell}^\mr{ord} \times_{\overline{\mcM}_\mr{ell}} \mcO p_G^{^\mr{Zzz...}}$.
    This completes the proof of assertion (ii).
  \end{proof}
\SSP

\bco \label{Cor22}
Both $\mcO p_G^{^\mr{Zzz...}}$ and $\mcM\mcO p_G^{^\mr{Zzz...}}$ are connected.
Similarly, both $\mcM_\mr{ell} \times_{\overline{\mcM}_\mr{ell}}\mcO p_G^{^\mr{Zzz...}}$ and $\mcM_\mr{ell} \times_{\overline{\mcM}_\mr{ell}} \mcM\mcO p_G^{^\mr{Zzz...}}$ are connected.
\eco
\begin{proof}
By similarity of arguments, we only consider the first assertion.
Moreover,  
since $\Xi_G$ is surjective (cf. the commutativity of \eqref{Eq114}), it suffices to consider the connectedness of $\mcM\mcO p_G^{^\mr{Zzz...}}$.
By the irreducibility of $\overline{\mcM}_\mr{ell}$ together with Theorem \ref{Th7}, (i) and (ii), 
there exists  a connected component $\mcQ$ of $\mcM\mcO p_G^{^\mr{Zzz...}}$ that dominate $\overline{\mcM}_{\mr{ell}}$.
Since  $\widehat{\Pi}_G$  is finite, the projection $\mcQ \rightarrow \overline{\mcM}_{\mr{ell}}$ is surjective.
In particular, the unique point in the fiber of $\widehat{\Pi}_G$ over every point of $\overline{\mcM}_\mr{ell}^{\mr{ss}}$  (cf. Theorem \ref{Th12})  lies in $\mcQ$.
Now, suppose 
 that there exists another connected component $\mcQ' \subseteq \mcM\mcO p_G^{^\mr{Zzz...}}$.
 The component  $\mcQ$ contains $\widehat{\Pi}_G^{-1} (\overline{\mcM}_\mr{ell}^\mr{ss})$, so
 the image of $\mcQ'$ via $\widehat{\Pi}_G$ is contained in $\overline{\mcM}_\mr{ell}^\mr{ord}$ (i.e., $\mcQ' \subseteq \widehat{\Pi}_G^{-1}(\overline{\mcM}_\mr{ell}^\mr{ord})$).
 Since the restriction of $\widehat{\Pi}_G$ over $\overline{\mcM}_\mr{ell}^\mr{ord}$ is \'{e}tale and $\overline{\mcM}_\mr{ell}^\mr{ord}$ is dense in $\overline{\mcM}_{\mr{ell}}$, $\mcQ'$ dominates $\overline{\mcM}_{\mr{ell}}$.
It follows that $\widehat{\Pi}_G |_{\mcQ'} : \mcQ' \rightarrow \overline{\mcM}_{\mr{ell}}$ is surjective because of its finiteness.
However, it contradicts the fact that $\mcQ' \cap \widehat{\Pi}_G^{-1} (\overline{\mcM}_\mr{ell}^{\mr{ss}}) = \emptyset$.
This completes the proof of the assertion.
\end{proof}

\SSP
\begin{rem} \label{Rem5}
We have so far essentially been discussing the classification of dormant (generic Miura) $G$-opers on unpointed genus-$1$ curves. Then, what about the remaining case where the Euler characteristic is zero, i.e., the case of a $2$-pointed  projective line?

The  curve under consideration is now
 the projective line $\mbP := \mr{Proj}(k[x, y])$  over $k$, which  contains 
 the $k$-rational points $[0]$, $[1]$, and  $[\infty]$  determined by the values $0$, $1$,  and $\infty$, respectively.
 In particular, $\msP := (\mbP, \{ [0], [1], [\infty]\})$ is a unique (up to isomorphism) $3$-pointed curve of genus $0$.
 We equip  $\mbP$ with the log structure determined by the divisor $[0] + [\infty]$, and 
   denote the resulting log curve by $\mbP^\mr{log}/k$.
  The section 
 $\partial := s \cdot \frac{d}{ds} \left(= t \cdot \frac{d}{dt} \right)$, where $s := x/y$, $t := y/x$,  specifies  a  global generator of $\mcT_{\mbP^\mr{log}/k}$ and satisfies  $\partial^{[p]} = \partial$.
 
 We shall  refer to any  $G$-oper (resp., Miura $G$-oper)  on $\mbP^\mr{log}/k$ as a {\bf $G$-oper on $\msP_\mr{ell}$} (resp., a {\bf Miura $G$-oper on $\msP_\mr{ell}$}).
  The set of isomorphism classes of  dormant $G$-opers (resp., dormant generic Miura $G$-opers) on $\msP_\mr{ell}$ is denoted by 
 \begin{align}
 \mcO p^{^\mr{Zzz...}}_{G, \msP} \ \left(\text{resp.,} \  \mcM\mcO p^{^\mr{Zzz...}}_{G, \msP}\right). 
\end{align}
 Just as in the previous discussion (cf. the proof of Proposition \ref{Prop19}), 
 the collection of data
  \begin{align}
  \msE^\spadesuit_{\partial^\vee, \kappa^{-1}(\rho), \msP} := (\mcE_{B, \mr{triv}}, \nabla_{\partial, q_{-1}+ \kappa^{-1}(\rho)}) \ \left(\text{resp.,} \  \widehat{\msE}^\spadesuit_{\partial^\vee, v, \msP} := (\mcE_{B, \mr{triv}}, \mcE_{B^{-}, \mr{triv}}, \nabla_{\partial, q_{-1}+ v}) \right)
  \end{align}
  associated to  each  $\rho \in \mft_\mr{reg} (\mbF_p)/W \subseteq \mfc (k)$ 
  (resp., $v \in \mft_\mr{reg}(\mbF_p) \subseteq \mft (k)$)
  turns out to define a dormant $G$-oper (resp., a dormant generic Miura $G$-oper) on $\msP_\mr{ell}$.
  The resulting assignment $\rho \mapsto \msE^\spadesuit_{\partial^\vee, \kappa^{-1}(\rho), \msP}$ (resp., $v \mapsto \widehat{\msE}^\spadesuit_{\partial^\vee, v, \msP}$) gives a bijection
 \begin{align} \label{Eq144}
 \mft_\mr{reg} (\mbF_p)/W \xrightarrow{\sim}  \mcO p^{^\mr{Zzz...}}_{G, \msP} \ \left(\text{resp.,} \ 
  \mft_\mr{reg} (\mbF_p) \xrightarrow{\sim} \mcM\mcO p^{^\mr{Zzz...}}_{G, \msP}
  \right).
 \end{align}
The  projection $\Xi_{G, \msP} : \mcM\mcO p^{^\mr{Zzz...}}_{G, \msP} \rightarrow \mcO p^{^\mr{Zzz...}}_{G, \msP}$ fits into  the  commutative square diagram 
\begin{align} \label{Eq166}
\vcenter{\xymatrix@C=46pt@R=36pt{
 \mft_\mr{reg} (\mbF_p) \ar[r]^-{\eqref{Eq144}}_-{\sim} \ar[d]_-{\chi |_{\mft_\mr{reg}(\mbF_p)}} &   \mcM\mcO p^{^\mr{Zzz...}}_{G, \msP}\ar[d]^-{\Xi_{G, \msP}} \\
 \mft_\mr{reg} (\mbF_p)/W \ar[r]_-{\eqref{Eq144}}^-{\sim} & \mcO p^{^\mr{Zzz...}}_{G, \msP},
 }}
\end{align}
which  is regarded as an analogue of \eqref{Eq123}.
This allows us to describe dormant (generic Miura) $G$-opers in a simple fashion.

Moreover, note that the stable elliptic curve $\msX := (X/k, \sigma)$ classified by  the unique $k$-rational point $b$ in $\partial \overline{\mcM}_\mr{ell}$
  can be constructed  from $\msP$ by attaching $[0]$ and $[\infty]$ to form a  nodal point.
For any (logarithmic) $k$-connection $\nabla$ on the trivial $G$-bundle $\msE_{G, \mr{triv}}$ on $\mbP$,
its residues at the marked points satisfies  $\mr{Res}_{[0]} (\nabla) = - \mr{Res}_{[\infty]} (\nabla)$.
Hence, 
each dormant $G$-oper (resp., dormant generic Miura $G$-oper) on $\msP_\mr{ell}$ induces that on $\msX_\mr{ell}$ by gluing the fibers over
the points $[0]$ and $[\infty]$ (cf. ~\cite[\S\,7.2]{Wak5}).
This construction  yields  a bijection between $\mcO p^{^\mr{Zzz...}}_{G, \msP}$ (resp., $\mcM \mcO p^{^\mr{Zzz...}}_{G, \msP}$) and  the fiber of $\Pi_G$ (resp., $\widehat{\Pi}_G$) over $b$ (cf. ~\cite[Proposition 2.6.1]{Wak6}).
 In this way, one can identify \eqref{Eq166}  with the fiber of \eqref{Eq123} over the point $(\msX, \partial)$ in $\mcN_\mr{inv}$. 
\end{rem}

\LSP
\section{Canonical diagonal liftings for  ordinary  elliptic curves} \label{S545}
\LSP

This section extends the study of dormant $\mr{PGL}_n$-opers on elliptic curves in two directions.
The first extension considers  the case of  {\it prime-power characteristic}, while  the second explores  a {\it higher-level} generalization using  the sheaf $\mcD^{(m)}$ of differential operators of level $m \in \mbZ_{\geq 0}$  introduced in ~\cite{PBer1} and ~\cite{PBer2}.
In particular, 
dormant $\mr{PGL}_n$-opers in prime-power characteristic and those of  higher level (along with their  corresponding Miura oper variants) are explicitly classified by using  extensions of  \eqref{Eq144} and   \eqref{Eq166};
 see Propositions \ref{T39} and \ref{T34}.
As a consequence of this classification,  we prove   that these extensions are, in fact, equivalent.
To be more precise, they can be transformed into each other  via  {\it diagonal reduction/lifting} in the sense of ~\cite{Wak12} (cf. Theorem \ref{T47}).

\LSP
\subsection{$\mcD$-module structures  of finite level} \label{SS58}

Throughout this section, we fix 
an element  $\N$ of $\mbZ_{> 0}$, and let  $(\LL, \MM)$ be either $(\N, 1)$ or $(1, \N)$.
For each perfect field $k'$ in characteristic $p$, we denote by $W_{\LL} (k')$ the ring of Witt vectors of length $\LL$ over $k'$.

Next, let us consider  a smooth curve 
$f_{\LL} : X_{\LL} \rightarrow \mr{Spec}(W_{\LL})$ over $W_{\LL} := W_{\LL}(k)$, and denote 
 $\Omega_{\LL} := \Omega_{X_{\LL}/W_{\LL}}$ and $\mcT_\LL := \mcT_{X_\LL/W_\LL}$ for simplicity.
 
Following  ~\cite[\S\,2.2]{PBer1},
the sheaf of differential operators  $\mcD^{(\MM -1)} := \mcD_{X_{\LL}/W_{\LL}}^{(\MM -1)}$   on $X_{\LL}/W_{\LL}$ of level $\MM -1$ is defined, where 
$\mr{Spec}(W_{\LL})$ is equipped with an $(\MM -1)$-PD structure extending  to $X_{\LL}$ via $f_{\LL}$.
For each nonnegative integer $\LL'$ with  $\LL' \leq \LL$, we employ a subscript ``$\LL'$" to denote the reductions of objects over $W_{\LL}$ modulo $p^{\LL'}$ (e.g., $X_{\LL'} := W_{\LL'} \times_{W_\LL} X_\LL$).
We also denote by $\iota_{\LL'}$ the natural closed immersion $X_{\LL'} \hookrightarrow X_{\LL}$.

For each $j \in \mbZ_{\geq 0} \sqcup \{ \infty \}$,
we write $\mcD_{\leq j}^{(\MM -1)}$ for the subsheaf  of $\mcD^{(\MM -1)}$ consisting of differential operators of order $\leq j$, where $\mcD_{\leq \infty}^{(\MM -1)} := \mcD^{(\MM -1)}$.
Hence, we have $\mcD^{(\MM -1)} = \bigcup_{j \in \mbZ_{\geq 0}} \mcD_{\leq j}^{(\MM -1)}$.
Note that $\mcD_{\leq j}^{(\MM -1)}$ admits two different structures of $\mcO_{X_{\LL}}$-module, i.e., 
one as given by left multiplication, where we denote this $\mcO_{X_{\LL}}$-module by ${^L}\mcD_{\leq j}^{(\MM -1)}$, and the other given by right multiplication, where we denote this $\mcO_{X_{\LL}}$-module by ${^R}\mcD_{\leq j}^{(\MM -1)}$.
Given an $\mcO_{X_{\LL}}$-module $\mcF$,
we equip the tensor product $\mcD_{\leq j}^{(\MM -1)} \otimes \mcF := {^R}\mcD_{\leq j}^{(\MM -1)} \otimes \mcF$ with the $\mcO_{X_{\LL}}$-module structure given by left multiplication.

A {\bf (left) $\mcD^{(\MM -1)}$-module structure} on  $\mcF$ is a left $\mcD^{(\MM -1)}$-action $\nabla : {^L}\mcD^{(\MM -1)}  \left(:= {^L}\mcD_{\leq \infty}^{(\MM -1)} \right)\rightarrow \mcE nd_{W_{\LL}} (\mcF)$ on $\mcF$ extending its $\mcO_{X_{\LL}}$-module structure.
A {\bf  $\mcD^{(\MM -1)}$-bundle} is  a pair $(\mcF, \nabla)$ consisting of a vector bundle  $\mcF$ on $X_\LL$ and a $\mcD^{(\MM -1)}$-module structure $\nabla$ on $\mcF$.
A {\bf line  $\mcD^{(\MM -1)}$-bundle} is defined to be a $\mcD^{(\MM -1)}$-bundle  $(\mcL, \nabla)$ such that $\mcL$ is  a line bundle.
For instance, the structure sheaf $\mcO_{X_{\LL}}$ together with a natural left $\mcD^{(\MM -1)}$-action, denoted by $\nabla^{(\MM -1)}_\mr{triv}$,  specifies a trivial line $\mcD^{(\MM -1)}$-bundle.
For  two $\mcD^{(\MM -1)}$-bundles $(\mcF_\circ, \nabla_\circ)$, $(\mcF_\bullet, \nabla_\bullet)$,
their tensor product $\mcF_\circ \otimes \mcF_\bullet$ naturally  inherits 
 a $\mcD^{(\MM -1)}$-module structure,  denoted by  $\nabla_\circ \otimes \nabla_\bullet$.
The notion of a morphism between  $\mcD^{(\MM -1)}$-bundles can be defined in a natural manner.

The following assertion concerning line $\mcD^{(\MM -1)}$-bundles will be used in the proofs of Propositions \ref{T39} and \ref{T34}.

\SSP
\bpr \label{Prop449}
Let $\mcL$ be a line bundle on $X_{\LL}$ and for each $\Box \in \{ \circ, \bullet \}$, let $\nabla_\Box$  be a $\mcD^{(\MM -1)}$-module structures on $\mcL$.
Then, 
there exists a surjective morphism  of  $\mcD^{(\MM -1)}$-modules 
$(\mcL, \nabla_\circ) \rightarrow  (\mcL, \nabla_\bullet)$ if and only if $\nabla_\circ = \nabla_\bullet$.
\epr
\begin{proof}
The assertion follows from the fact that 
 any $\mcO_{X_\LL}$-linear  endomorphism of $\mcL$ must be the multiplication by an element of   $W_{\LL}$.
\end{proof}
\SSP

\begin{rem} \label{Rem442}
Suppose that $(\LL, \MM) = (\N, 1)$.
A {\bf connection} on an $\mcO_{X_{\N}}$-module $\mcF$ is a $W_{\N}$-linear morphism $\nabla : \mcF \rightarrow \Omega_{\N} \otimes \mcF$ satisfying the Leibniz rule, i.e., $\nabla (a \cdot v) = d a \otimes v + a \cdot \nabla (v)$ for any local sections $a \in \mcO_{X_{\N}}$, $v \in \mcF$.
It is well-known that giving a $\mcD^{(0)}$-module structure on $\mcF$ amounts to giving a connection on it (cf. ~\cite[Theorem 4.8]{BeOg}).
In the subsequent discussion, we will not distinguish between these two structures on $\mcF$.
\end{rem}
\SSP

Now, consider the case where  $(\LL, \MM) = (1, \N)$.
Let $\nabla$
be   a $\mcD^{(\N -1)}$-module structure of an $\mcO_{X_1}$-module $\mcF$.
As discussed in ~\cite[\S\,2.5]{Wak12},   the  
short exact sequence  $0 \rightarrow \mcD^{(\N -1)}_{\leq p^{\N}-1} \rightarrow \mcD^{(\N -1)}_{\leq p^\N} \rightarrow \mcT^{\otimes p^\N}_1 \rightarrow 0$ admits a canonical split injection $\mcT^{\otimes p^\N}_1 \hookrightarrow \mcD^{(\N -1)}_{\leq p^\N}$.
Using this, we define the {\bf $p^{\N}$-curvature} of $\nabla$ (or, of $(\mcF, \nabla)$) as 
 the  composite
\begin{align} \label{Eq1030}
\psi^\N (\nabla) :  \mcT^{\otimes p^\N}_1 \left(= F_X^{\N *}(\mcT_1) \right) \rightarrow {^L}\mcD^{(\N -1)} \xrightarrow{\nabla} \mcE nd_{k} (\mcF).
\end{align}
For a slightly different description of higher-level $p$-curvature, see ~\cite[Definition 3.1.1]{LSQ}.
When $\N =1$, this  coincides with the classical definition of $p$-curvature.
In particular, the $p$-curvature  $\psi^1 (\nabla) : \mcT_1^{\otimes p} \rightarrow \mcE nd_k (\mcF)$ of a connection $\nabla$ is determined  by  $\partial^{\otimes p} \mapsto (\nabla_\partial)^p -\nabla_{\partial^p}$ for any local section $\partial \in \mcT_1$, where
we define $\nabla_\partial := (\partial \otimes \mr{id}_\mcF) \circ \nabla$, and 
 $\partial^p$ denotes the section of $\mcT_1$ corresponding to the $p$-th iterate $\partial \circ \cdots \circ \partial$ of $\partial$.
If $\mcF$ is a vector bundle of rank $n \in \mbZ_{>0}$, then its $p$-curvature can be identified with the $p$-curvature of the corresponding $\mr{GL}_n$-bundle (cf. \eqref{Eq3990}). 

Additionally,
let $\mcS ol (\nabla)$ denote the subsheaf of $\mcF$ on which $\mcD_+^{(\N -1)}$ acts as zero via $\nabla$, where $\mcD_+^{(\N -1)}$ denotes the kernel of the canonical projection $\mcD^{(\N -1)} \twoheadrightarrow \mcO_{X_1}$.
 This sheaf can be regarded as an {\it $\mcO_{X_1}$-submodule} of $F_{X_1 *}^{\N}(\mcF)$ through  the underlying homeomorphism of $F_{X_1}^{\N}$.

On the other hand, 
 for an $\mcO_{X}$-module $\mcG$, there exists a canonical $\mcD^{(N-1)}$-module structure 
\begin{align} \label{E445}
\nabla_{\mcG, \mr{can}}^{(N-1)} : {^L}\mcD^{(N-1)} \rightarrow \mcE nd_k (F_{X_1}^{N*}(\mcG))
\end{align}
on the pull-back $F_{X_1}^{N*}(\mcG)$ with vanishing $p^\N$-curvature (cf. ~\cite[Definition 3.1.1 and Corollary 3.2.4]{LSQ}).
The  assignments $(\mcF, \nabla) \mapsto \mcS ol (\nabla)$ and $\mcG \mapsto (F_{X_1}^{\N *} (\mcG), \nabla_{\mcG, \mr{can}}^{(\N -1)})$ establishes  an equivalence of categories
\begin{align} \label{EQ1039}
\left(\begin{matrix}\text{the category of $\mcD^{(\N -1)}$-bundles} \\ \text{with vanishing $p^\N$-curvature}\end{matrix} \right) \xrightarrow{\sim} 
\left(\begin{matrix}\text{the category of} \\ \text{vector bundles on $X_1$}\end{matrix} \right).
\end{align}

We now return to the  case of a general $(\LL, \MM)$.

\SSP
\bde[cf. ~\cite{Wak12}, Definition 3.4] \label{Def500}
Let $(\mcF, \nabla)$ be a $\mcD^{(\MM -1)}$-bundle.
\begin{itemize}
\item[(i)]
Suppose that $(\LL, \MM) = (1, \N)$.
Then, we  say that $(\mcF, \nabla)$ (or, $\nabla$) is {\bf dormant} if it  has vanishing $p^\N$-curvature.
\item[(ii)]
Suppose that $(\LL, \MM) = (\N, 1)$ and $\mcF$ has rank $n \in \mbZ_{\geq 0}$.
Then, we say that $(\mcF, \nabla)$ (or, $\nabla$) is {\bf dormant} if it is, Zariski locally on $X_{\N}$, isomorphic to $(\mcO_{X_{\N}}, \nabla_\mr{triv}^{(0)})^{\oplus n}$.
\end{itemize}
\ede
\SSP

\begin{rem} \label{Rem739}
\begin{itemize}
\item[(i)]
In the case of $(\LL, \MM) = (\N, 1)$, 
the dormancy condition introduced in ~\cite{Wak12} was defined  for logarithmic  connections and    differs from the definition given  above.
However, as shown in  ~\cite[Proposition 3.1]{Wak12}, these two definitions are actually  equivalent in   the non-logarithmic setting.
Moreover,  when $\N =1$, 
  Cartier's classical theorem (cf. ~\cite[Theorem 5.1]{Kat1}) implies that the dormancy conditions defined in (i) and (ii) are equivalent.
\item[(ii)]
If $(\mcF_\circ, \nabla_\circ)$ and $(\mcF_\bullet, \nabla_\bullet)$ are both dormant $\mcD^{(\MM -1)}$-bundles,
then their tensor product $\nabla_\circ \otimes \nabla_\bullet$ is also dormant.
Additionally, 
for a dormant $\mcD^{(\MM -1)}$-bundle $(\mcF, \nabla)$,
 its   dual $(\mcF^\vee, \nabla^\vee)$ also defines  a dormant $\mcD^{(\MM -1)}$-bundle.
\end{itemize}
\end{rem}

\LSP
\subsection{Dormant $\mcD$-module structures on the structure sheaf} \label{SS166}
We denote by
\begin{align}
\mr{Conn}_{\LL, \MM}^{^\mr{Zzz...}}
\end{align}
 the set of dormant $\mcD^{(\MM -1)}$-module structures  on $\mcO_{X_{\LL}}$.

This set forms  an abelian group  under the operation of taking tensor products, 
where  the trivial line $\mcD^{(\MM -1)}$-bundle $(\mcO_{X_\LL}, \nabla^{(\MM -1)}_\mr{triv})$ serves as  the identity element (cf. Remark \ref{Rem739}, (ii)).
On the other hand, let $J_{X_1}$  denote the Jacobian variety of the mod $p$ reduction $X_1$ of $X_\LL$, which classifies line bundles of degree $0$ on $X_1$.
We write  $J_{X_1} [p^\N]$ for  the kernel of the endomorphism of $J_{X_1}$ given by
  $\mcL \mapsto \mcL^{\otimes p^\N}$ for each line bundle $\mcL$.
  In particular, we obtain an abelian group $J_{X_1} [p^\N] (k)$ consisting of its $k$-rational points.
  In what follows, we establish a correspondence between $\mr{Conn}_{\LL, \MM}^{^\mr{Zzz...}}$ and $J_{X_1} [p^\N] (k)$.

First, the case of  $(\LL, \MM) = (1, \N )$ follows directly from  \eqref{EQ1039}.
Indeed, take an element of $J_{X_1} [p^\N] (k)$, represented by a line bundle $\mcL$ on $X_1$ such that there exists  an isomorphism $s : \mcL^{\otimes p^\N} \left(= F_{X_1}^{\N *}(\mcL) \right)\xrightarrow{\sim}\mcO_{X_1}$.
Then, via $s$, $\nabla_{\mcL, \mr{can}}^{(\N -1)}$ can be regarded as  
a $\mcD^{(\N -1)}$-module structure  on $\mcO_{X_1}$.
This structure  does not depend on the choice of $s$, as   $\mr{End}_{\mcO_{X_1}} (\mcL) = \mr{End}_{\mcO_{X_1}} (\mcO_{X_1}) \cong k$.
By the equivalence of categories \eqref{EQ1039},
the resulting assignment $\mcL \mapsto \nabla_{\mcL, \mr{can}}^{(\N -1)}$ defines an isomorphism of  abelian groups
\begin{align} \label{Eq1011}
J_{X_1} [p^\N] (k) \xrightarrow{\sim}\mr{Conn}_{1, \N}^{^\mr{Zzz...}}.
\end{align}

Next, suppose that $(\LL, \MM) = (\N, 1)$, and  take a connection  $\nabla$  on $\mcO_U$ for a nonempty open subscheme $U$ of $X_\N$.
This connection  can be expressed as $\nabla = d + c (\nabla)$ for a unique $c (\nabla) \in H^0 (U, \Omega_\N)$, where $d$ denotes the universal derivation.
Note that  $c (\nabla \otimes \nabla') = c (\nabla) + c (\nabla')$ and $c (\nabla_\mr{triv}^{(0)}) = 0$.
Moreover, performing a gauge transformation by  multiplication with  an element 
$b \in H^0 (U, \mcO_U^\times)$,   the connection $\nabla$ is transformed  into 
$d + b^{-1} db + c (\nabla)$.
Thus,  $\nabla$ is  trivialized by a gauge transformation if and only if
$c (\nabla) = b^{-1} db$ for some $b \in H^0 (U, \mcO_U)$, i.e., $c (\nabla)$ belongs to the image of $H^0 (U, d \mr{log}) : H^0 (U, \mcO^\times_{X_\N}) \rightarrow H^0 (U, \Omega_{\N})$ induced from  $d \mr{log} : \mcO_{X_\N}^\times \rightarrow \Omega_\N$
 given by $a \mapsto a^{-1} da$.
If $\mcB^\mr{log}_\N$ denotes the image of $d \mr{log}$, then 
the assignment $\nabla \mapsto c (\nabla)$ yields  an isomorphism of abelian groups
\begin{align} \label{Eq1014}
\mr{Conn}_{\N, 1}^{^\mr{Zzz...}} \xrightarrow{\sim} H^0 (X_\N, \mcB^\mr{log}_\N).
\end{align}
We shall set
$\mcQ_\N := \mr{Ker} (d \mr{log}) \left(\subseteq \mcO_{X_\N}^\times \right)$.
The natural short exact sequence
$0 \rightarrow \mcQ_\N \rightarrow \mcO_{X_\N}^\times \rightarrow \mcB^\mr{log}_\N \rightarrow 0$ induces the following long exact sequence of cohomology groups:
\begin{align}
H^0 (X_\N, \mcQ_\N) \rightarrow H^0 (X_\N, \mcO^\times_{X_\N}) \rightarrow H^0 (X_\N, \mcB^\mr{log}_\N)
\rightarrow
H^1 (X_\N, \mcQ_\N) \xrightarrow{\tau} H^1 (X_\N, \mcO^\times_{X_\N}) \rightarrow  \cdots
\end{align}
Since $H^0 (X_\N, \mcQ_\N) =  H^0 (X_\N, \mcO_{X_\N}^\times) = W_\N^\times$,
this sequence yields an isomorphism
\begin{align} \label{Eq1050}
H^0 (X_\N, \mcB^\mr{log}_\N) \xrightarrow{\sim} \mr{Ker} (\tau).
\end{align}
On the other hand, by observing  the local description of $d \mr{log}$, one can see  that the composite $\mcQ_\N \hookrightarrow \mcO^\times_{X_\N} \xrightarrow{\iota^*_1} \mcO_{X_1}^\times$ factors through $F^{\N*}_{X} : \mcO_{X_1}^\times \rightarrow \mcO_{X_1}^\times$. 
(We here identify the underlying topological space of $X_\N$ with that of $X_1$ via $\iota_1$.)
If  $\eta$ denotes the resulting morphism $\mcQ_\N \rightarrow \mcO_{X_1}^\times$, then the associated morphism $H^1 (\eta) : H^1 (X_\N, \mcQ_\N) \rightarrow H^1 (X_1, \mcO_{X_1}^\times)$ makes the following square diagram commute:
\begin{align} \label{Eq1040}
\vcenter{\xymatrix@C=46pt@R=36pt{
H^1 (X_\N, \mcQ_\N) \ar[r]^-{\tau} \ar[d]_-{H^1 (\eta)} & H^1 (X_\N, \mcO_{X_\N}^\times) \ar[d]^-{H^1 (\iota^*_{1})} \\
H^1 (X_1, \mcO_{X_1}^\times) \ar[r]_-{H^1 (F^{\N *}_{X_1})} & H^1 (X_1, \mcO_{X_1}^\times). 
 }}
\end{align}
This induces  a morphism between the kernels of the horizontal arrows
\begin{align} \label{Eq1053}
\mr{Ker}(\tau) \rightarrow J_{X_1} [p^\N] (k) \left(= \mr{Ker} (H^1 (F^{\N *}_{X_1})) \right).
\end{align}

 \SSP
\bpr \label{Prop468}
Suppose that $X_1$ is ordinary in the usual sense.
Then, the morphism \eqref{Eq1053} is  an isomorphism.
In particular,  the inverse of 
the composition of  \eqref{Eq1014}, \eqref{Eq1050},  and  \eqref{Eq1053}  specifies  an isomorphism of abelian groups
\begin{align} \label{Eq1060}
J_{X_1} [p^\N] (k)  \xrightarrow{\sim} \mr{Conn}_{\N, 1}^{^\mr{Zzz...}}.
\end{align}
\epr
\begin{proof}
Let $M$ be an arbitrary positive integer satisfying  $M \leq N$.
Denote by  $\mcQ'_{M}$  the image of the mod $p^M$ reduction  $\mcQ_M \rightarrow \mcO_{X_M}^\times$  of the inclusion $\mcQ_\N \hookrightarrow \mcO_{X_\N}^\times$.
The inclusion  $\mcQ'_M \hookrightarrow \mcO_{X_M}^\times$ induces  a morphism between $1$-st cohomology groups
\begin{align}
\tau_M : H^1 (X_M, \mcQ'_M) \rightarrow H^1 (X_M, \mcO^\times_{X_M}).
\end{align}
For  $M \geq 2$, 
the short exact sequence
$0 \rightarrow \mcO_{X_1} \rightarrow \mcO^\times_{X_M} \rightarrow \mcO^{\times}_{X_{M-1}} \rightarrow 0$, where the second arrow is given by $a \mapsto 1 + p^{M-1} \cdot a$, yields  a short exact sequence
\begin{align} \label{Eq1012}
0 \rightarrow H^1 (X_1, \mcO_{X_1}) \rightarrow H^1 (X_{M}, \mcO^\times_{X_M}) \rightarrow H^1 (X_{M-1}, \mcO_{X_{M-1}}^\times) \rightarrow 0.
\end{align}
Similarly, the morphism $\mcO_{X_1} \rightarrow \mcQ'_M$ given by $a \mapsto 1 + p^{M-1} a^{p^{\N-M+1}}$ fits into the short exact sequence $0 \rightarrow \mcO_{X_1} \rightarrow \mcQ'_M \rightarrow \mcQ'_{M-1} \rightarrow 0$,
where the third arrow arises from 
the mod $p^{M-1}$ reduction $\mcO_{X_{M}}^\times \twoheadrightarrow \mcO_{X_{M-1}}^\times$.
This induces another  short exact sequence
\begin{align} \label{Eq1070}
0 \rightarrow H^1 (X, \mcO_{X}) \rightarrow H^1 (X_{M}, \mcQ'_M) \rightarrow H^1 (X_{M-1}, \mcQ'_{M-1}) \rightarrow 0.
\end{align}
The short exact sequences \eqref{Eq1012} and \eqref{Eq1070}  fit into the following commutative diagram:
\begin{align} \label{Eq203}
\vcenter{\xymatrix@C=26pt@R=36pt{
0 \ar[r] &  H^1 (X, \mcO_{X}) \ar[r] \ar[d]^-{H^1 (F_X^{\N -M +1})} & H^1 (X_M, \mcQ'_M) \ar[r] \ar[d]^-{\tau_M} &H^1 (X_{M-1}, \mcQ'_{M-1}) \ar[r] \ar[d]^-{\tau_{M-1}} & 0 \\
0 \ar[r] & H^1 (X, \mcO_{X_{}}) \ar[r] &
 H^1 (X_M, \mcO^\times_{X_M}) 
  \ar[r] & H^1 (X_{M-1}, \mcO^\times_{X_{M-1}})\ar[r] & 0.
 }}
\end{align}
Since $X$ is  assumed to be ordinary, the left-hand vertical arrow is an isomorphism.
Applying  the snake lemma to this diagram,  we deduce that
the morphism $H^1 (X_M, \mcQ'_M) \rightarrow H^1 (X_{M-1}, \mcQ'_{M-1})$ restricts to an isomorphism $\mr{Ker}(\tau_M) \xrightarrow{\sim} \mr{Ker} (\tau_{M-1})$.
In particular, we obtain a composite isomorphism
\begin{align} \label{er38}
\mr{Ker}(\tau_\N) \xrightarrow{\sim} \mr{Ker}(\tau_{\N -1}) \xrightarrow{\sim} \cdots \xrightarrow{\sim} \mr{Ker}(\tau_1).
\end{align}
It   coincides with \eqref{Eq1053} under the identifications 
 $\mr{Ker}(\tau) = \mr{Ker}(\tau_\N)$ and  $J_{X_1} [p^\N](k) \left(= \mr{Ker} (H^1 (F_{X_1}^{\N *})) \right)$ $= \mr{Ker} (\tau_1)$.
 Thus, we conclude that  \eqref{Eq1053} is bijective, thereby completing  the proof of this proposition.
\end{proof}
\SSP

\begin{rem} \label{Rem912}
Recall that the sequence $0 \rightarrow \mcO_{X_1}^\times \xrightarrow{F_{X_1}^*} \mcO_{X_1}^\times \xrightarrow{d \mr{log}} \Omega_1 \xrightarrow{C-\mr{id}} \Omega_{1}$ is exact, where $C$ denotes the Cartier morphism, considered as an endomorphism of $\Omega$ via the identification $\Omega_{X_1^{(1)}/k} = \Omega$ given by the base-change along $F_{\mr{Spec}(k)}$.
Hence, $H^0 (X_1, \mcB^\mr{log}_1)$ coincides with the subspace of $H^0 (X_1, \Omega)$ consisting of elements $\delta$ with $C (\delta) = \delta$.
(If $X_1$ has genus $1$, then it follows from the discussion in \S\,\ref{SS35} that the condition  $C (\delta) = \delta$ is equivalent to requiring  that the dual $\delta^\vee$ satisfies $(\delta^\vee)^{[p]} = \delta^\vee$.)

Now, suppose that $X_1$ is ordinary.
The square diagram 
\begin{align} \label{Eq1082}
\vcenter{\xymatrix@C=46pt@R=36pt{
H^0 (X_\N, \mcB^\mr{log}_\N) \ar[d]_-{\mr{mod}\,p\,\mr{reduction}} \ar[r]^-{\eqref{Eq1053}\circ \eqref{Eq1050}} & J_{X_1} [p^\N] (k) \ar[d]^-{\mcL \mapsto \mcL^{\otimes p^{\N-1}}} \\
H^0 (X_1, \mcB_1^\mr{log}) \ar[r]_-{\eqref{Eq1053}\circ \eqref{Eq1050}} &  J_{X_1}[p] (k)
 }}
\end{align}
is commutative and  the right-hand vertical arrow $J_{X_1} [p^\N] (k) \mapsto J_{X_1} [p] (k)$ is surjective.
Both the upper and lower horizontal arrows in this square are isomorphisms by  Proposition \ref{Prop468}.
It follows that  the mod $p$ reduction
\begin{align}
H^0 (X_\N, \mcB_\N^\mr{log}) \rightarrow H^0 (X_1, \mcB_1^\mr{log})
\end{align}
is surjective.
In other words, every  global section $\delta \in H^0 (X_1, \Omega)$ with $C (\delta) = \delta$ can be lifted to an element of $H^0 (X_\N, \mcB_\N^\mr{log})$.
\end{rem}

\LSP
\subsection{$\mr{GL}_n^{(\MM)}$-opers and generic Miura $\mr{GL}_n^{(\MM)}$-opers} \label{SS123}

We fix an integer $n$ with $1 < n < p$.
Consider a collection of data
\begin{align}
\msF^\heartsuit := (\mcF, \nabla, \{ \mcF^j \}_{j=0}^n)
\end{align}
where
\begin{itemize}
\item
$\mcF$ is a vector bundle on $X_{\LL}$ of rank $n$;
\item
$\nabla$ is a $\mcD^{(\MM -1)}$-module structure on $\mcF$;
\item 
$\{ \mcF^j \}_{j=0}^n$ is an $n$-step decreasing filtration
\begin{align}
0 = \mcF^n \subseteq \mcF^{n-1} \subseteq \cdots \subseteq \mcF^0 = \mcF
\end{align}
of $\mcF$ consisting of subbundles such that the subquotients $\mcF^j/\mcF^{j+1}$ are line bundles.
\end{itemize}

\SSP
\bde[cf. ~\cite{Wak12}, Definition 5.6] \label{Def55}
\begin{itemize}
\item[(i)]
We shall say that $\msF^\heartsuit$ is a {\bf $\mr{GL}_n^{(\MM)}$-oper} (or a {\bf $\mr{GL}_n$-oper of level $\MM$}) on $X_{\LL}$ if, for every $j=0, \cdots, n-1$, the $\mcO_{X_\LL}$-linear morphism
$\mcD^{(\MM -1)} \otimes \mcF \rightarrow \mcF$ induced  by  $\nabla$ restricts to an isomorphism
\begin{align} \label{Eq1100}
\mcD^{(\MM -1)}_{\leq n-j-1} \otimes \mcF^{n-1} \xrightarrow{\sim} \mcF^j.
\end{align}
\item[(ii)]
Let $\msF_\circ^\heartsuit := (\mcF_\circ, \nabla_\circ, \{ \mcF^j_\circ \}_{j})$ and $\msF_\bullet^\heartsuit := (\mcF_\bullet, \nabla_\bullet, \{ \mcF^j_\bullet \}_j)$ be $\mr{GL}_n^{(\MM)}$-opers on $X_\LL$.
An {\bf isomorphism of $\mr{GL}_n^{(\MM)}$-opers} from $\msF^\heartsuit_\circ$ to $\msF^\heartsuit_\bullet$ is an isomorphism of $\mcD^{(\MM -1)}$-bundles $(\mcF_\circ, \nabla_\circ) \xrightarrow{\sim} (\mcF_\bullet, \nabla_\bullet)$ 
compatible with the respective filtrations $\{ \mcF_\circ^j \}_j$ and $\{ \mcF_\bullet^j \}_j$.
\end{itemize}
\ede 
\SSP

Next, suppose further that we are given another $n$-step decreasing filtration $\{ \mcF^{(j)} \}_{j = 0}^n$ on $\mcF$
such that, for each $j =0, \cdots, n$, the composite $\alpha_j : \mcF^{(j)} \hookrightarrow \mcF \twoheadrightarrow \mcF/\mcF^{n-j}$ is an isomorphism, or equivalently, the composite $\beta_j : \mcF^j \hookrightarrow \mcF \twoheadrightarrow \mcF/\mcF^{(n-j)}$ is an isomorphism.
The composite
\begin{align} \label{Eq1000}
\mcF^{(j-1)} \xrightarrow{\alpha_{j-1}} \mcF/\mcF^{n-j+1} \twoheadrightarrow \mcF/\mcF^{n-j} \xrightarrow{\alpha_j^{-1}} \mcF^{(j)}
\end{align}
(where $j \in \{1, \cdots, n\}$) determines a split surjection of the short exact sequence $0 \rightarrow \mcF^{(j)} \rightarrow \mcF^{(j-1)} \rightarrow \mcF^{(j-1)}/\mcF^{(j)} \rightarrow 0$, which gives rise to a decomposition
\begin{align}
\gamma_j : \mcF^{(j-1)} \xrightarrow{\sim} \mcF^{(j)} \oplus \mcF^{(j-1)}/\mcF^{(j)}.
\end{align}
This implies that the kernel of the surjection $\mcF/\mcF^{n-j+1} \twoheadrightarrow \mcF/\mcF^{n-j}$ (appearing in the second arrow of \eqref{Eq1000})
is isomorphic to $\mcF^{(j-1)}/\mcF^{(j)}$.
That is, we have a canonical isomorphism
\begin{align}
\delta_{j} : \mcF^{n-j}/\mcF^{n-j+1} \xrightarrow{\sim} \mcF^{(j-1)}/\mcF^{(j)}.
\end{align}
Furthermore,   $\mcF$ decomposes into the direct sum of $n$ line bundles $\{ \mcF^{(j)}/\mcF^{(j+1)} \}_{j=0}^{n-1}$ by means of the composite isomorphism
\begin{align}
\gamma : \mcF \xrightarrow{\sim} \mcF^{(1)}\oplus \mcF^{(0)}/\mcF^{(1)} \xrightarrow{\sim} 
\mcF^{(2)} \oplus (\mcF^{(1)}/\mcF^{(2)}) \oplus (\mcF^{(0)}/\mcF^{(1)})
\xrightarrow{\sim} \cdots \xrightarrow{\sim}
\bigoplus_{j=0}^{n-1} \mcF^{(j)}/\mcF^{(j+1)},
\end{align}
where the $j$-th isomorphism (for each $j \in \{ 1, \cdots, n \}$)
arises from $\gamma_j$.
Hereinafter, we consider $\mcF$ as being equipped with a grading (indexed by $\{0, \cdots, n-1 \}$) with respect to  $\gamma$, i.e., a grading whose $j$-th component is $\mcF^{(j)}/\mcF^{(j+1)}$. 
 
\SSP
\bde \label{Dek8}
\begin{itemize}
\item[(i)]
The collection
\begin{align} \label{Eq225}
\widehat{\msF}^\heartsuit := (\mcF, \nabla, \{ \mcF^j \}_j, \{ \mcF^{(j)}\}_{j})
\end{align}
is 
called a {\bf generic Miura $\mr{GL}_n^{(\MM)}$-oper} on $X_\LL$  if  $\mcF^{(j)}$ is closed under the $\mcD^{(\MM -1)}$-action $\nabla$ for every $j=0, \cdots, n$.
\item[(ii)]
Let $\widehat{\msF}_\circ^\heartsuit := (\mcF_\circ, \nabla_\circ, \{ \mcF^j_\circ \}_j, \{ \mcF^{(j)}_\circ\}_{j})$ and $\widehat{\msF}_\bullet^\heartsuit := (\mcF_\bullet, \nabla_\bullet, \{ \mcF^j_\bullet \}_j, \{ \mcF^{(j)}_\bullet\}_{j})$
be generic Miura $\mr{GL}_n^{(\MM)}$-opers.
An {\bf isomorphism of generic Miura $\mr{GL}_n^{(\MM)}$-opers} from  $\widehat{\msF}_\circ^\heartsuit$ to $\widehat{\msF}_\bullet^\heartsuit$ is an isomorphism of $\mr{GL}_n^{(\MM)}$-opers $(\mcF_\circ, \nabla_\circ, \{ \mcF^j_\circ \}_j) \xrightarrow{\sim} (\mcF_\bullet, \nabla_\bullet, \{ \mcF_\bullet^j \}_j)$ compatible with the respective  second filtrations $\{ \mcF^{(j)}_\circ \}_j$ and $\{ \mcF^{(j)}_\bullet \}_j$. 
\end{itemize}
\ede

\SSP
\bde[cf. ~\cite{Wak12}, Definition 5.7, for the case of $\mr{GL}_n^{(\MM)}$-opers] \label{Def59}
Let $\msF^\heartsuit := (\mcF, \nabla, \{ \mcF^j \}_j)$ (resp., $\widehat{\msF}^\heartsuit := (\mcF, \nabla, \{ \mcF^j \}_j, \{ \mcF^{(j)}\}_{j})$) be a $\mr{GL}_n^{(\MM)}$-oper (resp., a generic Miura $\mr{GL}_n^{(\MM)}$-oper) on $X_\LL$.
Then, we shall say that $\msF^\heartsuit$ (resp., $\widehat{\msF}^\heartsuit$) is {\bf dormant} if $(\mcF, \nabla)$ is dormant in the sense of Definition \ref{Def500}, (i) and (ii).
\ede 
\SSP

\begin{rem} \label{Rem4891}
Suppose that $(\LL, \MM) = (1, \N)$, and  let $\N'$ be a positive integer with $\N' \leq \N$.
Consider a dormant $\mr{GL}_n^{(\N)}$-oper $\msF^\heartsuit := (\mcF, \nabla^{(\N)}, \{ \mcF^{j}\}_{j})$ on $X_1$.
Since there exists a natural morphism $\mcD^{(\N')} \rightarrow \mcD^{(\N)}$  (cf. ~\cite[\S\,2.2.1]{PBer1}), $\nabla^{(\N)}$ induces a $\mcD^{(\N')}$-module structure $\nabla^{(\N) \Rightarrow (\N')}$ on  $\mcF$ obtained by composing  this morphism.
The resulting collection $(\mcF, \nabla^{(\N) \Rightarrow (\N')}, \{ \mcF^j \}_j)$, i.e., the truncation of $\msF^\heartsuit$ to level $\N'$,  forms a dormant $\mr{GL}_n^{(\N')}$-oper.
The same holds for dormant generic Miura $\mr{GL}_n^{(\N)}$-opers.
\end{rem}

\LSP
\subsection{$\mr{PGL}_n^{(\MM)}$-opers and generic Miura $\mr{PGL}_n^{(\MM)}$-opers} \label{SS100}

We shall define equivalence relations in the sets of $\mr{GL}_n^{(\MM)}$-opers and generic Miura $\mr{GL}_n^{(\MM)}$-opers, respectively (cf. ~\cite[Definition 5.10]{Wak12}).
Let $\msF^\heartsuit := (\mcF, \nabla, \{ \mcF^j \}_j)$ (resp., $\widehat{\msF}^\heartsuit := (\mcF, \nabla, \{ \mcF^j \}_j, \{ \mcF^{(j)}\}_{j})$) be a $\mr{GL}_n^{(\MM)}$-oper (resp., a generic Miura $\mr{GL}_n^{(\MM)}$-oper) on $X_\LL$ and $\msL := (\mcL, \nabla_\mcL)$ a line  $\mcD^{(\MM -1)}$-bundle.
The tensor product $\mcF \otimes \mcL$  naturally inherits a  $\mcD^{(\MM -1)}$-module structure $\nabla \otimes \nabla_\mcL$ induced by  $\nabla$ and $\nabla_\mcL$, and 
 the collection
\begin{align}
\msF^\heartsuit_{\otimes \msL} := (\mcF \otimes \mcL, \nabla \otimes \nabla_\mcL, \{ \mcF^j \otimes \mcL \}_{j=0}^n)  \hspace{25mm}\\
\left(\text{resp.,} \  \widehat{\msF}^\heartsuit_{\otimes \msL} := (\mcF \otimes \mcL, \nabla \otimes \nabla_\mcL, \{ \mcF^j \otimes \mcL \}_{j=0}^n, \{ \mcF^{(j)} \otimes \mcL \}_{j=0}^n) \right) \notag
\end{align}
forms a $\mr{GL}_n^{(\MM)}$-oper (resp., a generic Miura $\mr{GL}_n^{(\MM)}$-oper) on $X_\LL$.
If both $\msF^\heartsuit$ (resp., $\widehat{\msF}^\heartsuit$) and $\msL$ are  dormant,
then $\msF^\heartsuit_{\otimes \msL}$ (resp., $\widehat{\msF}^\heartsuit_{\otimes \msL}$) is also  dormant.

Now, let $\msF^\heartsuit_{\circ}$ and $\msF^\heartsuit_\bullet$ (resp., $\widehat{\msF}^\heartsuit_{\circ}$ and $\widehat{\msF}^\heartsuit_\bullet$) be $\mr{GL}_n^{(\MM)}$-opers (resp., generic Miura $\mr{GL}_n^{(\MM)}$-opers) on $X_\LL$.
We shall say that $\msF^\heartsuit_\circ$ (resp., $\widehat{\msF}^\heartsuit_\circ$) is {\bf equivalent  to $\msF^\heartsuit_\bullet$} (resp., $\widehat{\msF}^\heartsuit_\bullet$)
if there exists a line $\msD^{(\MM -1)}$-bundle $\msL$ satisfying  $(\msF^\heartsuit_\circ)_{\otimes \msL} \cong \msF^\heartsuit_\bullet$ (resp.,  $(\widehat{\msF}^\heartsuit_\circ)_{\otimes \msL} \cong \widehat{\msF}^\heartsuit_\bullet$).
This  binary relation defines equivalence relations  on the  set of $\mr{GL}_n^{(\MM)}$-opers  (resp.,  generic Miura $\mr{GL}_n^{(\MM)}$-opers), as well as  on the corresponding set of   dormant objects.
For each   element $\msF^\heartsuit$ (resp., $\widehat{\msF}^\heartsuit$) in one of these   sets, 
we denote by  $[\msF^\heartsuit]$ (resp., $[\widehat{\msF}^\heartsuit]$) its equivalence class.

\SSP
\bde \label{Def79}
By a {\bf $\mr{PGL}_n^{(\MM)}$-oper} (resp., a {\bf dormant $\mr{PGL}_n^{(\MM)}$-oper}) on $X_\LL$,
we mean an equivalence class  of a $\mr{GL}_n^{(\MM)}$-oper (resp., a dormant $\mr{GL}_n^{(\MM)}$-oper) on $X_\LL$.
Similarly,  a {\bf generic Miura $\mr{PGL}_n^{(\MM)}$-oper} (resp., a {\bf dormant generic Miura $\mr{PGL}_n^{(\MM)}$-oper}) on $X_\LL$ is defined as  an equivalence class  of a generic Miura $\mr{GL}_n^{(\MM)}$-oper (resp., a dormant generic Miura $\mr{GL}_n^{(\MM)}$-oper) on $X_\LL$.
\ede 
\SSP

In the remainder  of this section,
 let us fix   a section $\sigma_\LL : \mr{Spec}(W_\LL) \rightarrow X_\LL$ of $f_\LL$ such that  the pair $\msX_\LL := (X_\LL, \sigma_\LL)$ defines  a smooth elliptic curve, and suppose that its mod $p$ reduction is {\it ordinary}.
Additionally, we  fix 
an element $\widetilde{\delta}$ of $H^0 (X_\LL, \mcB_\LL^\mr{log})$ whose reduction modulo $p$, denoted by $\delta$, belongs to $H^0 (X_1, \mcB^\mr{log}_1)$, or equivalently, satisfies $(\delta^\vee)^{[p]} = \delta^\vee$ (cf. Remark \ref{Rem912}).

We shall write
\begin{align} \label{Eq1093}
\mcO p_{\LL, \MM}^{^\mr{Zzz...}} \ \left(\text{resp.,} \  \mcM \mcO p_{\LL, \MM}^{^\mr{Zzz...}} \right)
\end{align}
for the set of dormant $\mr{PGL}_n^{(\MM)}$-opers (resp., dormant generic Miura $\mr{PGL}_n^{(\MM)}$-opers) on $X_\LL$.
By forgetting the second  filtration in the data of  each generic Miura $\mr{PGL}_n^{(\MM)}$-oper, we obtain a map of sets
\begin{align} \label{Eq1001}
\mcM \mcO p_{\LL, \MM}^{^\mr{Zzz...}} \rightarrow \mcO p_{\LL, \MM}^{^\mr{Zzz...}}.
\end{align}

\LSP
\subsection{Comparison with dormant $\mr{PGL}_n$-opers} \label{SS8e8}

We now return to the setting of the previous section and 
consider the case where $G = \mr{PGL}_n$
with $2n < p$.
Let  $B$ (resp., $T$)
  be the Borel subgroup (resp., the maximal torus) of $\mr{PGL}_n$ 
defined as  the images, via the natural quotient $\mr{GL}_n \twoheadrightarrow \mr{PGL}_n$,
 of 
 invertible  upper-triangular matrices (resp., invertible diagonal matrices).

For each commutative ring $R$, let us define 
$\tau_\mr{reg} (R)$ as
 the  set of $n$-tuples $(a_1, \cdots, a_n)$ of elements in $R$ such that their reductions  modulo  $p$,  denoted by  $\vec{a} := (\overline{a}_1, \cdots, \overline{a}_n)$, are pairwise distinct  in $\mbF_p \otimes R$.
For two elements $\vec{a}_\circ := (a_{\circ 1}, \cdots, a_{\circ n})$ and  $\vec{a}_\bullet := (a_{\bullet 1}, \cdots, a_{\bullet n})$ in $\tau_{\mr{reg}}(R)$,
we say that $\vec{a}_\circ$ is {\bf equivalent to $\vec{a}_\bullet$}
if there exists an element $c \in R$ with $\vec{a}_\circ = \vec{a}_\bullet + c$, where $\vec{a}_\bullet + c := (a_{\bullet 1}+c, \cdots, a_{\bullet n} + c)$.
The set  of equivalence classes of elements in $\tau_\mr{reg}(R)$ is denoted by 
\begin{align} \label{eqw22}
\overline{\tau}_\mr{reg}(R).
\end{align}
Note that $\tau_\mr{reg}(R)$ carries  the action of the symmetric group  of $n$ letters $\mfS_n$  by permutation.
 This action induces a well-defined $\mfS_n$-action on $\overline{\tau}_\mr{reg}(R)$, allowing us to define the quotient set
\begin{align}
\overline{\tau}_\mr{reg} (R)/\mfS_n.
\end{align}

Since the Lie algebra $\mft$ of $T$ can be defined over $\mbF_p$,
 it makes sense to consider  the set of $\mbF_p$-rational points $\mft_\mr{reg} (\mbF_p)$ in $\mft_\mr{reg}$.
 The sets $\mft_\mr{reg} (\mbF_p)$ and $\mft_\mr{reg} (\mbF_p)/W$ can be naturally  identified with
 $\overline{\tau}_\mr{reg} (\mbF_p)$ and $\overline{\tau}_\mr{reg} (\mbF_p)/\mfS_n$, respectively.

Now, let us specialize to the case 
$\N =1$, i.e., $(\LL, \MM) = (1, 1)$, and 
  take a $\mr{PGL}_n^{(1)}$-oper $[\msF^\heartsuit]$ on $X_1$, where $\msF^\heartsuit := (\mcF, \nabla, \{ \mcF^j \}_j)$.
The vector bundle $\mcF$ determines, via  $\mr{GL}_n \twoheadrightarrow \mr{PGL}_n$,
a $\mr{PGL}_n$-bundle $\mcE_{\mr{PGL}_n}$.
The filtration $\{ \mcF^j \}_j$ corresponds to a $B$-reduction $\mcE_B$ of $\mcE_{\mr{PGL}_n}$ and the $\mcD^{(0)}$-module structure $\nabla$ corresponds to a connection $\nabla_\mcE$ on $\mcE_{\mr{PGL}_n}$.
The resulting pair $\msE^\spadesuit := (\mcE_B, \nabla_\mcE)$ defines  a $\mr{PGL}_n$-oper on $\msX_1$ in the sense of Definition \ref{Def1}, and its  isomorphism class  depends only on the equivalence class $[\msF^\heartsuit]$.
The  assignment $[\msF^\heartsuit] \mapsto \msE^\spadesuit$ gives a well-defined bijection  between dormant $\mr{PGL}_n$-opers and dormant $\mr{PGL}_n^{(1)}$-opers.
That is to say,
 there exists a canonical bijection
\begin{align} \label{Eq184}
\mcO p_{1, 1}^{^\mr{Zzz...}} \xrightarrow{\sim} \Pi^{-1}_{\mr{PGL}_n} (q)   \ \left(\text{resp.,} \ \mcM\mcO p_{1, 1}^{\mr{Zzz...}}  \xrightarrow{\sim}  \widehat{\Pi}^{-1}_{\mr{PGL}_n} (q) \right),
\end{align}
where $q$ denotes the $k$-rational point of $\overline{\mcM}_\mr{ell}$ classifying $\msX_1$.
By means of   this bijection and   the natural identification $\mft_\mr{reg}(\mbF_p) /W = \overline{\tau}_\mr{reg}(\mbF_p)/\mfS_n$ (resp., $\mft_\mr{reg}(\mbF_p) = \overline{\tau}_\mr{reg}(\mbF_p)$),
the fiber of $\Theta_G$ (resp., $\widehat{\Theta}_G$) over the point $(\msX_1, \delta) \in \mcN_\mr{inv}$ induces 
 a bijection
\begin{align} \label{Eq1004}
\Theta_{1, 1} : \overline{\tau}_{\mr{reg}}(\mbF_p) /\mfS_n \xrightarrow{\sim} \mcO p_{1, 1}^{^\mr{Zzz...}} \ \left(\text{resp.,} \  \widehat{\Theta}_{1, 1} : \overline{\tau}_{\mr{reg}}(\mbF_p)\xrightarrow{\sim} \mcM \mcO p_{1, 1}^{^\mr{Zzz...}} \right).
\end{align}
In the following discussion, we generalize $\Theta_{1, 1}$ (resp., $\widehat{\Theta}_G$)
to prime-power characteristic, as well as to higher level.

\LSP
\subsection{The case of $(\LL, \MM) = (\N, 1)$} \label{SS58}

For each $a \in W_\N (\mbF_p)\left(= \mbZ/p^\N \mbZ \right)$, 
let $\nabla_{a}^{(0)}$ denote  the dormant connection  on $\mcO_{X_{\N}}$
corresponding to the section  $a \cdot \widetilde{\delta} \in H^0 (X_\N, \mcB_\N^\mr{log})$ 
   via \eqref{Eq1014}.
 This means that   $\nabla_a^{(0)} = d + a \cdot \widetilde{\delta}$ and $c (\nabla_a^{(0)}) = a \cdot \widetilde{\delta}$.

Let us take  an $n$-tuple $\vec{a} := (a_1, \cdots, a_n)$ of elements in $W_\N (\mbF_p)$.
The direct sum $\nabla_{\vec{a}}^{(0)} := \bigoplus_{i=1}^n \nabla_{a_i}^{(0)}$ defines  a dormant  connection on $\mcG_{\N} := \mcO_{X_{\N}}^{\oplus n}$.
Given each $j=0, \cdots, n$, we denote by  $\mcG_{\N, \vec{a}}^j$  the $\mcO_{X_{\N}}$-submodule of $\mcG_\N$ generated by the sections $(\nabla_{\vec{a}}^{(0)})_{\widetilde{\delta}^\vee}^{n-j-1} (v)$, where $v$ runs over all sections  in the image of the diagonal embedding $\Delta_{\N} : \mcO_{X_{\N}} \hookrightarrow \mcG_{\N}$.
Here, we set   
$\mcG_{\N, \vec{a}}^n := 0$ by convention.
Additionally, let
 $\mcG_{\N}^{(j)}$ be  the image of the standard  inclusion into the first $n-j$ factors $\mcO_{X_{\N}}^{\oplus (n-j)} \hookrightarrow \mcO_{X_{\N}}^{\oplus n}$.
With these definitions,  we obtain the  collection of data
\begin{align} \label{Eq1020}
\msG^\heartsuit_{\N, \vec{a}} := (\mcG_{\N}, \nabla^{(0)}_{\vec{a}}, \{ \mcG_{\N, \vec{a}}^j \}_{j=0}^n)  \ 
\left(\text{resp.,} \ \widehat{\msG}^\heartsuit_{\N, \vec{a}} := (\mcG_{\N}, \nabla^{(0)}_{\vec{a}}, \{ \mcG_{\N, \vec{a}}^j \}_{j=0}^n, \{\mcG^{(j)}_{\N}\}_{j=0}^n) \right).
\end{align}

\SSP
\bpr \label{Prop1001}
We continue with the notation introduced above.
Then,
 $\vec{a}$ belongs to $\tau_\mr{reg} (W_\N (\mbF_p))$ if and only if $\msG^\heartsuit_{\N, \vec{a}}$  (resp., $\widehat{\msG}^\heartsuit_{\N, \vec{a}}$) defines a dormant $\mr{GL}_n^{(1)}$-oper (resp., a dormant generic Miura $\mr{GL}_n^{(1)}$-oper).
\epr
\begin{proof}
We first prove the non-resp'd assertion.
For each $j=0, \cdots, n-1$, consider the column vector $\vec{v}_j := (\nabla_{\vec{a}}^{(0)})_{\widetilde{\delta}^\vee}^{n-j-1} (\Delta_{\N} (1) )$.
By explicit computation, we have 
  $\vec{v}_j = {^t}(a_1^{n-j-1}, \cdots, a_n^{n-j-1})$.
The determinant of the matrix formed by these vectors is given by the Vandermonde determinant formula:
\begin{align}
\mr{det} \left(\vec{v}_0, \cdots, \vec{v}_{n-1} \right) = (-1)^{n-1} \cdot \prod_{i < j} (a_i -a_j).
\end{align}
Thus, the vectors $\vec{v}_0, \cdots, \vec{v}_{n-1}$ forms a basis (or equivalently  
$\{ \mcG^j_{\N, \vec{a}} \}_{j=0}^n$ defines a complete flag on $\mcG_{\N}$) if and only if 
$\prod_{i < j} (a_i -a_j)$ lies in $W_\N (\mbF_p)^\times$.
This condition  is precisely the requirement  that $\vec{a} \in \tau_\mr{reg}(W_\N (\mbF_p))$, which  completes the proof.

The resp'd portion   follows  immediately from  the fact just proved  and  the definition of  the second filtration $\{ \mcG_{\N}^{(j)} \}_{j}$.
\end{proof}
\SSP

Now, suppose that $\vec{a} \in \tau_\mr{reg} (W_\N (\mbF_p))$, and let $c \in W_\N (\mbF_p)$ be given.
Then, we have  $(\msG_{\N, \vec{a}}^\heartsuit)_{\otimes \msO_c^{(0)}} \cong \msG^\heartsuit_{\N, \vec{a}+c}$, where $\msO_c^{(0)} := (\mcO_{X_\N}, \nabla^{(0)}_{c})$ and $\vec{a} + c := (a_1 +c, \cdots, a_n + c)$.
This implies  $[\msG_{\N, \vec{a}}^\heartsuit] = [\msG_{\N, \vec{a}+c}^\heartsuit]$.
Furthermore, for  
any permutation  $\sigma \in \mfS_n$, 
the corresponding permutation of the direct summands in $\msG_\N \left(= \mcO_{X_\N}^{\oplus n} \right)$ induces  an isomorphism
$\msG^\heartsuit_{\N, \vec{a}} \xrightarrow{\sim}\msG^\heartsuit_{\N, \sigma (\vec{a})}$, where $\sigma (\vec{a}) := (a_{\sigma (1)}, \cdots, a_{\sigma (n)})$.
These observations  lead to   the following conclusion.

 \SSP
 \bpr \label{T39}
 The assignment $\vec{a} \mapsto [\msG^\heartsuit_{\N, \vec{a}}]$ (resp., $\vec{a} \mapsto [\widehat{\msG}^\heartsuit_{\N, \vec{a}}]$)  for each  $\vec{a} \in \tau_\mr{reg}(W_\N (\mbF_p))$
 determines a well-defined bijection
 \begin{align}
\Theta_{\N, 1} :  \overline{\tau}_\mr{reg}(W_\N (\mbF_p))/\mfS_n \xrightarrow{\sim}  \mcO p_{\N, 1}^{^\mr{Zzz...}} \ \left(\text{resp.,} \ \widehat{\Theta}_{\N, 1} :  \overline{\tau}_\mr{reg}(W_\N (\mbF_p)) \xrightarrow{\sim}  \mcM \mcO p_{\N, 1}^{^\mr{Zzz...}} \right).
 \end{align}
  Moreover, 
 the following square diagram commutes:
 \begin{align} \label{Eq203}
\vcenter{\xymatrix@C=46pt@R=36pt{
 \overline{\tau}_\mr{reg} (W_\N (\mbF_p)) \ar[r]^-{\widehat{\Theta}_{\N, 1}}_-{\sim} \ar[d]_-{\mr{quotient}} &   \mcM\mcO p^{^\mr{Zzz...}}_{\N, 1}\ar[d]^-{\eqref{Eq1001}} \\
 \overline{\tau}_\mr{reg} (W_\N (\mbF_p))/\mfS_n \ar[r]_-{\Theta_{\N, 1}}^-{\sim} & \mcO p^{^\mr{Zzz...}}_{\N, 1}.
 }}
\end{align}
 \epr
\begin{proof}
The well-definedness of $\Theta_{\N, 1}$ (resp., $\widehat{\Theta}_{\N, 1}$) follows from the  preceding discussion.
Since the proof of the bijecitivity of $\widehat{\Theta}_{\N, 1}$ is simpler than that of $\Theta_{\N, 1}$,
we only consider the latter one.

First,  to consider the injectivity of $\Theta_{\N, 1}$,
let us take two elements $\vec{a}_\circ := (a_{\circ1}, \cdots, a_{\circ n})$,
$\vec{a}_\bullet := (a_{\bullet 1}, \cdots, a_{\bullet n})$ with $[\msG_{\vec{a}_{\circ}}^\heartsuit] = [\msG^\heartsuit_{\vec{a}_\bullet}]$.
Then, there exists an element $c \in W_\N (\mbF_p)$ admitting   an isomorphism $\msG^\heartsuit_{\vec{a}_\circ} \xrightarrow{\sim} \msG^\heartsuit_{\vec{a}_\bullet + c}$.
By Proposition \ref{Prop449},
the set $\{ a_{\circ 1}, \cdots a_{\circ n} \}$ coincides with $\{ a_{\bullet 1}+c, \cdots a_{\bullet n} +c \}$.
This means  that $\vec{a}_{\circ} = \sigma (\vec{a}_{\bullet} + c)$ for some permutation $\sigma \in \mfS_n$, which implies $\vec{a}_\circ = \vec{a}_\bullet$ in $\overline{\tau}_\mr{reg}(W_\N (\mbF_p))/\mfS_n$.
Thus, $\Theta_{\N, 1}$ is injective.

Next, we shall  consider the surjectivity of $\Theta_{\N, 1}$.
Let $\msF^\heartsuit := (\mcF, \nabla, \{ \mcF^j \}_j)$ be a dormant $\mr{GL}_n^{(1)}$-oper representing an element of $\mcO p_{\N, 1}^{^\mr{Zzz...}}$.
Note that the morphism $\mcD_{\leq n-1}^{(0)} \otimes \mcF^{n-1} \rightarrow \mcF$  (cf. \eqref{Eq1100})
obtained by restricting $\nabla : \mcD^{(0)} \otimes \mcF \rightarrow \mcF$ is an isomorphism.
Since  $\widetilde{\delta}$ gives a trivialization ${^R}\mcD^{(0)}_{\leq n-1} \cong \mcO_{X_{\N}}^{\oplus n}$,
we have  
\begin{align} \label{Eq1101}
\mcF \cong \mcD_{\leq n-1}^{(0)} \otimes \mcF^{n-1} \cong (\mcF^{n-1})^{\oplus n}.
\end{align}
This implies $\mr{det}(\mcF) \cong \mr{det}((\mcF^{n-1})^{\oplus n}) \cong (\mcF^{n-1})^{\otimes n}$.
Hence, the determinant  $\mr{det}(\nabla)$  of $\nabla$ can be regarded as a connection on $(\mcF^{n-1})^{\otimes n}$.
It is verified that there exists a unique connection $\nabla_{\mcF^{n-1}}$ on $\mcF^{n-1}$ whose $n$-fold  tensor product coincides with $\mr{det}(\nabla)$ (cf. ~\cite[Proposition 4.22, (i)]{Wak5}).
 After possibly replacing $\msF^\heartsuit$ with $\msF^\heartsuit_{\otimes (\mcF^{n-1}, \nabla_{\mcF^{n-1}})}$, we may assume that
 $\mcF^{n-1} = \mcO_X$, leading to $\mcF = \mcO_X^{\oplus n}$ by \eqref{Eq1101}.
 Then,  $\nabla$ takes the form  $\nabla = d + \widetilde{\delta} \otimes A$ for some $A \in M_n (W_N (k)) = \mr{End}_{\mcO_{X_{\N }}} (\mcO_{X_{\N }}^{\oplus n})$.
 Since  $\nabla$ is dormant and $(\delta^\vee)^{[p]} = \delta^{\vee p} = \delta^\vee$, the mod $p$ reduction $\overline{A}$ of $A$ satisfies  
 \begin{align}
 0 = \psi^1 (\nabla) ((\delta^\vee)^{\otimes p})= (\nabla_{\delta^\vee})^p - \nabla_{{\delta^{\vee p}}} = \overline{A}^p -  \overline{A}.
 \end{align} 
 Let   $f (t) \in W_\N  [t]$ be 
 the characteristic polynomial of $A$.
 Since  the roots of  its mod $p$ reduction  $\overline{f} (t)$ lies in the set of roots of   $t^p-t = \prod_{i \in \mbF_p} (t -i)$, we have  $\overline{f} (t) = \prod_{i=1}^n (t - \overline{a}_i)$ for some $\overline{a}_1, \cdots, \overline{a}_n \in \mbF_p$.
Just as in the proof of  Proposition \ref{Prop19},
 $\overline{A}$ turns out to be  diagonalized.
Thus, there exists  an invertible  matrix $P \in M_n (k)$ such that $P^{-1} \overline{A} P$ coincides with the diagonal 
  matrix $D$ with diagonal entries $\overline{a}_1, \cdots, \overline{a}_n$.
 Applying  the gauge transformation associated to  $P$, the mod $p$ reduction $\overline{\nabla}$ of $\nabla$ transforms  into
  \begin{align} \label{Eq1023}
  {^P}\overline{\nabla} = P^{-1} (d + \delta \otimes \overline{A}) P = P^{-1}dP + \delta \otimes P^{-1} \overline{A} P = \widetilde{\delta} \otimes D.
  \end{align}
  It follows that $(\mcO_{X_1}^{\oplus n}, \nabla)$ is isomorphic to $\bigoplus_{i=1}^n (\mcO_{X_1}, \nabla_{\overline{a}_i})$.
 By Proposition \ref{Prop1001}, the elements $\overline{a}_1, \cdots, \overline{a}_n$ are mutually distinct.
 In particular, the polynomial $f (t)$ is separable, so Hensel's lemma implies that
 $f (t)$ decomposes as  the product $\prod_{i=1}^n (t - a_i)$, where each $a_i$ denotes an element of $W_\N (\mbF_p)$ whose mod $p$ reduction coincides with $\overline{a}_i$.
Similarly to    \eqref{Eq1023},
there exists  an isomorphism of $\mcD^{(0)}$-bundles $\alpha : (\mcO_{X_{\N}}^{\oplus n}, \nabla) \xrightarrow{\sim}\bigoplus_{i=1}^n (\mcO_{X_{\N}}, \nabla_{a_i})$.
Since each connection $\nabla_{a_i}$ remains  invariant under gauge transformations by elements of $\mr{Aut}_{\mcO_{X_{\N}}} (\mcO_{X_{\N}}) \left(\cong W_\N^\times \right)$,  we can choose  $\alpha$ such that  $\alpha (\mr{Im} (\Delta_{\N})) = \mr{Im} (\Delta_{\N})$.
Consequently, we obtain  $\msF^\heartsuit \cong \msG^\heartsuit_{\N, \vec{a}}$, where $\vec{a} := (a_1, \cdots, a_n)$.
This proves  the required  surjectivity of $\Theta_{\N, 1}$.
 \end{proof}

\LSP
\subsection{The case of $(\LL, \MM) = (1, \N)$} \label{SS158}

For each $a \in W_\N (\mbF_p)$,
denote by $\nabla_a^{(\N -1)}$ the  $\mcD^{(\N -1)}$-module structure on $\mcO_{X_1}$ corresponding to the section $a \cdot \widetilde{\delta} \in H^0 (X_\N, \mcB_\N^\mr{log})$ via  
the composite bijection
\begin{align}
H^0 (X_\N, \mcB_\N^\mr{log}) \xrightarrow{\eqref{Eq1050}} \mr{Ker} (\tau) \xrightarrow{\eqref{Eq1053}} J_{X_1} [p^\N] (k) \xrightarrow{\eqref{Eq1011}} \mr{Conn}_{1, \N}^{^\mr{Zzz...}}. 
\end{align}

Now, consider  an $n$-tuple
 $\vec{a} := (a_1, \cdots, a_n)$ of elements in 
 $W_\N (\mbF_p)$.
 The direct sum $\nabla_{\vec{a}}^{(\N -1)} := \bigoplus_{i=1}^n \nabla_{a_i}^{(\N -1)}$ specifies a $\mcD^{(\N -1)}$-module structure on  $\mcG_{1}\left(= \mcO_{X_1}^{\oplus n}\right)$ with vanishing $p^\N$-curvature.
For each  $j = 0, \cdots, n$,
 define   $\mcG_{1, \vec{a}}^j$ as   the image of 
the composite
\begin{align}
\mcD^{(\N -1)}_{\leq n-j -1} \otimes \mr{Im} (\Delta_1) \xrightarrow{\mr{inclusion}} \mcD^{(\N -1)} \otimes \mcG_1 \xrightarrow{\nabla^{(\N -1)}_{\vec{a}}} \mcG_1.
\end{align}
Moreover, let $\mcG_{1}^{(j)}$ be   the image of the inclusion into the first $n-j$ factors $\mcO_{X_1}^{\oplus (n-j)} \hookrightarrow \mcO_{X_1}^{\oplus n}$.
This gives rise to  the  collection of data
\begin{align} \label{Eq1109}
\msG_{1, \vec{a}}^\heartsuit := (\mcG_1, \nabla^{(\N -1)}_{\vec{a}}, \{ \mcG_{1, \vec{a}}^{j} \}_{j=0}^n) \ \left(\text{resp.,} \  \widehat{\msG}_{1, \vec{a}}^\heartsuit := (\mcG_1, \nabla^{(\N -1)}_{\vec{a}}, \{ \mcG_{1, \vec{a}}^{j} \}_{j=0}^n, \{ \mcG_{1}^{(j)}\}_{j=0}^n) \right).
\end{align}
By applying  Proposition  \ref{Prop1001} to the level $1$ truncation of this data,
we see  that
$\msG_{1, \vec{a}}^\heartsuit$ (resp., $\widehat{\msG}_{1, \vec{a}}^\heartsuit$) defines  a dormant $\mr{GL}_n^{(\N)}$-oper (resp., a dormant generic Miura $\mr{GL}_n^{(\N)}$-oper) if and only if $\vec{a} \in \tau_\mr{reg} (W_\N (\mbF_p))$.

Furthermore, for each $c \in W_\N (\mbF_p)$,
we have $(\msG^\heartsuit_{1, , \vec{a}})_{\otimes \msO^{(\N -1)}_{c}} \cong \msG^\heartsuit_{1, \vec{a} + c}$, where $\msO^{(\N -1)}_{c} := (\mcO_{X_1}, \nabla^{(\N -1)}_{c})$, and 
 this  implies  $[\msG^\heartsuit_{1, \vec{a}}] = [\msG^\heartsuit_{1, \vec{a} + c}]$.
 Additionally, 
 any permutation $\sigma \in \mfS_n$  of the  direct  summands in $\msG_1 \left(=  \mcO_{X_1}^{\oplus n} \right)$ yields  an isomorphism $\msG^\heartsuit_{1, \vec{a}} \xrightarrow{\sim} \msG^\heartsuit_{1, \sigma (\vec{a})}$.
 Therefore,
the following assertion holds.

 \SSP
 \bpr \label{T34}
 The assignment $\vec{a} \mapsto [\msG^\heartsuit_{1, \vec{a}}]$ (resp., $\vec{a} \mapsto [\widehat{\msG}_{1, \vec{a}}^\heartsuit]$) defined for each $\vec{a} \in \tau_\mr{reg}(W_\N (\mbF_p))$ determines a well-defined bijection
 \begin{align}
 \Theta_{1, \N} : \overline{\tau}_\mr{reg}(W_\N (\mbF_p))/\mfS_n \xrightarrow{\sim} \mcO p_{1, \N}^{^\mr{Zzz...}} \ \left(\text{resp.,} \ 
  \widehat{\Theta}_{1, \N} : \overline{\tau}_\mr{reg}(W_\N (\mbF_p)) \xrightarrow{\sim}  \mcM\mcO p_{1, \N}^{^\mr{Zzz...}}.
  \right)
 \end{align}
 Moreover,   the following square diagram commutes:
 \begin{align} \label{Eq211}
\vcenter{\xymatrix@C=46pt@R=36pt{
 \overline{\tau}_\mr{reg} (W_\N (\mbF_p)) \ar[r]^-{\widehat{\Theta}_{1, \N}}_-{\sim} \ar[d]_-{\mr{quotient}} &   \mcM\mcO p^{^\mr{Zzz...}}_{1, \N}\ar[d]^-{\eqref{Eq1001}} \\
 \overline{\tau}_\mr{reg} (W_\N (\mbF_p))/\mfS_n \ar[r]_-{\Theta_{1, \N}}^-{\sim} & \mcO p^{^\mr{Zzz...}}_{1, \N}.
 }}
\end{align}
 \epr
\begin{proof}
The well-definedness of $\Theta_{1, \N}$ (resp., $\widehat{\Theta}_{1, \N}$) follows from  the preceding discussion.
Since the proof of the bijectivity of $\widehat{\Theta}_{1, \N}$ is simpler than that of $\Theta_{1, \N}$, we only consider the latter one.

The injectivity of 
 $\Theta_{1, \N}$ follows from an argument analogous  to the proof of the  injectivity of $\Theta_{\N, 1}$ (cf. the proof of Proposition \ref{T39}).

 Next, let us consider the surjectivity.
Let $\msF^\heartsuit := (\mcF, \nabla, \{ \mcF^j \}_j)$ be  a dormant $\mr{GL}_n^{(\N)}$-oper   on $X_1$ representing an element of $\mcO p_{1, \N}^{^\mr{Zzz...}}$.
 Just as in the proof of the surjectivity of  $\Theta_{\N, 1}$ (cf. the proof of Proposition \ref{T39}),
 we may  assume, without loss of generality,  that
$\mcF = \mcO_{X_1}^{\oplus n}$.
Under this assumption, we have  $F^{\N*}_{X_1} (\mcS ol (\nabla)) \cong \mcO_{X_1}^{\oplus n}$, and it follows from 
  ~\cite[Theorem 2.16]{Oda}  that  every  indecomposable component of $\mcS ol (\nabla)$ must have rank $1$.
 Hence,  there exist line bundles $\mcL_1, \cdots, \mcL_n$ with $\mcS ol (\nabla) \cong \bigoplus_{i=1}^n \mcL_i$.
 It follows that $F_{X_1}^{\N*} (\mcL_i) \cong \mcO_{X_1}$ for every $i$.
 By the bijection \eqref{Eq1011},
 there exists an $n$-tuple  $\vec{a} := (a_1, \cdots, a_n)\in W_\N (\mbF_p)^{\oplus n}$ admitting  an isomorphism of $\mcD^{(\N -1)}$-bundles $\beta : (\mcG_1, \bigoplus_{i=1}^n\nabla_{a_i}^{(\N -1)}) \xrightarrow{\sim} (\mcF, \nabla)$.
 Since each connection $\nabla_{a_i}^{(\N -1)}$ is invariant under  gauge transformations by elements of $\mr{Aut}_{\mcO_{X_1}}(\mcO_{X_1}) \left(\cong k^\times \right)$,  we can choose $\beta$ such that  $\beta (\mr{Im}(\mr{\Delta_1})) = \mr{Im} (\Delta_1)$.
 This implies $[\msF^\heartsuit ] = [\msG^\heartsuit_{1, \vec{a}}]$.
 By applying Proposition \ref{Prop1001}  to the level $1$  truncaton of $\msG^\heartsuit_{1, \vec{a}}$, we conclude  that $\vec{a} \in \tau_\mr{reg}(W_\N (\mbF_p))$, thus completing 
  the proof of the required  surjectivity.
\end{proof}

\LSP
\subsection{Diagonal reduction/lifting} \label{SS258}

We recall the notion of  diagonal reduction  (cf.  ~\cite[\S\,3]{Wak12}) in the non-logarithmic setting.
Let $\msF := (\mcF, \nabla_\N^{(0)})$ be a dormant $\mcD^{(0)}$-bundle on $X_{\N}$.
For each open subscheme $U$ of $X_1 \left(= k \times_{W_\N} X_\N \right)$, we define  $\mcV' (U)$ to be the set of sections $v \in \mcF_1 (U)$ admitting  a lifting $\widetilde{v}$ in $\mr{Ker} (\nabla_\N^{(0)}) (U)$, where we regard $U$ as an open subscheme of $X_\N$ via the underlying homeomorphism of  $\iota_1 : X  \hookrightarrow X_\N$.
  The Zariski sheaf $\mcV$ on $X_1$ associated to the presheaf $U \mapsto \mcV' (U)$  is endowed with  an $\mcO_{X_1}$-module structure via push-forward along   
$F^{\N}_{X_1}$.
  The morphism  $F_{X_1}^{\N*}(\mcV) \rightarrow \mcF_1$ corresponding to the inclusion $\mcV \hookrightarrow F_{X_1*}^{\N} (\mcF_1)$ via  the adjunction relation $F_{X_1}^* (-) \dashv F_{X_1 *} (-)$ turns out to be an isomorphism.
Through  this isomorphism, 
$\nabla_{\mcV, \mr{can}}^{(\N -1)}$ (cf. \eqref{E445}) is  transferred to  a dormant $\mcD^{(\N -1)}$-module structure
 ${^{\Diag\!\!}}\nabla$   on $\mcF_1$.
Thus,
we have obtained  a $\mcD^{(\N -1)}$-bundle
\begin{align}
{^{\Diag\!\!}}\msF := (\mcF_1, {^{\Diag\!\!}}\nabla),
\end{align}
which is called   the {\bf diagonal reduction} of $\msF$ (cf. ~\cite[Definition 3.7]{Wak12}).

\SSP
\bpr \label{Prop429}
Consider the map $\mr{Conn}_{\N, 1}^{^\mr{Zzz...}} \rightarrow \mr{Conn}_{1, \N}^{^\mr{Zzz...}}$ given by taking the diagonal reductions $\msL \mapsto {^{\Diag\!\!}}\msL$.
Then, this map defines an isomorphism of abelian groups   and makes the following diagram commute:
 \begin{align} \label{Eq215}
\vcenter{\xymatrix@C=46pt@R=36pt{
&  J_{X_1} [p^\N] (k)\ar[ld]^-{\sim}_-{\eqref{Eq1060}} \ar[rd]_-{\sim}^-{\eqref{Eq1011}}& \\
\mr{Conn}_{\N, 1}^{^\mr{Zzz...}} \ar[rr]_-{\sim}  &&\mr{Conn}_{1, \N}^{^\mr{Zzz...}}
 }}
\end{align}
That is to say,  in the notation of \S\,\ref{SS58}-\ref{SS158}
the equality
\begin{align}
{^{\Diag\!\!}}\nabla_{a}^{(0)} = \nabla_{a}^{(\N -1)}
\end{align}
holds for every $a \in W_\N (\mbF_p)$.
\epr
\begin{proof}
The assertion follows from the various definitions involved.
\end{proof}
\SSP

For a dormant $\mr{GL}^{(1)}_n$-oper $\msF^\heartsuit := (\mcF, \nabla, \{ \mcF^j \}_j)$ on $X_\N$,
the collection
\begin{align} \label{Eq1112}
{^{\Diag\!\!}}\msF^\heartsuit := (\mcF_1,  {^{\Diag\!\!}}\nabla, \{ \mcF_1^j \}_j)
\end{align}
forms  a dormant $\mr{GL}_n^{(\N)}$-oper on $X_1$.
If, moreover, we are given a dormant line $\mcD^{(0)}$-bundle $\msL := (\mcL, \nabla_{\mcL})$ on $X_\N$, then
${^{\Diag\!\!}}(\msF^\heartsuit_{\otimes \msL})$ is isomorphic to $({^{\Diag\!\!}}\msF^\heartsuit)_{\otimes {\tiny{^{\Diag\!\!}}\msL}}$.
It follows that 
the dormant $\mr{PGL}_n^{(\N)}$-oper ${^{\Diag\!\!}}[\msF^\heartsuit] :=  [{^{\Diag\!\!}}\msF^\heartsuit]$ on $X_1$ depends only on the equivalence class $[\msF^\heartsuit]$.
 We refer to ${^{\Diag\!\!}}\msF^\heartsuit$ (resp., ${^{\Diag\!\!}}[\msF^\heartsuit]$) as the {\bf diagonal reduction} of  $\msF^\heartsuit$ (resp., $[\msF^\heartsuit]$).

In a similar vein,
the diagonal reduction ${^{\Diag\!\!}}\widehat{\msF}^\heartsuit$ (resp., ${^{\Diag\!\!}}[\widehat{\msF}^\heartsuit]$) of a dormant generic Miura $\mr{GL}_n^{(1)}$-oper $\widehat{\msF}^\heartsuit$ (resp., a dormant generic Miura $\mr{PGL}_n^{(1)}$-oper $[\widehat{\msF}^\heartsuit]$) is well-defined.
The resulting assignments $[\msF^\heartsuit] \mapsto {^{\Diag\!\!}}[\msF^\heartsuit]$ and $[\widehat{\msF}^\heartsuit] \mapsto {^{\Diag\!\!}}[\widehat{\msF}^\heartsuit]$
determine maps of sets 
\begin{align} \label{Eq1123}
\Diag_{\!\!\N} :  \mcO p_{\N,  1}^{^\mr{Zzz...}} \rightarrow  \mcO p_{1, \N}^{^\mr{Zzz...}} \ \ \ \text{and} \ \ \ 
\widehat{\Diag}_{\!\!\N} : \mcM\mcO p_{\N, 1}^{^\mr{Zzz...}} \rightarrow \mcM \mcO p_{1, \N}^{^\mr{Zzz...}},
 \end{align}
 respectively, and 
 the following square diagram commutes:
 \begin{align} \label{Eq222}
\vcenter{\xymatrix@C=46pt@R=36pt{
 \mcM \mcO p_{\N, 1}^{^\mr{Zzz...}} \ar[r]^-{\widehat{\Diag}_{\!\!\N}} \ar[d]_-{\eqref{Eq1001}} &   \mcM\mcO p^{^\mr{Zzz...}}_{1, \N}\ar[d]^-{\eqref{Eq1001}} \\
\mcO p_{\N, 1}^{^\mr{Zzz...}} \ar[r]_-{\Diag_{\!\!\N}} & \mcO p^{^\mr{Zzz...}}_{1, \N}.
 }}
\end{align}

 \SSP
 \bt \label{T47}
 The following diagrams are commutative:
  \begin{align} \label{Eq1126}
\vcenter{\xymatrix@C=6pt@R=36pt{
 & \overline{\tau}_\mr{reg} (W_\N (\mbF_p))/\mfS_n\ar[ld]^-{\sim}_-{\Theta_{\N, 1}}\ar[rd]^-{\Theta_{1, \N}}_-{\sim}& \\
 \mcO p_{\N, 1}^{^\mr{Zzz...}}\ar[rr]_-{\Diag_{\!\!\N}}&& \mcO p_{1, \N}^{^\mr{Zzz...}},
 }}
 \hspace{5mm}
 \vcenter{\xymatrix@C=6pt@R=36pt{
 & \overline{\tau}_\mr{reg} (W_\N (\mbF_p))\ar[ld]^-{\sim}_-{\widehat{\Theta}_{\N, 1}}\ar[rd]^-{\widehat{\Theta}_{1, \N}}_-{\sim}& \\
 \mcM\mcO p_{\N, 1}^{^\mr{Zzz...}}\ar[rr]_-{\widehat{\Diag}_{\!\!\N}}&&\mcM \mcO p_{1, \N}^{^\mr{Zzz...}}.
 }}
\end{align}
  In particular, the two maps displayed in  \eqref{Eq1123} are bijective.
  \et
\begin{proof}
Since the equality ${^{\Diag\!\!}}\msG_{\N, \vec{a}}^\heartsuit = \msG^\heartsuit_{1, \vec{a}}$ holds for every $\vec{a} \in \tau_\mr{reg}(W_\N (\mbF_p))$,
the commutativity of the diagrams in \eqref{Eq1126} follows from the bijectivity of 
$\Theta_{\N, 1}$, $\Theta_{1, \N}$, $\widehat{\Theta}_{\N, 1}$,  and $\widehat{\Theta}_{1, \N}$, as established in  Theorems \ref{T39} and \ref{T34}.
Moreover,  this implies  that the maps in \eqref{Eq1123} are bijective.
\end{proof}
\SSP

Now, let $[\msF^\heartsuit]$ be a dormant $\mr{PGL}_n^{(\N)}$-oper on $X_1$.
According to the above theorem, there exists a unique (up to isomorphism)
dormant $\mr{PGL}_n^{(1)}$-oper $[\msG^\heartsuit]$ on $X_\N$ with ${^{\Diag\!\!}}[\msG^\heartsuit] = [\msF^\heartsuit]$.

\SSP
\bde \label{Def432}
We refer to $[\msG^\heartsuit]$  as the {\bf canonical diagonal lifting} of $[\msF^\heartsuit]$.
\ede

\LSP
\subsection*{Acknowledgements}
The first author was supported by Yoshimura Foundation.
The second author was partially supported by 
 JSPS KAKENHI Grant Number 21K13770.

\end{document}